\documentclass[11pt,a4paper]{article}
\usepackage{amsfonts,amsmath,amssymb,accents,amsthm}
\usepackage{verbatim}

\usepackage{hyperref}
\usepackage[normalem]{ulem}

\usepackage{graphicx}
\usepackage{xcolor}

\newcommand{\dis}{\displaystyle}

\newcommand{\noi}{\noindent}
\newcommand{\halmos}{\rule{1ex}{1.4ex}}
\newcommand{\QED}{\nopagebreak{\hspace*{\fill}$\halmos$\medskip}}
\newcommand{\med}{\medskip}
\newcommand{\quand}{\quad\mbox{and}\quad}

\newtheoremstyle{mythm}
  {}
  {}
  {\itshape}
  {}
  {\bfseries}
  {}
  {.5em}
  {#1 #2 \thmnote{(#3)}}

\theoremstyle{mythm}
\newtheorem{theorem}{Theorem}
\newtheorem{proposition}[theorem]{Proposition}
\newtheorem{lemma}[theorem]{Lemma}

\newtheorem{exercise}[theorem]{Exercise}
\newtheorem{corollary}[theorem]{Corollary}
\newtheorem{conjecture}[theorem]{Conjecture}

\newtheorem{counterex}[theorem]{Counterexample}
\newtheorem{question}{Question}
\newtheorem{problem}[question]{Problem}

\newcommand{\bt}{\begin{theorem}}
\newcommand{\et}{\end{theorem}}
\newcommand{\bl}{\begin{lemma}}
\newcommand{\el}{\end{lemma}}
\newcommand{\bp}{\begin{proposition}}
\newcommand{\ep}{\end{proposition}}
\newcommand{\bcor}{\begin{corollary}}
\newcommand{\ecor}{\end{corollary}}
\newcommand{\br}{\begin{remark}\rm}
\newcommand{\er}{\end{remark}}
\newcommand{\bcon}{\begin{conjecture}}
\newcommand{\econ}{\end{conjecture}}
\newcommand{\bex}{\begin{exercise}}
\newcommand{\eex}{\end{exercise}}
\newcommand{\bcou}{\begin{counterex}}
\newcommand{\ecou}{\end{counterex}}

\theoremstyle{definition}

\newtheorem{remark}[theorem]{Remark}

\newenvironment{Proof}[1][Proof]{\noi\textbf{#1} }{\QED}
\newcommand{\bpro}{\begin{Proof}}
\newcommand{\epro}{\end{Proof}}

\newcommand{\be}{\begin{equation}}
\newcommand{\ee}{\end{equation}}
\newcommand{\ba}{\begin{array}}
\newcommand{\ea}{\end{array}}
\newcommand{\bac}{\begin{array}{r@{\,}c@{\,}l}}
\newcommand{\bc}{\be\begin{array}{r@{\,}c@{\,}l}}
\newcommand{\ec}{\end{array}\ee}

\newcommand{\al}{\alpha}

\newcommand{\ga}{\gamma}

\newcommand{\de}{\delta}
\newcommand{\De}{\Delta}
\newcommand{\eps}{\varepsilon}
\newcommand{\la}{\lambda}
\newcommand{\La}{\Lambda}
\newcommand{\sig}{\sigma}

\newcommand{\om}{\omega}
\newcommand{\Om}{\Omega}

\newcommand{\si}{\ensuremath{\sigma}}

\newcommand{\Ai}{{\cal A}}

\newcommand{\Fi}{{\cal F}}

\newcommand{\Pc}{{\cal P}}

\newcommand{\Ti}{{\cal T}}
\newcommand{\Ui}{{\cal U}}

\newcommand{\R}{{\mathbb R}}
\newcommand{\N}{{\mathbb N}}
\newcommand{\Z}{{\mathbb Z}}

\renewcommand{\S}{{\mathbb S}}
\newcommand{\T}{{\mathbb T}}
\newcommand{\U}{{\mathbb U}}
\newcommand{\E}{{\mathbb E}}
\newcommand{\F}{{\mathbb F}}

\renewcommand{\P}{{\mathbb P}}
\newcommand{\Ob}{{\mathbb O}}
\newcommand{\A}{{\mathbb A}}
\newcommand{\B}{{\mathbb B}}

\newcommand{\Pb}{{\mathbf P}}

\newcommand{\volgt}{\ensuremath{\Rightarrow}}
\newcommand{\up}{\uparrow}

\newcommand{\sub}{\subset}
\newcommand{\beh}{\backslash}

\newcommand{\isd}{\stackrel{\scriptstyle{\rm d}}{=}}

\newcommand{\Asto}[1]{\underset{{#1}\to\infty}{\Longrightarrow}}

\newcommand{\ti}{\tilde}
\newcommand{\ch}{\check}

\newcommand{\ov}{\overline}
\newcommand{\un}{\underline}

\newcommand{\lvec}[1]{\accentset{\leftarrow}{#1}}

\newcommand{\pa}{\partial}
\newcommand{\ffrac}[2]{{\textstyle\frac{{#1}}{{#2}}}}
\newcommand{\dif}[1]{\ffrac{\partial}{\partial{#1}}}

\newcommand{\difif}[2]{\ffrac{\partial^2}{\partial{#1}\partial{#2}}}

\newcommand{\nab}{\nabla}

\newcommand{\di}{\mathrm{d}}
\newcommand{\half}{{[0,\infty)}}

\newcommand{\ha}{\ffrac{1}{2}}

\newcommand{\pbf}{\mathbf{p}}
\newcommand{\rbf}{\mathbf{r}}

\newcommand{\ibf}{\mathbf{i}}
\newcommand{\jbf}{\mathbf{j}}
\newcommand{\kbf}{\mathbf{k}}

\newcommand{\wurz}{\varnothing}
\newcommand{\Hh}{H}

\newcommand{\nnu}{\nu}
\newcommand{\rro}{\rho}
\newcommand{\ngu}{\mu}
\newcommand{\nggu}{\mu'}
\newcommand{\rgo}{\mu}

\newcommand{\percol}[1]{\overset{{#1}}{\longrightarrow}}

\newcommand{\Fx}{\F}

\begin{document}

\makeatletter\@addtoreset{equation}{section}
\makeatother\def\theequation{\thesection.\arabic{equation}}

\renewcommand{\labelenumi}{{\rm (\roman{enumi})}}
\renewcommand{\theenumi}{\roman{enumi}}

\title{Frozen percolation on the binary tree is nonendogenous}
\author{Bal\'azs R\'ath\footnote{MTA-BME Stochastics Research Group,
 Budapest University of Technology and Economics,
 Egry J\'ozsef u.\ 1, 1111 Budapest, Hungary.
 rathb@math.bme.hu},
Jan~M.~Swart\footnote{The Czech Academy of Sciences,
  Institute of Information Theory and Automation,
  Pod vod\'arenskou v\v e\v z\' i 4,
  18200 Praha 8,
  Czech Republic.
  swart@utia.cas.cz},
and Tam\'as Terpai\footnote{E\"otv\"os Lor\'and University,
 Faculty of Science, Department of Analysis,
P\'azm\'any P\'eter s\'et\'any 1/c, 1117 Budapest,  Hungary.
 terpai@math.elte.hu}}

\date{\today}

\maketitle

\begin{abstract}\noi
In frozen percolation, i.i.d.\ uniformly distributed activation times are
assigned to the edges of a graph. At its assigned time, an edge opens provided
neither of its endvertices is part of an infinite open cluster; in the
opposite case, it freezes. Aldous (2000) showed that such a process can be
constructed on the infinite 3-regular tree and asked whether the event that a
given edge freezes is a measurable function of the activation times assigned
to all edges. We give a negative answer to this question, or, using an
equivalent formulation and terminology introduced by Aldous and Bandyopadhyay
(2005), we show that the recursive tree process associated with frozen
percolation on the oriented binary tree is nonendogenous. An essential role in
our proofs is played by a frozen percolation process on a continuous-time
binary Galton Watson tree that has nice scale invariant properties.
\end{abstract}
\vspace{.5cm}

\noi
%
{\it Acknowledgements:} The work of B.~R\'ath is partially supported by
Postdoctoral Fellowship NKFI-PD-121165 and grant NKFI-FK-123962 of NKFI
(National Research, Development and Innovation Office), the Bolyai Research
Scholarship of the Hungarian Academy of Sciences and the \'UNKP-19-4-BME-85
New National Excellence Program of the Ministry for Innovation and Technology.
J.M.~Swart is supported by grant 19-07140S of the Czech Science Foundation (GA
CR). T.~Terpai is partially supported by the National Research, Development
and Innovation Office NKFIH Grant K 120697. We thank James Martin for useful
discussions and M\'arton Sz\H{o}ke for help with the simulations.

\newpage

{\setlength{\parskip}{-2pt}\tableofcontents}

\newpage


\section{Introduction}

\subsection{Frozen percolation on the 3-regular tree}\label{S:unor}

Let $(T,E)$ be a regular tree where each vertex has degree 3, and let
$\Ui=(U_e)_{e\in E}$ be an i.i.d.\ collection of uniformly distributed
$[0,1]$-valued random variables, indexed by the edges of the tree.
We write $E_t:=\{e\in E:U_e\leq t\}$ $(t\in[0,1])$. Aldous \cite{Ald00} has
proved the following theorem.

\bt[Frozen percolation on the 3-regular tree]
It\label{T:frz} is possible to couple $\Ui$ to a random subset $F\sub E$ with
the following properties:
\begin{enumerate}
\item $e\not\in F$ if and only if no endvertex of $e$ is part of an
  infinite cluster of $E_{U_e}\beh(F\cup\{e\})$.
\item The law of $(\Ui,F)$ is invariant under automorphisms of the tree.
\end{enumerate}
\et

At time $t\in[0,1]$, we call edges in $E_t\beh F$ \emph{open}, edges in
$E_t\cap F$ \emph{frozen}, and all other edges \emph{closed}. Then property
(i) can be described in word as follows. Initially all edges are closed. At
its \emph{activation time} $U_e$, the edge $e$ opens provided neither of its
endvertices is at that moment part of an infinite open cluster; in the
opposite case, it freezes.

It is not known if properties (i) and (ii) uniquely determine the joint law of
$(\Ui,F)$. However, it is possible to obtain an object that is unique in law
by adding one natural additional property. To formulate this, we view $T$ as
an oriented graph $(T,\vec E)$ where $\vec E:=\big\{(v,w),(w,v):\{v,w\}\in
E\big\}$ contains two oriented edges for every unoriented edge in $E$.
A \emph{ray} is an infinite sequence of oriented edges
$(v_n,w_n)_{n\geq 0}$ such that $v_n=w_{n-1}$ and $w_n\neq v_{n-1}$ $(n\geq
1)$. We let
\be\ba{r@{\,}l}\label{vecEvw}
\dis\vec E_{(v,w)}:=\big\{(v',w'):&\dis
\exists\mbox{ a ray }(v_n,w_n)_{n\geq 0}\mbox{ and }m\geq 0\\
&\dis\mbox{ s.t.\ }(v_0,w_0)=(v,w)\mbox{ and }(v_m,w_m)=(v',w')\big\}
\ec
denote the union of all rays that start with $(v,w)$, and we let
$E_{(v,w)}:=\big\{\{v',w'\}:(v',w')\in\vec E_{(v,w)}\big\}$ denote the
associated set of unoriented edges. For each subset $S$ of $T$, we let
\be
\pa S:=\big\{(v,w)\in\vec E:v\in S,\ w\in T\beh S\big\}
\ee
denote the collection of oriented edges pointing out of $S$, and we let
$E_S:=\big\{\{v,w\}\in E:v\in S\mbox{ and }w\in S\big\}$ denote the set of
edges induced by $S$. We say that $S$ is a \emph{subtree} if its induced
subgraph $(S,E_S)$ is a tree.

Let $\Ui=(U_e)_{e\in E}$ be as before and let $\vec E_t:=\{(v,w)\in\vec
E:U_{\{v,w\}}\leq t\}$ $(t\in[0,1])$. The existence part of the following
theorem was proved in \cite{Ald00}, but the uniqueness part is~new.

\bt[Frozen percolation on the oriented 3-regular tree]
It\label{T:orfrz} is possible to couple $\Ui$ to a random subset $\vec
F\sub\vec E$ with the following properties:
\begin{enumerate}
\item $(v,w)\in\vec F$ if and only if there exists a ray $(v_n,w_n)_{n\geq 0}$
  with $(v_0,w_0)=(v,w)$ and $(v_n,w_n)\in\vec E_{U_{\{v,w\}}}\beh\vec F$ for all
  $n\geq 1$.
\item The law of $(\Ui,\vec F)$ is invariant under automorphisms of the tree.
\item Let $\Ui_{(v,w)}:=(U_e)_{e\in E_{(v,w)}}$ and $\vec F_{(v,w)}:=\vec
  F\cap\vec E_{(v,w)}$. Then, for each finite subtree $S\sub T$, the random
  variables $(\Ui_{(v,w)},\vec F_{(v,w)})_{(v,w)\in\pa S}$ are independent of
  each other and of $(\Ui_e)_{e\in E_S}$.
\end{enumerate}
These properies uniquely determine the joint law of $(\Ui,\vec F)$.
Moreover, setting $F:=\big\{\{v,w\}\in E:(v,w)\in\vec F\mbox{ or }(w,v)\in\vec
F\big\}$ defines a pair $(\Ui,F)$ with properties (i) and (ii) of
Theorem~\ref{T:frz}.
\et

In this paper, our main interest is not in uniqueness in law but rather in
almost sure uniqueness. In \cite[Section~5.7]{Ald00}, Aldous asked whether the
set $F$ of frozen edges is measurable w.r.t.\ the \si-field generated by
$\Ui$, and cautiously conjectured that this might indeed be the case. In
\cite[Thm~55]{AB05}, an apparent proof of this conjecture by Bandyopadhyay was
announced that appeared on the arXiv \cite{Ban04} but turned out to contain an
error. In the last posted update of \cite{Ban04} from 2006, Bandyopadhyay
reported on numerical simulations (similar to those shown in
Figure~\ref{fig:rhoiter} below) that suggested nonuniqueness, and from this
moment on this seems to have been the generally held belief. We finally settle
the issue by proving this.


\bt[Frozen percolation is not almost sure unique]
Let\label{T:main} $(\Ui,F)$ be the pair defined in Theorem~\ref{T:orfrz} and
let $F'$ be a copy of $F$, conditionally independent of $F$ given $\Ui$.
Then $F\neq F'$ a.s. In particular, the random variable $F$ is not measurable
w.r.t.\ the \si-field generated by $\Ui$.
\et

The proofs of Theorems~\ref{T:orfrz} and \ref{T:main} will be given in
Subsection~\ref{S:unorient}.

\subsection{Frozen percolation on the oriented binary tree}\label{S:binar}

For a given oriented edge $(v,w)\in\vec E$ of the 3-regular tree, the set
$\vec E_{(v,w)}$ of oriented edges that lie on rays starting with $(v,w)$ can
naturally be labeled with the space $\T$ of all finite words $\ibf=i_1\cdots
i_n$ $(n\geq 0)$ made up from the alphabet $\{1,2\}$. We call $|\ibf|:=n$ the
length of the word $\ibf$ and denote the word of length zero by $\wurz$, which
we distinguish notationally from the empty set $\emptyset$. The concatenation
of two words $\ibf=i_1\cdots i_n$ and $\jbf=j_1\cdots j_m$ is denoted by
$\ibf\jbf:=i_1\cdots i_nj_1\cdots j_m$.

Apart from using $\T$ to label oriented edges as above, we can also interpret
$\T$ as labeling the vertices of a binary tree with root $\wurz$, in which
each vertex $\ibf$ has two descendants $\ibf 1,\ibf 2$ and each vertex
$\ibf=i_1\cdots i_n$ $(n\geq 1)$ except the root has a unique predecessor
$\lvec\ibf:=i_1\cdots i_{n-1}$. By definition, a \emph{ray} starting at $\ibf$
is a sequence $(\ibf_n)_{n\geq 0}$ such that $\ibf_0=\ibf$ and
$\lvec\ibf_n=\ibf_{n-1}$ $(n\geq 1)$. For any $A\sub\T$ and $\ibf\in\T$, we
write $\ibf\percol{A}\infty$ if there exists a ray $(\ibf_n)_{n\geq 0}$ with
$\ibf_0=\ibf$ and $\ibf_n\in A$ $(n\geq 0)$.

We write $\ibf\prec\jbf$ if $\jbf=\ibf\kbf$ for some $\kbf\in\T$. By
definition, a \emph{rooted subtree} of $\T$ is a set $\U\sub\T$ with the
property that $\ibf\prec\jbf\in\U$ implies $\ibf\in\U$. For each nonempty
rooted subtree $\U$ of $\T$, we let
$\pa\U:=\{\ibf\in\T\beh\U:\lvec\ibf\in\U\}$ denote the \emph{boundary} of $\U$
relative to $\T$, and we use the convention that $\pa\U=\{\wurz\}$ if
$\U=\emptyset$.

Let $\tau=(\tau_\ibf)_{\ibf\in\T}$ be an i.i.d.\ collection of uniformly
distributed $[0,1]$-valued random variables. In the picture where elements of
$\T$ label oriented edges in $\vec E_{(v,w)}$, this corresponds to the
collection of activation times $(U_e)_{e\in E_{(v,w)}}$.  Using the same
picture, let $(X_\ibf)_{\ibf\in\T}$ be a collection of real random variables,
which correspond to the first time when there is an infinite open ray of edges
starting with a given oriented edge, with $X_\ibf:=\infty$ if this never
happens. Note that $X_\ibf$ takes values in $I:=[0,1]\cup\{\infty\}$. By
properties~(ii) and (iii) of Theorem~\ref{T:orfrz}, for each finite rooted
subtree $\U\sub\T$,
\be\label{XU}
\mbox{the r.v.'s $(X_\ibf)_{\ibf\in\pa\U}$ are i.i.d.\ and independent of
$(\tau_\ibf)_{\ibf\in\U}$.}
\ee
Using also property~(i), it is easy to see that the random variables
$(X_\ibf)_{\ibf\in\T}$ satisfy the inductive relation (compare
\cite[formula (65)]{AB05})
\be\label{Xind}
X_\ibf=\ga[\tau_\ibf](X_{\ibf 1},X_{\ibf 2})\qquad(\ibf\in\T),
\ee
where $\ga:[0,1]\times I^2\to I$ is defined as
\be\label{gadef}
\ga[t](x,y):=\left\{\ba{ll}
x\wedge y\quad&\mbox{if }x\wedge y>t,\\[5pt]
\infty\quad&\mbox{otherwise.}\ea\right.
\ee
Generalising from the set-up of Theorem~\ref{T:orfrz}, we will more generally
be interested in collections of random variables
$(\tau_\ibf,X_\ibf)_{\ibf\in\T}$ such that $(\tau_\ibf)_{\ibf\in\T}$ are
i.i.d.\ uniformly distributed on $[0,1]$ and (\ref{XU}) and (\ref{Xind})
hold. As will be explained in the next subsection, in the terminology of
\cite{AB05}, such a collection forms a \emph{Recursive Tree Process}
(RTP). The theory of RTPs provides us with a convenient general framework to
reformulate and prove Theorems~\ref{T:orfrz} and \ref{T:main}.

\subsection{Recursive Tree Processes}\label{S:RTP}

Roughly speaking, a \emph{Recursive Tree Process} (RTP) is a stationary Markov
chain in which time has a tree-like structure and flows in the direction of
the root. The state at each node of the tree is a function of the states of
its descendants and i.i.d.\ randomness attached to the nodes. Following
\cite{AB05}, we call an RTP \emph{endogenous} if the state at the root is
measurable w.r.t.\ the \si-field generated by the i.i.d.\ randomness attached
to the nodes. It has been shown in \cite[Thm~11]{AB05} that endogeny is
equivalent to \emph{bivariate uniqueness}. We first explain these concepts in
a general setting and then specialise to frozen percolation.

Slightly generalising our previous notation, let $\T$ denote the space of all
finite words $\ibf=i_1\cdots i_n$ $(n\geq 0)$ made up from the alphabet
$\{1,\ldots,d\}$, where $d\geq 1$ is some fixed integer.
All previous notation involving the binary tree generalizes in a
straightforward manner to the $d$-ary tree $\T$.
Let $I$ and $\Om$ be Polish spaces, let $\ga:\Om\times I^d\to I$ be a
measurable function, and let $(\om_\ibf)_{\ibf\in\T}$ be i.i.d.\ $\Om$-valued
random variables with common law $\pbf$. Let $\nu$ be a probability law on $I$
that solves the \emph{Recursive Distributional Equation} (RDE)
\be\label{RDE}
X_\wurz\isd\ga[\om_\wurz](X_1,\ldots,X_d),
\ee
where $\isd$ denotes equality in distribution, $X_\wurz$ has law $\nu$, and
$X_1,\dots, X_d$ are copies of $X_\wurz$, independent of each other and of
$\om_\wurz$. A simple argument based on Kolmogorov's extension theorem (see
\cite[Lemma~8]{MSS20}) tells us that the i.i.d.\ random variables
$(\om_\ibf)_{\ibf\in\T}$ can be coupled to $I$-valued random
variables $(X_\ibf)_{\ibf\in\T}$ in such a way that:
\begin{enumerate}
\item For each finite rooted subtree $\U\sub\T$, the r.v.'s
  $(X_\ibf)_{\ibf\in\pa\U}$ are i.i.d.\ with common law $\nu$ and independent
  of $(\om_\ibf)_{\ibf\in\U}$.
\item $\dis X_\ibf=\ga[\om_\ibf](X_{\ibf 1},\ldots,X_{\ibf d})\qquad(\ibf\in\T)$.
\end{enumerate}
Moreover, these conditions uniquely determine the joint law of
$(\om_\ibf,X_\ibf)_{\ibf\in\T}$. We call the latter the
\emph{Recursive Tree Process} (RTP) corresponding to the maps $\ga$ and
solution $\nu$ of the RDE (\ref{RDE}). By definition, the
RTP $(\om_\ibf,X_\ibf)_{\ibf\in\T}$ is \emph{endogenous} if the
random variable $X_\wurz$ is measurable w.r.t.\ the \si-field generated by the
random variables $(\om_\ibf)_{\ibf\in\T}$. It has been shown in
\cite[Thm~11]{AB05} that this is equivalent to \emph{bivariate uniqueness}, as
we now explain.

Let $\Pc(I)$ denote the space of all probability measures on $I$. We can
define a map $T:\Pc(I)\to\Pc(I)$ by
\be\label{Tdef}
T(\mu):=\mbox{ the law of }\ga[\om_\wurz](X_1,\ldots,X_d),
\ee
where $X_1,\ldots,X_d$ are i.i.d.\ with law $\mu$ and independent of
$\om_\wurz$. In particular, solutions to the RDE (\ref{RDE}) correspond to
fixed points of $T$. Similarly, we can define a \emph{bivariate map}
$T^{(2)}:\Pc(I^2)\to\Pc(I^2)$ by
\be\label{bivdef}
T^{(2)}(\mu^{(2)}):=\mbox{ the law of }
\big(\ga[\om_\wurz](X_1,\ldots,X_d),\ga[\om_\wurz](X^*_1,\ldots,X^*_d)\big),
\ee
where $(X_1,X^*_1),\ldots,(X_d,X^*_d)$ are i.i.d.\ with common law
$\mu^{(2)}$ and independent of $\om_\wurz$. A trivial way to construct a
fixed point of $T^{(2)}$ is to set
\be\label{ovnu}
\ov\nu^{(2)}:=\P\big[(X_\wurz,X_\wurz)\in\,\cdot\,\big]
\ee
where the law $\nu$ of $X_\wurz$ is a fixed point of $T$. We will refer to
$\ov\nu^{(2)}$ as the \emph{trivial} fixed point or as the \emph{diagonal}
fixed point of $T^{(2)}$ with marginal distribution $\nu$. A more
interesting way to construct a fixed point of $T^{(2)}$ goes as follows. Let
$(\om_\ibf,X_\ibf)_{\ibf\in\T}$ be the RTP corresponding to the
map $\ga$ and a fixed point $\nu$ of $T$, and
let $(X'_\ibf)_{\ibf\in\T}$ be a copy of $(X_\ibf)_{\ibf\in\T}$, conditionally
independent given $(\om_\ibf)_{\ibf\in\T}$. It follows from
\cite[Lemma~2 and Prop~4]{MSS18} that
\be\label{unnu}
\un\nu^{(2)}:=\P\big[(X_\wurz,X'_\wurz)\in\,\cdot\,\big]
\ee
is also a fixed point of $T^{(2)}$. Let us denote by
 $(T^{(2)})^n$  the
$n$-th iterate of the bivariate map $T^{(2)}$.
 By \cite[Lemma~2 and Prop.~3]{MSS18},
one has
\be\label{prodtoun}
(T^{(2)})^n(\nu\otimes\nu)\Asto{n}\un\nu^{(2)}.
\ee

The following theorem links endogeny to bivariate uniqueness. The essential
idea goes back to \cite[Thm~11]{AB05}. In its present form, it follows from
\cite[Thms~1 and 5 and Lemma~14]{MSS18}. Below, $\Pc(I^2)_\nu$ denotes the
space of all probability measures on $I^2$ whose one-dimensional marginals are
given by $\nu$. Note that condition~(ii) below and formula (\ref{prodtoun})
suggest a method to numerically investigate whether an RTP is endogenous,
compare Figure~\ref{fig:rhoiter} below.

\bt[Endogeny and bivariate uniqueness]
Let\label{T:bivar} $(\om_\ibf,X_\ibf)_{\ibf\in\T}$ be the RTP
corresponding to a map $\ga$ and a solution $\nu$ of the corresponding
RDE (\ref{RDE}). Then the following statements are equivalent:
\begin{enumerate}
\item The RTP $(\om_\ibf,X_\ibf)_{\ibf\in\T}$ is endogenous.
\item $\un\nu^{(2)}=\ov\nu^{(2)}$.
\item The bivariate map $T^{(2)}$ has a unique fixed point in $\Pc(I^2)_\nu$.
\item $\dis(T^{(2)})^n(\mu^{(2)})\Asto{n}\ov\nu^{(2)}$ for all
  $\mu^{(2)}\in\Pc(I^2)_\nu$.
\end{enumerate}
\et

Note that since we know that $\un\nu^{(2)}$ and $\ov\nu^{(2)}$ are fixed points,
the implications (iv)$\volgt$(iii)$\volgt$(ii) are trivial. The
implication (ii)$\volgt$(i) follows from our characterisation of
$\un\nu^{(2)}$ in (\ref{unnu}), so the essential claim is that (i) implies (iv).

\subsection{Nonendogeny}\label{S:nonend}

Specialising from the general set-up of the previous subsection, we set
$d:=2$ and as our i.i.d.\ randomness $(\om_\ibf)_{\ibf\in\T}$ we use an
i.i.d.\ collection $(\tau_\ibf)_{\ibf\in\T}$ of uniformly distributed
$[0,1]$-valued random variables. We set $I:=[0,1]\cup\{\infty\}$, and choose for
$\ga:\Om\times I^2\to I$ the map defined in (\ref{gadef}).
Using these objects, we define a map $T:\Pc(I)\to\Pc(I)$ as in
(\ref{Tdef}). The associated RDE $T(\ngu)=\ngu$ then takes the form (compare
(\ref{RDE}))
\be\label{frzRDE}
X_\wurz\isd\ga[\tau_\wurz](X_1,X_2),
\ee
where $\isd$ denotes equality in distribution, $X_\wurz$ has law $\ngu$, and
$X_1,X_2$ are copies of $X_\wurz$, independent of each other and of
$\tau_\wurz$. Solutions to the RDE (\ref{frzRDE}) are not unique. We will
describe all solutions of (\ref{frzRDE}) in Lemma~\ref{L:equivRDE} and
Proposition~\ref{P:genRDE} below.

Let $(\tau_\ibf,X_\ibf)_{\ibf\in\T}$ be an RTP corresponding to the map $\ga$
in (\ref{gadef}) and an arbitrary solution to the RDE (\ref{frzRDE}). We
set
\be\label{Fxdef}
\T^t:=\big\{\ibf\in\T:\tau_\ibf\leq t\big\}
\quand
\Fx:=\big\{\ibf\in\T:\tau_\ibf\geq X_{\ibf 1}\wedge X_{\ibf 2}\big\},
\ee
and define $I$-valued random variables $(X^\up_\ibf)_{\ibf\in\T}$ by
\be\label{Xupdef}
X^\up_\ibf:=\inf\big\{t\in[0,1]:\ibf\percol{\T^t\beh\Fx}\infty\big\},
\ee
with $X^\up_\ibf:=\infty$ if the set on the right-hand side is empty. In line
with our interpretation where elements of $\T$ represent oriented edges in
$\vec E_{(v,w)}$ (with $(v,w)$ fixed), we say that at time $t\in[0,1]$, points
in $\T^t\beh\Fx$ are \emph{open}, points in $\T^t\cap\Fx$ are
\emph{frozen}, and all other points in $\T$ are \emph{closed}. We call
$\tau_\ibf$ the \emph{activation time} of $\ibf$ and refer to $X_\ibf$ and
$X^\up_\ibf$ as its \emph{burning time} and \emph{percolation time},
respectively. Note that since the subset of $[0,1]$ on the right-hand side of (\ref{Xupdef})
is a.s.\ closed (in the topological sense), $\ibf$ percolates at time $t$ if and only if $X^\up_\ibf\leq
t$. Formula (\ref{Fxdef}) says that initially, all points $\ibf\in\T$ are
closed. At its activation time $\tau_\ibf$, the point $\ibf$ freezes if at
that moment one of its descendants is burnt, and opens otherwise.

It follows from the inductive relation (\ref{Xind}) that $X_\ibf>\tau_\ibf$
a.s., i.e., a point $\ibf$ can only burn after its activation time. Comparing
the definition of $\Fx$ in (\ref{Fxdef}) with the definition of the map
$\ga$ in (\ref{gadef}), we observe that if $\ibf$ burns at some time
$X_\ibf\in[0,1]$, then $\ibf$ must be open at that time. Moreover, by
(\ref{gadef}), if $\ibf$ burns at some time $X_\ibf\in[0,1]$, then there must
be a ray starting at $\ibf$ of points that burn at the same time as $\ibf$. In
particular, such a ray must be open, which proves that
\be\label{XupX}
X^\up_\ibf\leq X_\ibf\quad{\rm a.s.}\quad(\ibf\in\T).
\ee
We will prove Theorem~\ref{T:orfrz} by showing that the opposite inequality
holds a.s.\ if and only if $(\tau_\ibf,X_\ibf)_{\ibf\in\T}$ is the RTP
corresponding to one particular solution of the RDE (\ref{frzRDE}). This
solution is described by the following lemma, which we cite from
\cite[Lemma~3]{Ald00}. Note that (\ref{nudef}) below implies that
$\nnu(\{\infty\})=\ha$.

\bl[Special solution to the RDE]
Let\label{L:frzRDE} $\nnu$ denote the probability measure on $I$ defined by
\be\label{nudef}
\nnu\big((t,1]\cup\{\infty\}\big):=1\wedge\frac{1}{2t}\qquad\big(t\in[0,1]\big).
\ee
Then $\nnu$ solves the RDE (\ref{frzRDE}).
\el

We will deduce Theorem~\ref{T:orfrz} from the following, more precise theorem.
Aldous proved the ``if'' part of the statement below in \cite{Ald00}, but the
``only if'' part is new. Theorem~\ref{T:dirfrz} is proved in
Subsection~\ref{S:binfrz}.

\bt[Frozen percolation on the oriented binary tree]
Consider\label{T:dirfrz} an RTP $(\tau_\ibf,X_\ibf)_{\ibf\in\T}$ corresponding
to the map $\ga$ in (\ref{gadef}) and an arbitrary solution $\mu$ to the RDE
(\ref{frzRDE}). Let $(X^\up_\ibf)_{\ibf\in\T}$ be defined as in
(\ref{Xupdef}). Then one has $X^\up_\wurz=X_\wurz$ a.s.\ if and only if
$\mu=\nnu$, the measure defined in (\ref{nudef}).
\et

Using the language of RTPs, we can formulate our main result as follows.
Theorem~\ref{T:main} will follow from the theorem below in a straightforward
manner using methods from \cite{Ald00}. Theorem~\ref{T:nonend} is proved in
Subsection~\ref{S:RDEequiv}.

\bt[Frozen percolation on the binary tree is nonendogenous]
The\label{T:nonend} RTP $(\tau_\ibf,X_\ibf)_{\ibf\in\T}$ corresponding to the
map $\ga$ from (\ref{gadef}) and the law $\nnu$ from (\ref{nudef}) is
nonendogenous.
\et

We prove Theorem~\ref{T:nonend} using Theorem~\ref{T:bivar}, by explicitly
identifying the solution $\un\nnu^{(2)}$ to the bivariate RDE in terms of the
solution to a certain differential equation (see formula (\ref{explit}) below)
and showing that $\un\nnu^{(2)}$ is not equal to $\ov\nnu^{(2)}$. An essential
role in our proofs is played by frozen percolation on a continuum tree that we
will call the \emph{Marked Binary Branching Tree} (MBBT). The advantage of
working with the latter is that it enjoys a nice scaling property that will
significantly simplify our analysis. A simple trick then allows us to relate
frozen percolation on the oriented binary tree to frozen percolation on the
MBBT and also prove Theorem~\ref{T:nonend}.

\subsection{The Marked Binary Branching Tree}

Roughly speaking, the \emph{marked binary branching tree} is the family tree
of a continuous-time, rate one binary branching process, equipped with a marked
Poisson point process of intensity one and uniformly distributed
$[0,1]$-valued marks. We now introduce this object more formally.

Let $(A_h)_{h\geq 0}$ be a continuous-time branching process, started with a
single particle, where each particle splits into two new particles with rate
one. We view $A_h$ as an evolving set. In particular, the cardinality $|A_h|$
is a Markov process in $\N$ that jumps from $a$ to $a+1$ with rate $a$, and
$A_0=\{x_0\}$ is a set containing a single element $x_0$. In the next
subsection, we will make a more explicit choice for the labels of elements of
$A_h$. We choose $(A_h)_{h\geq 0}$ to be right-continuous and let
$(A_{h-})_{h\geq 0}$ denote its left-continuous modification.

For each pair of times $g,h\geq 0$ and individuals $x\in A_g$, $y\in A_h$ that
are alive at these times, let $\tau(x,y)$ denote the last time in
$[0,g\wedge h]$ that a common ancestor existed of $x$ and $y$, and let
\be\label{gendist}
d\big((x,g),(y,h)\big):=g+h-2\tau(x,y)
\ee
denote their genetic distance. Then the random set
\be\label{Tidef}
\Ti:=\big\{(x,h):x\in A_{h-},\ h\geq 0\big\}
\ee
equipped with the metric (\ref{gendist}) is a random continuum tree.
In pictures, we draw $x$ horizontally and $h$ vertically, and from now on, we
refer to $h$ as the \emph{height}, rather than time, of a point $(x,h)=z\in\Ti$.
We call $\wurz:=(x_0,0)$ the \emph{root} of $\Ti$.

\begin{figure}[htb!]
\begin{center}
\includegraphics{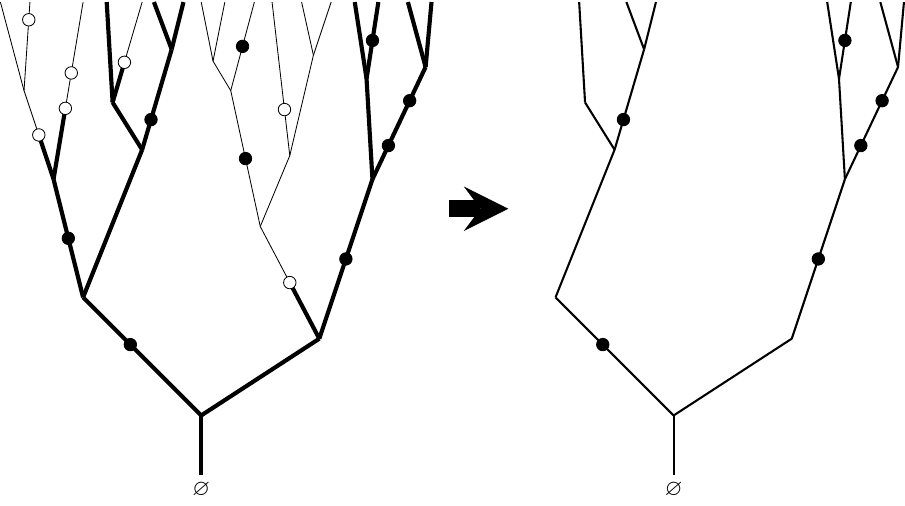}
\caption{Scaling of the marked binary branching tree. On the left:
at time $t$, points in $\Pi_t$ are still closed and marked with white circles,
while points in $\Pi_0\beh\Pi_t$ have already opened and are marked with
black circles. On the right: removing the loose ends from the open cluster at
the root yields a stretched version of the original marked binary branching
tree.}
\label{fig:MBBT}
\end{center}
\end{figure}

Conditional on $\Ti$, we let $\Pi_0$ be a Poisson point set on $\Ti$
whose intensity measure is the length measure on $\Ti$, and conditional on
$\Ti$ and $\Pi_0$, we let $(\tau_z)_{z\in\Pi_0}$ be i.i.d.\ uniformly distributed
$[0,1]$-valued marks. We think of $z=(x,h)\in\Pi_0$ as a hole on $\Ti$ that
disappears (i.e., gets filled in) at time $\tau_z$. We observe that
\be
\Pi=\big\{(z,\tau_z):z\in\Pi_0\big\}
\ee
is a Poisson set of intensity one on $\Ti\times[0,1]$ and that $\Pi_0$ as well
as the marks $(\tau_z)_{z\in\Pi_0}$ can be read off from $\Pi$. For lack of
better name, we call the pair $(\Ti,\Pi)$ the \emph{Marked Binary Branching
  Tree} (MBBT). See Figure~\ref{fig:MBBT} for an illustration.

We set
\be
\Pi_t:=\big\{z\in\Pi_0:\tau_z>t\big\}\quad(t\in[0,1]).
\ee
Intuitively, $\Pi_t$ is the set of holes on $\Ti$ which are still present at
time $t$. For any set $A\sub\Ti$ and point $z\in\Ti$, we write
$\wurz\percol{\Ti\beh A}z$ if $\wurz$ and $z$ are connected in $\Ti\beh A$ and
we write $z\percol{\Ti\beh A}\infty$ if there exists an infinite, continuous,
upward path through $\Ti\beh A$. We start with a simple observation. Note
that below, in contrast with our earlier notation $E_t$, points in $\Pi_t$
play the role of points that can \emph{not} be passed at time $t$.
The proof of the following lemma can be found in Subsection~\ref{S:frzMBBT}.

\bl[Oriented percolation on the marked binary branching tree]
One\label{L:perc} has
\begin{equation} \dis\P\big[\wurz\percol{\Ti\beh\Pi_t}\infty\big]
=t\quad(0\leq t\leq 1).
\end{equation}
\el

Indeed, if we cut $\Ti$ at points of $\Pi_t$, then the remaining connected
component of the root is the family tree of a branching process
where particles split into two with rate one and die with rate $1-t$. It is an
elementary exercise in branching theory to show that the survival probability
of such a branching process is $t$. The fact that the survival probability is
a linear function of $t$ reflects a scaling property of the marked binary
branching tree that will be important in our analysis. Below, we view
$(\Ti,\Pi)$ as a marked metric space, i.e., we consider two marked trees to be
equal if one can be mapped onto the other by an isometry that preserves the
marks. The following result is proved in Subsection~\ref{S:frzMBBT}.

\bp[Scaling]
Let\label{P:scale} $(\Ti,\Pi)$ be the marked binary branching tree. Fix
$0<t<1$ and define
\be\label{Tiac}
\Ti':=\big\{z\in\Ti:
\wurz\percol{\Ti\beh\Pi_t}z\percol{\Ti\beh\Pi_t}\infty\big\},\quad
\Pi':=\big\{(z,\tau_z)\in\Pi:z\in\Ti'\big\}.
\ee
Then the probability that $\Ti' \neq \emptyset$ is $t$ and conditional on this
event, the pair $(\Ti',\Pi')$, viewed as a marked metric space, is equally
distributed with the stretched marked binary branching tree $(\Ti'',\Pi'')$
defined as
\be
\Ti'':=\big\{(x,t^{-1}h):(x,h)\in\Ti\big\},\ 
\Pi'':=\big\{(x,t^{-1}h,t\tau_{(x,h)}):(x,h,\tau_{(x,h)})\in\Pi\big\}.
\ee
\ep

In words, this says that if we cut off all parts of $\Ti$ that lie above
points $z\in\Pi_t$, then remove the loose ends of the tree, and condition on
the event that the remaining tree is nonempty, then we end up with the family
tree of a branching process where particles split into two with rate $t$,
equipped with a marked Poisson point set with intensity $t$ and
i.i.d.\ uniformly distributed on $[0,t]$-valued marks. See
Figure~\ref{fig:MBBT} for an illustration.


\subsection{Frozen percolation on the MBBT}\label{S:MBBT}

In the previous subsection, we have been deliberately vague about the labeling
of elements of the evolving set-valued branching process $(A_h)_{h\geq 0}$. We
now make an explicit choice, which naturally leads to an RTP that is
closely related to, but different from the one introduced in
Subsection~\ref{S:nonend}.

We will construct the branching process $(A_h)_{h\geq 0}$ in such a way that
$A_0=\{\wurz\}$ and $A_h\sub\T$ for all $h\geq 0$. (Note that by a slight
abuse of notation, $\wurz$ now denotes both the root of the discrete tree $\T$
and of the continuum tree $\Ti$, the latter being defined as
$\wurz=(\wurz,0)$.) Each element $\ibf\in A_h$ branches with rate one into two
new elements labeled $\ibf 1$ and $\ibf 2$. In addition, we arrange things in
such a way that each element $\ibf\in A_h$ is with rate one replaced by a new
element labeled $\ibf 1$. The idea of this is to encode the Poisson point set
$\Pi_0$ from the MBBT in terms of the labels of elements of $A_h$, in such a
way that $\Pi_0$ is given by the collection of points $(\ibf,h)$ for which
$\ibf$ is at time $h$ replaced by $\ibf 1$.

We will give an explicit construction of the MBBT based on three collections
of i.i.d.\ random variables:
\begin{enumerate}
\item $(\tau_\ibf)_{\ibf\in\T}$ are i.i.d.\ uniformly distributed on $[0,1]$,
\item $(\kappa_\ibf)_{\ibf\in\T}$ are i.i.d.\ uniformly distributed on
  $\{1,2\}$,
\item $(\ell_\ibf)_{\ibf\in\T}$ are i.i.d.\ exponentially distributed with
  mean $1/2$.
\end{enumerate}
We interpret $\ell_\ibf$ as the lifetime of the individual $\ibf$ and let
\be\label{birthdeath}
b_{i_1\cdots i_n}:=\sum_{k=0}^{n-1}\ell_{i_1\cdots i_k}
\quand
d_{i_1\cdots i_n}:=\sum_{k=0}^{n}\ell_{i_1\cdots i_k}
\ee
with $b_\wurz:=0$ and $d_\wurz:=\ell_\wurz$ denote the birth and death
times of $\ibf$. The random variable $\kappa_i$ indicates what happens with the
individual $\ibf$ at the end of its lifetime. If $\kappa_\ibf=1$, then it is
replaced by a single new individual with label $\ibf 1$, and if
$\kappa_\ibf=2$, then it is replaced by two new individuals with labels $\ibf 1$
and $\ibf 2$. In line with this, we let $\S$ denote the random subtree of $\T$
defined by
\be\label{Sdef}
\S:=\big\{ \, i_1\cdots i_n\in\T:
i_m\leq\kappa_{i_1\cdots i_{m-1}}\ \forall\, 1\leq m\leq n \, \big\},
\ee
which is the collection of all individuals that will ever be born.
Recall that $\pa\U$ denotes the boundary of a rooted subtree $\U\sub\T$
relative to $\T$. Likewise, for any rooted subtree $\U\sub\S$ we let
$\nab\U:=\pa\U\cap\S$ denote the \emph{boundary} of $\U$ relative to $\S$.

For $h\geq 0$, we let
\be\ba{l@{\qquad}l}\label{nabS}
\dis\T_h:=\big\{\ibf\in\T: d_\ibf\leq h \big\},
&\dis\pa\T_h=\big\{\ibf\in\T: b_\ibf\leq h<d_\ibf \big\},\\[5pt]
\dis\S_h:=\T_h\cap\S,
&\dis\nab\S_h=\pa\T_h\cap\S
\ec
denote the sets of individuals that have died by time $h$ and those that are
alive at time $h$, respectively. Note that the former are a.s.\ finite rooted
subtrees of $\T$ and $\S$, respectively, and the latter are their boundaries.
Then
\be
(\nab\S_h)_{h\geq 0}=(A_h)_{h\geq 0}
\ee
gives an explicit construction of the branching process $(A_h)_{h\geq 0}$ we
have earlier described in words. Defining $\Ti$ as in (\ref{Tidef}) and
setting
\be\label{Pidef}
\Pi:=\big\{(\ibf,d_\ibf,\tau_\ibf):\ibf\in\S,\ \kappa_\ibf=1\big\}
\ee
yields an explicit construction of the MBBT $(\Ti,\Pi)$ based on
i.i.d.\ randomness.

Instead of giving a description of oriented frozen percolation on $(\Ti,\Pi)$
similar to Theorem~\ref{T:dirfrz}, we immediately jump to the corresponding
RTP for the percolation times. Letting $Y_\ibf$ denote the first time when
there is an infinite upwards open path in frozen percolation on $(\Ti,\Pi)$
starting from the point $(\ibf,b_\ibf)$, it is not hard to see that
$(Y_\ibf)_{\ibf\in\S}$ must satisfy the inductive relation
\be\label{Yind}
Y_\ibf=\chi[\tau_\ibf,\kappa_\ibf](Y_{\ibf 1},Y_{\ibf 2}),
\ee
where $\chi:[0,1]\times\{1,2\}\times I^2\to I$ is the function
\be\label{chidef}
\chi[\tau,\kappa](x,y):=\left\{\ba{ll}
x\quad&\mbox{if }\kappa=1,\ x>\tau,\\[5pt]
\infty\quad&\mbox{if }\kappa=1,\ x\leq\tau,\\[5pt]
x\wedge y\quad&\mbox{if }\kappa=2.\ea\right.
\ee
Note that in (\ref{Yind}), $Y_\ibf$ is a priori only defined for $\ibf\in\S$.
The definition of $\S$ in (\ref{Sdef}) is such, however, that in cases when
$\ibf\in\S$ but $\ibf 2\not\in\S$, the value of $Y_{\ibf 2}$ is irrelevant for
the outcome of the function $\chi$. The subtree $\S$ plays an important role
in the theory of continuous-time RTPs, see \cite[Sect.~1.4]{MSS20}.

Like in the case of the oriented binary tree (as discussed in
Subsection~\ref{S:nonend}) it is possible to go the other way, i.e., starting
from a solution to the RDE corresponding to the map $\chi$, one can construct
an RTP $(\tau_\ibf,\kappa_\ibf,Y_\ibf)_{\ibf\in\T}$ where now $Y_\ibf$ is
defined for all $\ibf\in\T$, and then restrict to $\S$ to construct oriented
frozen percolation on the MBBT. In the present setting, it turns out that the
``right'' solution to the corresponding RDE is given by the following
lemma. Since we will later (in Subsection~\ref{S:binfrz} below) see that frozen
percolation on the MBBT and on the oriented binary tree can be mapped into
each other, we will at this moment not explain why in the present setting,
Lemma~\ref{L:MRDE} describes the ``right'' solution to the RDE.

\bl[Special solution to the RDE]
Let\label{L:MRDE} $\rro$ denote the probability measure on $I$ defined by
\be\label{rhodef}
\rro\big([0,t]\big):=\ha t\qquad\big(t\in[0,1]\big), \qquad \rro(\{\infty\}):=\ha.
\ee
Then $\rro$ solves the RDE
\be\label{MRDE}
Y_\wurz\isd\chi[\tau_\wurz,\kappa_\wurz](Y_1,Y_2),
\ee
where $\isd$ denotes equality in distribution, $Y_\wurz$ has law $\rro$, and
$Y_1,Y_2$ are copies of $Y_\wurz$, independent of each other and of
$\tau_\wurz,\kappa_\wurz$.
\el

\bpro
Let $Y_1,Y_2$ be i.i.d.\ with law $\rro$, let $\tau$ and $\kappa$ be
independent r.v.'s that are uniformly distributed on $[0,1]$ and $\{1,2\}$,
respectively, and define $Y_\wurz:=\chi[\tau,\kappa](Y_1,Y_2)$, where $\chi$
is defined in \eqref{chidef}. We claim that $Y_\wurz$ has law $\rro$. Indeed,
for each $t\in[0,1]$, we have
\bc
\dis\P[Y_\wurz \leq t]
&=&\dis\ha\int_0^1\P[\chi[s,1](Y_1,Y_2) \leq t]\, \di s
+\ha\P[Y_1\wedge Y_2 \leq t]\\[5pt]
&=&\dis\ha\int_0^1 \P[s\leq Y_1 \leq t] \,  \di s
+\ha\big(1-\P[Y_1 > t]^2\big)\\[5pt]
&=&\dis\ha\int_0^t (\ha t-\ha s)\,  \di s
+\ha\big(1-(1-\ha t)^2\big)=\ha t.
\ec
\epro

We will prove that the RTP $(\tau_\ibf,\kappa_\ibf,Y_\ibf)_{\ibf\in\T}$
corresponding to the map $\chi$ from (\ref{chidef}) and law $\rro$ from
(\ref{rhodef}) is nonendogenous. We apply Theorem~\ref{T:bivar}. We will
explicitly identify the special solutions $\un\rro^{(2)}$ and $\ov\rro^{(2)}$
to the bivariate RDE and show that they are not equal.

It is clear from the definitions of $\un\rro^{(2)}$ and $\ov\rro^{(2)}$
in (\ref{ovnu}) and (\ref{unnu}) that both measures are symmetric under a
permutation of the two coordinates and that their one-dimensional marginals
equal $\rro$. The main advantage of working with the MBBT is that as a result
of the scaling property described in Proposition~\ref{P:scale}, the measures
$\un\rro^{(2)}$ and $\ov\rro^{(2)}$ are also scale invariant. We let
$\Pc_\ast(I^2)_\rro$ denote the space of symmetric measures $\mu^{(2)}$ on
$I^2$ whose one-dimensional marginals are given by $\rro$ and
that are moreover scale invariant in the sense that
\be\label{scinv}
\mu^{(2)}\big([0,tr]\times[0,ts]\big)
=t\mu^{(2)}\big([0,r]\times[0,s]\big)\qquad\big(r,s,t\in[0,1]\big).
\ee
The following lemma is proved in Subsection~\ref{S:frzMBBT}.

\bl[Scale invariance]
One\label{L:scinv} has $\un\rro^{(2)},\ov\rro^{(2)}\in\Pc_\ast(I^2)_\rro$.
\el

By Theorem~\ref{T:bivar}, to show that the RTP
$(\tau_\ibf,\kappa_\ibf,Y_\ibf)_{\ibf\in\T}$ is nonendogenous, it suffices to
show that apart from the trivial fixed point $\ov\rro^{(2)}$, the bivariate
map $T^{(2)}$ has at least one other fixed point in $\Pc(I^2)_\rro$. The
following theorem identifies all fixed points in $\Pc_\ast(I^2)_\rro$. Since
there are precisely two of them, by Lemma~\ref{L:scinv} we conclude that the
nontrivial fixed point is $\un\rro^{(2)}$. The following theorem is proved
in Subsection~\ref{subsection_main_line}.

\begin{theorem}[Nonendogeny]
The \label{T:scalefix} bivariate map $T^{(2)}$ associated with the map $\chi$
from (\ref{chidef}) has precisely two fixed points $\rho^{(2)}_1,\rho^{(2)}_2$
in $\mathcal{P}_\ast(I^2)_\rro$. For each $c \geq 0$, let $f_c:[0,1]\to[0,1]$
denote the continuous function given by the unique
solution to the Cauchy problem
\begin{equation}
\dif{r} f_c(r)= \frac{c r}{f_c(r)-r/2 }, \quad 0 \leq r < 1, \qquad f_c(0)=\ffrac{1}{2}.
\end{equation}
The equation
\begin{equation}
 f_c(1)^2- \ha f_c(1)= 2c
\end{equation}
is solved for precisely two values of $c$ in $\half$. Denoting these by
$c_1$ and $c_2$ with $c_1<c_2$, we have $c_1=0$ and $c_2\in(0,1/4)$.
The measures $\rro^{(2)}_i$ $(i=1,2)$ are uniquely characterised by
\begin{equation}
\label{finite_both_coord}
\rro_i^{(2)} \big( \{ [0,r] \times I \} \cup \{ I \times [0,s] \}  \big) =  (s \vee r) f_{c_i}\left( \frac{r \wedge s }{ r \vee s} \right),  \quad\big((r,s)\in[0,1]^2\beh\{(0,0)\}\big).
\end{equation}
One has $\rho^{(2)}_1=\ov\rro^{(2)}$, the trivial fixed point defined
as in (\ref{ovnu}).
\end{theorem}


\begin{figure}[htb!]
\begin{center}
\includegraphics{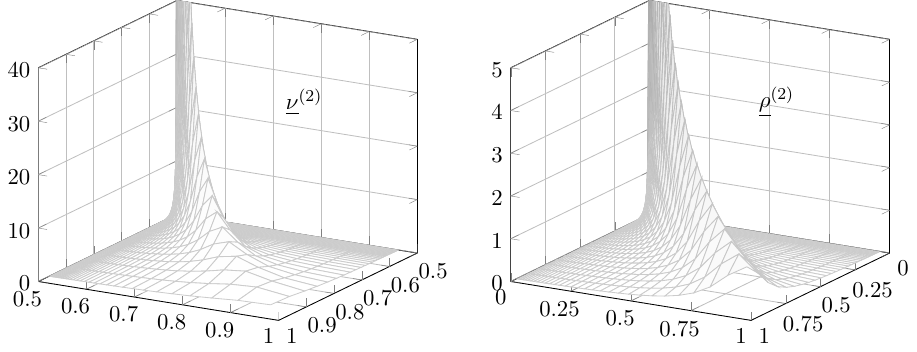}
\caption{The nontrivial solutions $\un\nu^{(2)}$ and $\un\rho^{(2)}$ of the
      bivariate RDE for frozen percolation on the oriented binary tree and the
      MBBT, respectively. Plotted are the densities of the restrictions of the
      measures to $[\ha,1]^2$ and $[0,1]^2$, respectively.}
\label{fig:unnu}
\end{center}
\end{figure}

Numerically, we find $c_2\approx 0.01770838$. The function $f_{c_2}$ is
increasing and convex with $f_{c_2}(0)=\ha$ and $f_{c_2}(1)\approx 0.5629165415$.
Lemma~\ref{L:scinv} allows us to identify $\rro^{(2)}_2$ as $\un\rro^{(2)}$,
the nontrivial fixed point defined in (\ref{unnu}). As a result of
Theorem~\ref{T:scalefix}, we also have an explicit expression for the
nontrivial solution $\un\nnu^{(2)}$ to the bivariate RDE for frozen
percolation on the oriented binary tree, see formula (\ref{explit}) below.
Numerical data for $\un\nnu^{(2)}$ and $\un\rro^{(2)}$ are plotted in
Figure~\ref{fig:unnu}.

\begin{figure}[htb!]
\begin{center}
\includegraphics{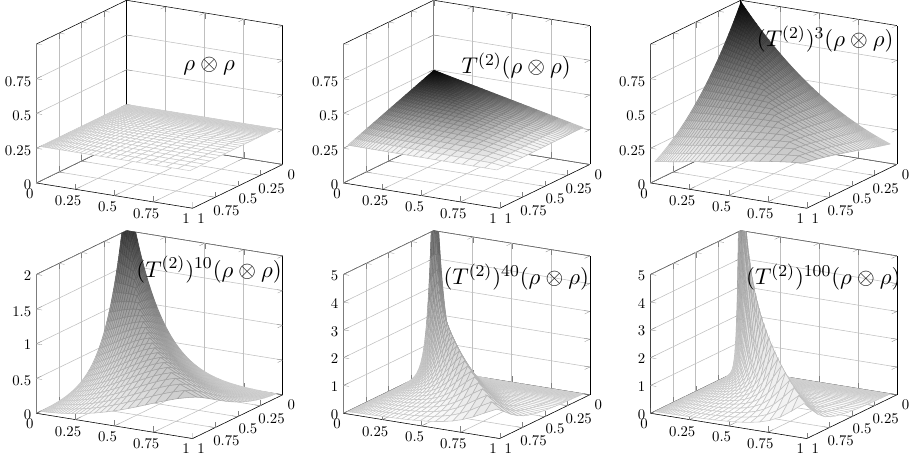}
\caption{Iterating the bivariate map $T^{(2)}$ on the product measure
  $\rho\otimes\rho$ produces a series of measures that by (\ref{prodtoun})
  converge to the nontrivial fixed point $\un\rho^{(2)}$. Plotted is the
  density of the restriction of $(T^{(2)})^n(\rho\otimes\rho)$ to the unit
  square for $n=0,1,3,10,40$, and $100$. The last plot is already very close
  to the theoretical limit.}
\label{fig:rhoiter}
\end{center}
\end{figure}

\subsection{Discussion}

\subsubsection*{Frozen percolation on finite graphs}

Let $G=(V,E)$ be a finite graph. Let $(U_e)_{e\in E}$ be i.i.d.\ uniformly
distributed on $[0,1]$ and let $(\La_v)_{v\in V}$ be an independent
i.i.d.\ collection of exponentially distributed random variables with mean
$\la^{-1}$. Now consider a process where edges and vertices can be in two
possible states: edges are either \emph{closed} or \emph{open}, and vertices
are either \emph{available} or \emph{frozen}. Initially, all edges are closed
and all vertices are available.  The evolution is as follows:
\begin{enumerate}
\item At time $U_e$, the edge $e$ becomes open, provided neither of its
  endvertices is frozen.
\item At time $\La_v$, all vertices of the open component containing $v$
  become frozen.
\end{enumerate}
We call such a process \emph{frozen percolation} on the finite graph $G$, and
by a certain analogy with forest fire models, we call $\lambda$ the
\emph{lightning rate}.

One is typically interested in the limit when $G$ is large. Let us therefore
consider a sequence $G_n=(V_n,E_n)$ of finite graphs with $|V_n|=n$ vertices
and with lightning rates $\la_n$, and make two assumptions:
\begin{itemize}
\item[(A1)] The graphs $G_n$ converge to a weak local limit $G$ in the sense of
  Benjamini and Schramm.
\item[(A2)] $n^{-1}\ll\la_n\ll 1$ as $n\to\infty$.
\end{itemize}
We recall that a sequence of graphs converge to a weak local limit if the
neighbourhood of a typical (uniformly chosen) vertex converges in law to a
(possibly random) rooted graph; see \cite{BS01} or \cite[Section~1.4]{Hof17II}.
Assumption (A2) guarantees that in the limit, small open clusters with size of
order one never freeze, but giant components that occupy a positive fraction
of all vertices freeze immediately.

We can think of this as a model for polymerisation, where open components
represent polymers that grow through merger with neighbours. Polymers that grow
too large become part of the ``gel'' and are unable to grow any further. In
the model we have just described, this is guaranteed by the lightning process,
which has certain mathematical advantages. However, one can also think about
alternative models where polymers are, e.g., prevented from growing when they
reach a certain deterministic size.

If $p_{\rm c}$ is the critical value for percolation in the local limit graph
$G$ from assumption (A1), then up to time $p_{\rm c}$, open clusters grow as
in normal percolation. Since beyond this time, large clusters are prevented
from growing further, one can expect the model to exhibit \emph{self-organised
  criticality} (SOC) in the sense of \cite{B96,J98}, i.e., in the whole time
regime beyond time $p_{\rm c}$ we can expect phenomena that are usually
associated with the behaviour of large systems at their critical
point. Statements of this form have indeed been proved. With the model
described above in mind, we will give a short overview of the literature and
mention some open problems.

\subsubsection*{Frozen percolation on the complete graph}

Although historically not the oldest, frozen percolation on the complete graph
is one of the most natural models to consider. Since in this case, the degree
of each vertex is $n$, it is more natural to take the $(U_e)_{e\in E}$ to be
uniformly distributed on $[0,n]$ instead of $[0,1]$. The complete graph does
not have a weak local limit, but one can take the local limit of the combined
object consisting of the complete graph and the edge activation times
$U_e$. The resulting limiting object is called the Poisson Weighted Infinite
Tree (PWIT) \cite[Sect~4.2]{AS04}.

Following a suggestion in \cite[Sect.~5.5]{Ald00}, one of us has studied frozen
percolation on the complete graph. In \cite{Rat09}, it was shown that
the fraction of clusters of sizes $k\in\N$ at time $t$ converges to a solution
of Smoluchowski's equations with multiplicative kernel, an infinite system of
differential equation that serves as a deterministic model of polymerisation,
and that is known to exhibit self-organised criticality (SOC).

\hspace{-0.4pt}The closely related forest fire model of \cite{RT09} is further studied in
\cite{CFT15,CRY18,C18}. In \cite{CRY18} it is shown that the asymptotic
distribution of a typical cluster is that of a critical multi-type
Galton-Watson tree after gelation.

Aldous \cite[Sect.~5.5]{Ald00} in fact suggested to study the variant of the
mean field frozen percolation model where clusters are frozen when their size
exceeds a deterministic threshold $1\ll\al(n)\ll n$. This model is studied in
in \cite{MN14}. Their Theorem~1.3 states that at any time $t\geq 1$, the
limiting distribution of a typical non-frozen cluster is that of a critical
Galton-Watson tree with Poisson offspring distribution, again establishing SOC.
As an open problem, we mention:
\begin{problem}\label{P:PWIT}
Construct frozen percolation on the PWIT and show that it is the local weak
limit of the models in \cite{Rat09,MN14}.
\end{problem}


\subsubsection*{Coagulation equations}

The relation of frozen percolation on the complete graph to Smoluchowski's
coagulation equations has already been mentioned. A remark of Stockmayer
\cite{S43} on these equations inspired Aldous' work for the 3-regular tree. In
\cite[Section 1.1]{Ald00} Aldous compares the post-gel behaviour of
Smoluchowski's coagulation equations to the self-similar behaviour of his
model. In \cite{Ald99} Aldous surveys the connections between variants of
Smoluchowski's coagulation equations and various stochastic models of
coagulation.

The configuration model \cite{Hof17I} is a well-studied random graph whose
weak local limit is well-known. In particular, one can choose the parameters
of the configuration model so that this limit is the 3-regular or more
generally any $d$-regular graph. The configuration model has a dynamical
construction where to vertices there are assigned ``half-edges'' or ``arms''
that are then randomly linked. In \cite{MN15} a variant of this model is
treated where components freeze once their size exceeds a fixed
threshold. They link the model to a variant of Smoluchowski's equations and it
is shown that after gelation, the asymptotic distribution of a typical
non-frozen cluster is that of a critical Galton-Watson tree.

The mathematical connection between more general stochastic models of
coalescence where clusters with a size above a certain threshold are frozen
and Smoluchowski's equation with more general kernels is established in
\cite{FL09}.

\subsubsection*{Frozen percolation the 3-regular tree}

Aldous' work on frozen percolation on the 3-regular tree is the first example
of a dynamically constructed random graph model that exhibits SOC. In
\cite[Prop~11 and Thm~14]{Ald00} it is proved that at any time $t\in[\ha,1]$,
a typical finite cluster is distributed as a critical percolation cluster on
the binary tree, and infinite clusters are distributed as the incipient
infinite cluster. As an open problem, we mention:
\begin{problem}\label{P:loclim}
Show that frozen percolation on the 3-regular tree is the weak local limit of
frozen percolation on a suitable sequence of finite graphs.
\end{problem}
When proving convergence, it is very useful to have a unique characterization of
the limit. A unique characterization of frozen percolation on the 3-regular
tree is provided by our Theorem~\ref{T:orfrz}. We do not know if
condition~(iii) is in fact needed for uniqueness. Likewise, the following
question is still open:
\begin{question}
Do conditions (i) and (ii) of Theorem~\ref{T:frz} uniquely determine the law
of $(\Ui,F)$?
\end{question}
We note that using Theorem~\ref{T:main}, it is not hard to show that (i) alone
is not sufficient for distributional uniqueness.

Let us note here that a variant of Aldous' frozen percolation model on the
binary tree where clusters with size greater than a large number $N$ are
frozen was introduced in \cite{vdBKN12}. The law of the cluster of the origin
at time $t \in [0,1]$ in the frozen percolation model with freezing threshold
$N$ locally converges to the corresponding law in the frozen percolation
model of Aldous \cite{Ald00} as $N \to \infty$ (see \cite[Theorem
  1]{vdBKN12}).

\subsubsection*{Nonendogeny}

In line with Problem~\ref{P:loclim}, we expect that for a suitable sequence of
finite graphs whose weak local limit is the 3-regular tree, if we couple two
frozen percolation processes on these graphs by using the same edge
activation times but independent lightning processes, then the weak local
limit should be the process $(\Ui,F,F')$ from Theorem~\ref{T:main}. In
particular, the local limit of such processes should be a.s.\ different
because of nonendogeny.

Even though the basic question of endogeny has now been settled for the binary
tree, more detailed questions remain open. In Section~\ref{S:RDEsol}, we
classify all solutions to RDE (\ref{MRDE}). This leads to the question:
\begin{question}
For which solutions of the RDE (\ref{MRDE}) is the corresponding RTP
nonendogenous?
\end{question}
Even for the RTP in Theorem~\ref{T:nonend}, one would like to understand
better what is going on.
\begin{question}
By Theorem~\ref{T:nonend}, the \si-field generated by
$(\tau_\ibf,X_\ibf)_{\ibf\in\T}$ is larger than the \si-field generated by
$(\tau_\ibf)_{\ibf\in\T}$. Give an explicit characterisation of the extra
randomness needed to construct $(X_\ibf)_{\ibf\in\T}$.
\end{question}
In this context, we mention that in \cite[Thm~1.2]{Ban06}, it is proved that
the tail \si-algebra of $(X_\ibf)_{\ibf\in\T}$ is trivial. Proposition 1.1 of
\cite{Ban06} states that generally, endogeny of a RTP implies its
tail-triviality, however our main result exemplifies that the converse
implication does not necessarily hold.

Related to our previous question is the following problem. Let $X'_\wurz$ denote
the first time when there is an infinite path of open or frozen edges starting
at the root. Then clearly $X'_\wurz\leq X_\wurz$ a.s. If the answer to the
following question is positive, then this is all that can be said with
certainty about $X_\wurz$ based on $(\tau_\ibf)_{\ibf\in\T}$.
\begin{question}
Let $\xi:=\P\big[X_\wurz\in\,\cdot\,\big|\,(\tau_\ibf)_{\ibf\in\T}\big]$.
Is it true that the support of $\xi$ is a.s.\ equal to $[X'_\wurz,\infty)$?
\end{question}

In Theorem~\ref{T:scalefix}, we have shown that the bivariate RDE has
precisely two scale invariant fixed points. We believe that there exist fixed
points that are not scale invariant. To see why, recall that we suggested
that $\un\nnu^{(2)}$ should describe the local limit of two finite frozen
percolation processes that use the same edge activation times but independent
lightning processes. We believe that the local limit of two processes that use
the same lightning process up to some time $\ha<s<1$ and independent lightning
processes thereafter should be described by a fixed point of $T^{(2)}$ that is
neither $\un\nnu^{(2)}$ nor $\ov\nnu^{(2)}$.

It has been shown in \cite[Thm~1]{MSS20} that for each initial state, the
differential equation
\be\label{difmu2}
\dif{h}\mu^{(2)}_h=T^{(2)}(\mu^{(2)}_h)-\mu^{(2)}_h\qquad(h\geq 0)
\ee
has a unique solution.
\begin{problem}\label{P:attract}
For frozen percolation on the oriented binary tree, find all fixed points of
(\ref{difmu2}) and their domains of attraction.
\end{problem}
In \cite[Prop~12]{MSS20} Problem~\ref{P:attract} is solved for a different
RTP, which is also nonendogenous. In that example, $\un\nu^{(2)}$ and
$\ov\nu^{(2)}$ turned out to be the only fixed points, where the trivial fixed
point $\ov\nu^{(2)}$ is unstable and the nontrivial fixed point
$\un\nu^{(2)}$ attracts all other initial states. One wonders if the situation
for frozen percolation is similar. In general, we ask:
\begin{question}
For a general RTP, can one prove nonendogeny by proving that the trivial fixed
point $\ov\nu^{(2)}$ is unstable?
\end{question}
Our proof of Theorem~\ref{T:nonend} is based on an explicit formula for
$\un\nu^{(2)}$. Ultimately, one would like to be able to prove nonendogeny
without having to solve the bivariate RDE.

\subsubsection*{Frozen percolation on regular trees}

In Problem~\ref{P:PWIT}, we have already mentioned frozen percolation on the
PWIT. Aldous \cite[Sect.~5.4]{Ald00} observed that his construction can be
carried out on any $d$-regular tree, and even gave a formula for the
distribution of freezing times on $d$-regular trees. This leads to:
\begin{question}\label{question_d_regular_nonendo}
Are frozen percolation on the PWIT or on general $d$-regular trees endogenous?
\end{question}
We conjecture the answer to this question to be negative, but this does not
follow from the methods of this paper. Our main results are for the MBBT and
essentially rely on the nice scaling property of the latter that simplifies
our formulas. The fact that we are also able to treat the oriented binary tree
and consequently the unoriented 3-regular tree depends on a trick that uses in
an essential way that the MBBT is a binary tree.

Nevertheless, we hope that our methods will be useful in answering
Question~\ref{question_d_regular_nonendo}. The reason for this optimism is
that the MBBT can be seen as the near-critical scaling limit of percolation on
a wide class of oriented trees, such as oriented $d$-ary trees or the PWIT.

Indeed, since edges with $U_e\leq p_{\rm c}$ belong to finite clusters when
they open, from the point of view of frozen percolation it does not matter
when they open. In view of this, let us focus only on those edges whose
activation times lie between $p_{\rm c}$ and $p_{\rm c}+\eps$ for some small
$\eps>0$. If we condition on the event that there is an infinite path starting
at the root along edges with activation times $U_e\leq p_{\rm c}+\eps$, and
cut off all parts of the tree that do not lie on such an infinite path, then
the scaling limit as $\eps \to 0$ of our tree $\Ti$, and the locations marked
with the (scaled) activation times of edges with $p_{\rm c}<U_e<p_{\rm
  c}+\eps$ converge to the marked Poisson process $\Pi$ on~$\Ti$.

In view of this, we expect that on a general class of oriented trees, frozen
percolation is nonendogenous and the nontrivial fixed point $\un\nu^{(2)}$ of
the bivariate RDE will in a small neighbourhood of the critical point look
similar to the nontrivial fixed point from Theorem~\ref{T:scalefix}.

\subsubsection*{Frozen percolation on integer lattices}

One can try to ``naively'' define frozen percolation on any infinite graph as
in property~(i) of Theorem~\ref{T:frz}, by specifying that clusters stop
growing as soon as they reach infinite size. It is an observation of Benjamini
and Schramm that such a process cannot be defined on the planar square lattice
(for a sketch of a proof, see \cite[Section 3]{vdBT01}). The following
question is open:
\begin{question}
For which $d\geq 3$ does there exists a frozen percolation process on the
nearest-neighbour lattice $\Z^d$ that satisfies property~(i) of
Theorem~\ref{T:frz}?
\end{question}

There exists an extensive literature for finite versions of frozen percolation
on the planar lattice. A model where clusters with diameter greater than a
large number $N$ are frozen was introduced in \cite{vdBLN12}.
The behaviour of this model is rather different from the
the analogous model of \cite{vdBKN12} on the binary tree that we have discussed after
Problem~\ref{P:loclim}, because in planar diameter-frozen percolation all
frozen clusters freeze in the critical time window around the Bernoulli
percolation threshold $p_c$, the frozen clusters look similar to critical
percolation clusters, and moreover macroscopic non-frozen clusters
asymptotically have full density as $N \to \infty$, c.f.\ \cite{vdBLN12, K15}.
In \cite{vdBLN17} it is shown that the particular mechanism to freeze clusters
(the ``boundary rules'') matters strongly, i.e., if we modify the
diameter-frozen site percolation model on the triangular lattice in a way that
the outer boundary of frozen connected components can become occupied (and
later freeze) then frozen clusters in the terminal configuration have
non-vanishing density as $N \to \infty$.

The percolation on the planar lattice where clusters with volume (cardinality) greater than a large number $N$ are frozen was introduced in \cite{vdBN17}, the main result being that if we restrict the process to a large box with side-length $n$, then the probability that the origin freezes depends on the relation between $N$ and $n$ in an oscillatory fashion. Thus the behaviour of the volume-frozen process is substantially different from that of the diameter-frozen process. In  \cite{vdBKN18} it is shown that in the volume-frozen model many frozen clusters surrounding the origin appear successively, each new cluster having a diameter much smaller than the previous one. In \cite{vdBKN18} it is also proved that in the full planar case ($n=\infty$) with high probability (as $N \to \infty$), the origin does not belong to a frozen cluster in the final configuration.
 In  \cite{vdBN18} it is proved that if the freezing mechanism in a box of
 size $n$ is governed by independent lightnings hitting the vertices then the
 density of frozen sites depends on the relation between the lightning rate
 and $n$  in an oscillatory fashion.



\subsubsection*{Self-destructive percolation and forest fire model on infinite graphs}

The ``naive'' definition of the forest fire model on an infinite graph  $G=(V,E)$ (dating back to  \cite{DS92}) is  as follows:
vacant sites become occupied at rate $1$ and infinite occupied clusters become vacant instantaneously. Similarly to the case of the frozen percolation model,
it is a highly non-trivial question whether such a process exists.

The model of \emph{self-destructive percolation} was introduced by \cite{vdBB04} in order to address this question on the planar lattice: given some $p>p_c$, let us switch all of the sites which are in an infinite occupied component
into vacant  state (destruction) and then turn any vacant site occupied with probability $\delta$ (enhancement). Denote by $\delta(p)$ the smallest enhancement needed
for the appearance of an infinite cluster in the enhanced configuration. Theorem 4.1 of \cite{vdBB04} states that if $\lim_{p \searrow p_c} \delta(p) >0$ then the forest fire process cannot be defined on the planar lattice.  Theorem 1 of \cite{KMS15} states that indeed $\lim_{p \searrow p_c} \delta(p) >0$ on the planar lattice.
However, we have $\lim_{p \searrow p_c} \delta(p) =0$ on  non-amenable graphs \cite{AST14} and $\mathbb{Z}^d$ for high enough $d$ \cite{ADKS15}.
Also note that in mean field percolation models we have $\lim_{p \searrow p_c} \frac{\delta(p)}{p-p_c}=1$. This asymptotic relation becomes an exact equality of the lengths of growth and recovery time intervals if one considers self-destructive (frozen) percolation on the MBBT, moreover  the general  solutions to the RDE (\ref{MRDE}) (c.f.\ Section \ref{S:RDEsol}) and the associated RTP's (c.f.\ Section \ref{S:MBBTfrz}) also exhibit time intervals of (supercritical) growth and (subcritical) recovery, which are of equal length.

Currently it is an open question whether it is possible to define a forest fire process on the nearest-neighbour lattice $\mathbb{Z}^d, \, d \geq 3$.
In \cite{vdBT01} a variant of the forest fire model (with site-dependent occupation rates) is constructed on the half-line.
The construction of the variant of the forest fire model with a positive rate of lightning per vertex on $\mathbb{Z}^d$ is given in \cite{D06a, D06b}: if a lightning hits a vertex $v$, then all of the sites in the occupied cluster of $v$ become vacant instantaneously.
In \cite{G14,G16} a variant of the forest fire model on the half-plane is defined where components that touch the boundary (or become infinite) are destroyed.
It is shown that before (and including) the critical time, the effect of the destruction mechanism is only felt locally near the boundary of the half-plane, whereas after the critical time, it is felt globally on the entire half-plane.


\subsubsection*{Outline}

The rest of the paper is devoted to proofs. We prove
Theorem~\ref{T:scalefix} in Section~\ref{S:RDE} and the
remaining results in Section~\ref{S:frz}. The paper concludes with a small
appendix on skeletal branching processes, which are related to the scaling
property of the MBBT described in Proposition~\ref{P:scale}.

Even though Theorem~\ref{T:nonend}, which is proved in
Subsection~\ref{S:RDEequiv}, is our main result, considerable extra effort is
needed to prove additional results, in particular, uniqueness of the
nontrivial fixed point in Theorem~\ref{T:scalefix} and its subsequent
identification as $\un\rho^{(2)}$ with the help of Lemma~\ref{L:scinv}, as
well as Theorem~\ref{T:orfrz}, which depends on the classification of general
solutions to the RDE (\ref{MRDE}) in Subsection~\ref{S:RDEsol}.

\section{The bivariate RDE}\label{S:RDE}

\subsection{Main line of the proof}\label{subsection_main_line}

In this section, we prove
Theorem~\ref{T:scalefix}. The
main steps of the proof
are summarised in the following
lemmas. We first need a convenient way to parametrise elements of the space
$\Pc_\ast(I^2)_\rro$.

\bl[Parametrisation of the space of interest]
For\label{L:rhoF} each $\rho^{(2)}\in\Pc_\ast(I^2)_\rro$, there exists a
unique continuous function $f:[0,1]\to\R$ such that
\be\label{rhoF}
\rho^{(2)} \big( \{ [0,r] \times I \} \cup \{ I \times [0,s] \}  \big) =  (s \vee r) f\left( \frac{r \wedge s }{ r \vee s} \right),  \quad\big((r,s)\in[0,1]^2\beh\{(0,0)\}\big),
\ee
and such a function $f$ uniquely characterizes $\rho^{(2)}$.
In particular, the trivial fixed point $\ov\rro^{(2)}$ corresponds to $\ov
f(r)=\ha$, $(r\in[0,1])$.
\el

There are a priori many ways of parametrising elements of
$\Pc_\ast(I^2)_\rro$. The para\-met\-risation in terms of the function $f$ from
(\ref{rhoF}) turns out to lead to
a particularly simple form of the bivariate RDE.

\bl[Bivariate RDE]
An\label{L:bifix} element $\rho^{(2)}\in\Pc_\ast(I^2)_\rro$ is a fixed point
of the bivariate map $T^{(2)}$ associated with the map $\chi$
from (\ref{chidef}) if and only if the function $f:[0,1]\to\R$ from
(\ref{rhoF}) is continuously differentiable on $[0,1)$ and satisfies
\be\ba{l}\label{bifix}
\dis{\rm(i)}\ \dif{r}f(r)=\frac{c r}{f(r)-r/2  } \qquad\big(r\in[0,1)\big),\\[8pt]
\dis{\rm(ii)}\ f(0)=\ha,\qquad
{\rm(iii)}\  f(1)^2- \ha f(1) =2c ,
\ec
for some $c\geq 0$.
\el

In particular, the trivial fixed point $\ov f(r)=\ha$ solves
(\ref{bifix}) with $c=\ov c:=0$.
The following lemma shows that there is exactly one
other, nontrivial solution.

\bl[Nontrivial solution of (\ref{bifix})]
For\label{L:fixsol} each $c\geq 0$, there exists a unique solution $f_c$ to
(\ref{bifix})~(i) and (ii). There exists a unique $c_2>0$ such that the
function $f_{c_2}$ also satisfies (\ref{bifix})~(iii). Moreover, we have
$c_2 \in (0, \ffrac{1}{4} )$.
\el

In Lemma~\ref{L:rhoF}, we have shown that a probability law
$\rho^{(2)}\in\Pc_\ast(I^2)_\rro$ is uniquely characterised by the
corresponding function $f$ from (\ref{rhoF}), but we have not given sufficient
conditions for a function $f:[0,1]\to\R$ to correspond to an element of
$\Pc_\ast(I^2)_\rro$. In view of this, to complete the proof of
Theorem~\ref{T:scalefix}, we need one more lemma.

\bl[Nontrivial solution of the bivariate RDE]
The\label{L:nontriv} function $f_{c_2}$ from Lemma~\ref{L:fixsol} defines
through \eqref{rhoF} a probability measure
$\rho_2^{(2)}\in\Pc_\ast(I^2)_\rro$. The restriction of $\rho_2^{(2)}$ to
$[0,1]^2$ has a density w.r.t.\ the Lebesgue measure. In particular,
$\rho_2^{(2)}$ puts no mass on the diagonal $\big\{(r,r):r\in[0,1]\big\}$.
\el

\bpro[Proof of Theorem \ref{T:scalefix}]
By Lemmas~\ref{L:rhoF}, \ref{L:bifix}, \ref{L:fixsol}, and \ref{L:nontriv},
the bivariate map $T^{(2)}$ has, apart from the trivial fixed point
$\ov\rho^{(2)}$, precisely one more fixed point in $\Pc_\ast(I^2)_\rho$, which
is given as in (\ref{finite_both_coord}) in terms of the function $f_{c_2}$.
\epro

We will prove Lemmas 
\ref{L:rhoF}, \ref{L:bifix}, \ref{L:fixsol} and \ref{L:nontriv} in Sections 
\ref{subsection_proof_of_lemma_rhof}, \ref{subsection_proof_of_lemma_bifix}, \ref{subsection_proof_of_lemma_fixsol} and \ref{subsection_p_of_l_nontriv}, respectively.

%
%
%

\subsection{Parametrisation of scale invariant measures}
\label{subsection_proof_of_lemma_rhof}

In this subsection we prove Lemma~\ref{L:rhoF}. We also prepare for the proof
of Lemma~\ref{L:nontriv} by giving sufficient conditions for a function
$f:[0,1]\to\R$ to define a measure $\rho^{(2)}\in\Pc_\ast(I^2)_\rho$ through
(\ref{rhoF}).

\bl[Encoding $\rho^{(2)}$ as a bivariate function] \label{lemma_fixed_point_G_char}
Any $\rho^{(2)} \in \mathcal{P}(I^2)_\rro $ is uniquely characterised by the
continuous function $F:[0,1]^2\to[0,1]$ defined as
 \begin{equation}\label{Frs}
 F(r,s):=\rho^{(2)} \big(  \{ [0,r] \times I\} \cup \{ I \times [0,s] \}  \big), \qquad\big(r,s\in(0,1]\big).
 \end{equation}
Moreover, $\rho^{(2)}\in\Pc_\ast(I^2)_\rho$ if and only if
$F$ is symetric in the sense that $F(r,s)=F(s,r)$ and
\be\label{Fscale}
F(tr,ts)=tF(r,s)\qquad\big(r,s,t\in[0,1]\big).
\ee
\el

\begin{proof}
Since both marginals of $\rho^{(2)}$ are equal to $\rho$, formula
(\ref{Frs}) is equivalent to
\be\ba{r@{\ }r@{\,}c@{\,}l}\label{Fro}
{\rm(i)}&
\dis\rho^{(2)}\big(\{\infty\}\times\{\infty\}\big)&=&\dis 1-F(1,1),\\[5pt]
{\rm(ii)}&
\dis\rho^{(2)}\big([0,r]\times\{\infty\}\big)&=&\dis F(r,1)-\ha,\\[5pt]
{\rm(iii)}&
\dis\rho^{(2)}\big(\{\infty\}\times[0,s]\big)&=&\dis F(1,s)-\ha,\\[5pt]
{\rm(iv)}&
\dis\rho^{(2)}\big([0,r]\times[0,s]\big)&=&\dis\ha r+\ha s-F(r,s).
\ec
Since these functions uniquely determine the restrictions of $\rho^{(2)}$ to
$\{(\infty,\infty)\}$, $[0,1]\times\{\infty\}$, $\{\infty\}\times[0,1]$, and
$[0,1]^2$, the function $F$ determines $\rho^{(2)}$ uniquely. Moreover, we see
from (\ref{Fro})~(iv) that $\rho^{(2)}$ is scale invariant in the sense of
(\ref{scinv}) if and only if (\ref{Fscale} holds.
Since the marginals of $\rho^{(2)}$ are equal to $\rho$, and $\rho$ has no
atoms in $[0,1]$, we see from (\ref{Fro})~(iv) that $F$ is a continuous
function.
\end{proof}

If a closed subset $C$ of $\R^d$ is the closure of its interior, then we say
that a function is $n$ times continuously differentiable on $C$ is all partial
derivatives up to $n$-th order exist on the interior of $C$ and can be
extended to continuous functions on $C$.

\bl[Sufficient conditions on $F$ corresponding to $\rho^{(2)} \in \mathcal{P}(I^2)_\rho$]
Let\label{L:suffP} $\De:=\big\{(r,s)\in[0,1]^2:0\leq r\leq s\big\}$ and let
$F:\De\to\half$ be a twice continuously differentiable function such that:
\[\ba{c}
{\rm(i)}\ F(1,1)\leq 1,\quad{\rm(ii)}\ F(0,s)=\ha s,
\quad{\rm(iii)}\ r\mapsto F(r,1)\mbox{ is nondecreasing,}\\[5pt]
{\rm(iv)}\ \dif{r}F(r,s)\big|_{r=s}=\dif{s}F(r,s)\big|_{r=s},\qquad
{\rm(v)}\ g(r,s):=-\difif{r}{s}F(r,s)\geq 0.
\ea\]
Extend $F$ and $g$ to $[0,1]^2$ by setting $F(s,r):=F(r,s)$ and $g(s,r):=g(r,s)$ for $\big((r,s)\in\De\big)$.
Then there exists a unique probability measure $\rho^{(2)}\in\Pc(I^2)_\rho$
such that (\ref{Frs}) holds, and the restriction of $\rho^{(2)}$ to $[0,1]^2$
has density $g$ with respect to the Lebesgue measure.
\el

\bpro
Uniqueness follows from Lemma~\ref{lemma_fixed_point_G_char}. By (\ref{Fro}),
condition~(i) guarantees that the mass at $(\infty,\infty)$ is nonnegative,
while conditions~(ii) and (iii) guarantee that the restrictions of $\rho^{(2)}$
to $[0,1]\times\{\infty\}$ and $\{\infty\}\times[0,1]$ are nonnegative
measures.

To complete the proof, we will show that conditions~(ii), (iv) and
(v) imply that (\ref{Fro})~(iv) defines a measure on $[0,1]^2$ with density $g$.
Equivalently, we must show that
\be
D(r,s):=\int_0^r\di r'\int_0^s\di s'\,g(r',s')-\ha r-\ha s+F(r,s)
\qquad\big((r,s)\in\De\big)
\ee
is identically zero. Conditions (ii), (iv) and (v) imply that
\be
D(0,s)=0,\quad\dif{r}D(r,s)\big|_{r=s}=\dif{s}D(r,s)\big|_{r=s},\quand
\difif{r}{s}D(r,s)=0
\ee
$((r,s)\in\De)$ The third equality implies that $D(r,s)=u(r)+v(s)$ for some
differentiable functions $u$ and $v$, but then $D(0,s) \equiv 0$ implies that
$u(0)+v(s) \equiv 0$, thus $v$ is constant and therefore
$u'(r)=\dif{r}D(r,s)\big|_{r=s}=\dif{s}D(r,s)\big|_{r=s}=v'(r)=0$, so $u$ is
also a constant, so $D(r,s)=0$ for any $0 \leq r \leq s \leq 1$.
\epro

\begin{proof}[Proof of Lemma \ref{L:rhoF}]
Given $\rho^{(2)}\in\Pc_\ast(I^2)_\rro$, let $F$ be as in (\ref{Frs}) and let
$f:[0,1]\to\R$ be the continuous function defined by
\be\label{f_rho_prob_def}
f(r):=F(r,1)=\rho^{(2)}\left( \left\{ [0,r]\times I   \right\} \cup \left\{ I \times [0,1]   \right\}   \right), \quad 0 \leq r \leq 1.
\ee
Then $F(r,s)=sf(r/s)$ $(s\neq 0)$ by (\ref{Fscale}) and hence \eqref{rhoF}
follows by symmetry. The fact that $f$ uniquely characterizes the measure $\rho^{(2)}$ follows from \eqref{rhoF} and Lemma \ref{lemma_fixed_point_G_char}.
The trivial fixed point $\overline{\rro}^{(2)}$ of $T^{(2)}$ is the distribution of $(Y,Y)$, where $Y \sim \rro$. In this case $\ov f(r)=  \mathbb{P}( \{ Y \leq r \} \cup \{ Y \leq 1 \} )= \ha$ for any $r \in [0,1]$.
\end{proof}

\bl[Sufficient conditions on $f$ corresponding to $\rho^{(2)} \in \mathcal{P}_*(I^2)_\rho$]
Let\label{L:suffPast} $f:[0,1]\to\R$ be a twice continuously differentiable
function such that
\[\ba{c}
{\rm(i)}\ f(1)\leq 1,\quad{\rm(ii)}\ f(0)=\ha,
\quad{\rm(iii)}\ r\mapsto f(r)\mbox{ is nondecreasing,}\\[5pt]
{\rm(iv)}\ 2f'(1)=f(1),\qquad
{\rm(v)}\ f''(r)\geq 0\quad(r\in[0,1]).
\ea\]
Then there exists a unique probability measure $\rho^{(2)}\in\Pc_\ast(I^2)_\rho$
such that (\ref{rhoF}) holds, and the restriction of $\rho^{(2)}$ to $[0,1]^2$
has a density with respect to the Lebesgue measure.
\el

\bpro
For $0\leq r\leq s$, define $F(r,s):=sf(r/s)$ if $s\neq 0$ and $:=0$
otherwise, and $F(s,r):=F(r,s)$. Then (\ref{rhoF}) is equivalent to
(\ref{Fro}) so uniqueness follows from Lemma~\ref{lemma_fixed_point_G_char}.
Since $F$ is symmetric and satisfies (\ref{Fscale}), the same lemma shows that
if $\rho^{(2)}$ exists, then $\rho^{(2)}\in\Pc_\ast(I^2)_\rho$.

To get existence, we apply Lemma~\ref{L:suffP}. We claim that
conditions~(i)--(v) of that lemma follow from the corresponding
conditions of the present lemma. This is trivial for conditions~(i)--(iii).
Condition~(iv) of Lemma~\ref{L:suffP} yields
\be
f'\big(\frac{r}{s}\big)=f\big(\frac{r}{s}\big)
-\frac{r}{s}f'\big(\frac{r}{s}\big)\qquad(r=s),
\ee
which corresponds to the present condition~(iv). Finally, condition~(v) of
Lemma~\ref{L:suffP} requires that
\be
-\difif{r}{s}sf\big(\frac{r}{s}\big)=-\dif{s}f'\big(\frac{r}{s}\big)
=\frac{r}{s^2}f''\big(\frac{r}{s}\big)\geq 0,
\ee
which corresponds to the present condition~(v).
\epro


\subsection{Bivariate RDE and controlled ODE }
\label{subsection_proof_of_lemma_bifix}

In this subsection we prove Lemma \ref{L:bifix}, i.e., we equivalently
reformulate the bivariate fixed point property $T^{(2)} \rho^{(2)}=\rho^{(2)}$
for a scale invariant measure $\rho^{(2)} \in \mathcal{P}_\ast(I^2)_\rro $ as
the controlled ODE problem \eqref{bifix} for the function $f$ defined in
\eqref{f_rho_prob_def}. We start by deriving an integral expression for the
map $T^{(2)}$. Equation (\ref{F_recursion}) below is an
adaptation of equations (11) and (12) of \cite{Ban04} to the MBBT.

\bl[Bivariate map]
Let\label{lemma_G_int_eq} $\rho^{(2)}\in\Pc(I^2)_\rho$, let $T^{(2)}$ denote
the bivariate map defined as in (\ref{bivdef}) for the map $\chi$ in
(\ref{chidef}), taking uniformly distributed $\tau_\wurz,\kappa_\wurz$ as its
input. Let $F$ be defined in terms of $\rho^{(2)}$ as in (\ref{Frs}) and let
$\widetilde{F}$ be defined similarly in terms of $T^{(2)}(\rho^{(2)})$.
Then
\begin{multline}\label{F_recursion}
 \widetilde{F}(r,s)=F(r,s) - \frac{F(r,s)^2}{2}+\frac{r^2}{8}+\frac{s^2}{8} \\ +\ha \int_0^{r \wedge s} \left( F(r,s)-F(t,s)-F(r,t)+F(t,t) \right) \, \mathrm{d}t, \quad
 r,s \in (0,1].
 \end{multline}
\el
\begin{proof}
Let $\tau_\varnothing$ and $\kappa_\varnothing$ denote independent random variables, where $\tau_\varnothing \sim \mathrm{Uni}[0,1]$ and
$\kappa_\varnothing$ is uniformly distributed on $\{1,2\}$.
Let  $(Y_1,Y_1^*)$ and $(Y_2,Y_2^*)$ denote $I^2$-valued random variables
with distribution $\rho^{(2)}$, independent from each other and of $\tau_\varnothing$, $\kappa_\varnothing$.
Let us define
\begin{equation}\label{def_recursion_of_random_variables}
  Y_\varnothing := \chi[\tau_\varnothing,\kappa_\varnothing](Y_1,Y_2), \qquad  Y^*_\varnothing := \chi[\tau_\varnothing,\kappa_\varnothing](Y^*_1,Y^*_2),
\end{equation}
where $\chi$ is defined in \eqref{chidef}. Then the distribution of
$(Y_\varnothing, Y^*_\varnothing)$ is $T^{(2)}(\rho^{(2)})$. It follows that
\begin{equation}\label{kappa_1_or_2}
\widetilde{F}(r,s)= \ha \mathbb{P}[\, Y_\varnothing \leq r\, \text{ or }  \, Y_\varnothing^* \leq s \, | \, \kappa_\varnothing=1  ]+
\ha \mathbb{P}[\, Y_\varnothing \leq r \, \text{ or } \,  Y_\varnothing^* \leq s \, | \, \kappa_\varnothing=2  ].
\end{equation}
Here
\begin{multline}\label{kappa_1_case}
  \mathbb{P}[\, Y_\varnothing \leq r \, \text{ or }\,  Y_\varnothing^* \leq s \, | \, \kappa_\varnothing=1  ]\\
  \stackrel{(\ref{chidef}) }{=}
 \int_0^1
  \mathbb{P}[\,  Y_1 \in(t, t \vee r] \, \text{ or } \,  Y_1^* \in(t, t \vee s]   \,) \,\mathrm{d}t\\
   =\int_0^1 \mathbb{P}[\, Y_1 \in(t, t \vee r]\,] \,\mathrm{d}t + \int_0^1
     \mathbb{P}[\, Y_1^* \in(t, t \vee s]\,] \,\mathrm{d}t \\
  - \int_0^1
  \mathbb{P}[\,  Y_1 \in(t, t \vee r], \, Y_1^* \in(t, t \vee s] \,  ]\,\mathrm{d}t \\
  \stackrel{(\ref{rhodef}) }{=}\int_0^r \ha (r-t)\, \mathrm{d}t +\int_0^s \ha (s-t)\, \mathrm{d}t- \int_0^{r \wedge s } \mathbb{P}[\,  Y_1 \in(t,  r], \, Y_1^* \in(t,   s] \,  ]\,\mathrm{d}t \\
  \stackrel{(*)}{=}  \frac{r^2}{4} + \frac{s^2}{4} - \int_0^{r \wedge s } \left(  F(t,s)+F(r,t) -F(t,t) -F(r,s) \right) \, \mathrm{d}t,
\end{multline}
where in $(*)$ we used \eqref{Frs} and inclusion-exclusion. Moreover
\begin{multline}\label{kappa_2_case}
\mathbb{P}[\, Y_\varnothing \leq r \,\text{ or }\,  Y_\varnothing^* \leq s \, | \, \kappa=2  ]
 \stackrel{(\ref{chidef}) }{=}
 1- \mathbb{P}[\, Y_1 \wedge Y_2 > r, \, Y_1^* \wedge Y^*_2  > s \,  ]\\
 =1- \mathbb{P}[\, Y_1  > r, \, Y_1^* > s \,  ]^2=
 1-  (1-F(r,s))^2= 2F(r,s)-F(r,s)^2.
\end{multline}
Now \eqref{F_recursion} follows as a combination of \eqref{kappa_1_or_2}, \eqref{kappa_1_case} and \eqref{kappa_2_case}.
\end{proof}



\bl[Scale invariant bivariate fixed point] \label{map_is_homogeneous} $\rho^{(2)} \in
  \mathcal{P}_*(I^2)_\rro$ satisfies $T^{(2)}(\rho^{(2)}) = \rho^{(2)}$ if and
  only if the function $f$ in (\ref{rhoF}) satisfies
 \begin{equation}\label{int_eq_simplified}
   f(r)^2= \frac{1}{4} +r f(r) - \int_0^r f(u) \, \mathrm{d}u + \left( \frac{1}{4} +\frac{f(1)}{2}  -\int_0^1 f(s)\, \mathrm{d}s \right) r^2, \;\;  r \in [0,1].
 \end{equation}
\el
\begin{proof}
If $\rho^{(2)}\in\Pc(I^2)_\rho$ then $\rho^{(2)}$ is symmetric, so
Lemmas~\ref{lemma_fixed_point_G_char} and \ref{lemma_G_int_eq} imply
that $T^{(2)}(\rho^{(2)}) = \rho^{(2)}$ holds if and only
if for any  $0<r \leq s \leq 1$
\begin{equation}\label{F_int_fixed}
   F(r,s)^2=\frac{r^2}{4}+\frac{s^2}{4} + \int_0^r \left( F(r,s)-F(t,s)-F(t,r)+F(t,t) \right) \, \mathrm{d}t.
 \end{equation}
If $\rho^{(2)}\in\Pc_\ast(I^2)_\rho$, then the function $F$ from
(\ref{Frs}) can be expressed in the function $f$ from (\ref{rhoF}) as
 \begin{equation}\label{F_scale_inv_char}
F(r,s)=s f\left( \frac{r}{s} \right),
\qquad 0<r \leq s \leq 1.
 \end{equation}
Plugging this into \eqref{F_int_fixed}, dividing both sides by $s^2$ and using the substitution
$u=t/s$ in the integral we obtain
 \begin{equation}\label{scale_inv_int_eq_complicated}
  f\left(\frac{r}{s}\right)^2= \frac{(r/s)^2}{4}+ \frac{1}{4}  +
  \int_0^{r/s} \left( f\left(\frac{r}{s}\right)-f(u)-\frac{r}{s}f\left( \frac{u}{r/s} \right) + uf(1) \right) \, \mathrm{d}u,
\end{equation}
which holds for all $0< r \leq s \leq 1$ if and only if
\begin{equation}\label{scale_inv_integral_eq}
  f(r)^2= \frac{1}{4} + \frac{r^2}{4}  +
  \int_0^r \left( f(r)-f(u)-rf\left( \frac{u}{r} \right) + uf(1) \right) \, \mathrm{d}u, \quad
 0 < r \leq 1.
\end{equation}
Evaluating the integrals, using the substitution $s=u/r$,
we arrive at (\ref{int_eq_simplified}), which also
holds for $r=0$ since $f$ is continuous.
 \end{proof}

\begin{remark}\label{remark_constant_f_triv}
For any $\rho^{(2)}\in\Pc_\ast(I^2)_\rho$, setting $s=1$ in (\ref{rhoF})
yields (\ref{f_rho_prob_def}), which shows that $f$ is nondecreasing. Since
the marginals of $\rho^{(2)}$ are $\rho$, we have $f(0)=\ha$.
If $f(1)=\ha$, then we must have $f(r)=\ha, \, r \in [0,1]$.
In this case \eqref{int_eq_simplified} holds. This is the
$\overline{f}$ associated to the (scale invariant) diagonal fixed point
$\overline{\rro}^{(2)}$ of $T^{(2)}$.
\end{remark}

\bl[Controlled ODE] \label{lemma_from_int_to_ode}
Let $f:[0,1]\to \half$ be continuous and nondecreasing with $f(1) > \ha$. Then $f$
satisfies (\ref{int_eq_simplified}) if and only if $f$ is continuously
differentiable and solves
\begin{equation}\label{diff_eq_for_f_fix}
{\rm(i)}\ f(0)=\frac{1}{2},
\quad{\rm(ii)}\  f'(r)=\frac{c r}{ f(r)-r/2}, \quad r \in [0,1],
\end{equation}
where
\begin{equation}\label{diff_bc}
c=  \frac{1}{4} +\frac{f(1)}{2}  -\int_0^1 f(s)\, \mathrm{d}s >0.
\end{equation}
\el
\begin{proof}
Plugging in $r=0$ into \eqref{int_eq_simplified} we obtain $f(0)=\ffrac{1}{2}$.
Using this, we have
\begin{equation}\label{stieltjes_int}
{\rm(i)}\  f(r)^2=\frac{1}{4} + \int_0^r 2 f(u)\mathrm{d}f(u),
\quad{\rm(ii)}\  rf(r) -\int_0^r f(u) \, \mathrm{d}u = \int_0^r u \, \mathrm{d}f(u),
\end{equation}
where both integrals in \eqref{stieltjes_int} are Stieltjes.
Inserting this into \eqref{int_eq_simplified} yields
\begin{equation}\label{int_eq_simp_stieltjes}
  \int_0^r (2 f(u)-u) \, \mathrm{d} f(u)=c r^2, \quad 0 \leq r \leq 1,
\end{equation}
with $c$ as in (\ref{diff_eq_for_f_fix}).
Since $f$ is nondecreasing with $f(0)=\ha$, we observe that $2 f(u)-u \geq 2
f(0)-u >0$ for all $u \in [0,1)$. Combining this with the assumption
$f(1)>\ha$ we get $m:=\min_{0\leq u \leq 1}\left( 2 f(u)-u \right) >0$, since
$u \mapsto 2 f(u)-u$ is a continuous function on the compact interval
$[0,1]$. Since the right-hand side of \eqref{int_eq_simp_stieltjes} is has
Lipschitz constant $2c$, we conclude that $f$ is
$2c/m$-Lipschitz-continuous on $[0,1]$. Thus, by the Radon-Nykodim theorem,
there exists a Lebesgue-a.s.\ unique measurable function $f^\circ:[0,1] \to
[0,2c/m]$ such that $\int_0^r f^\circ(u) \, \mathrm{d}u= f(r)-\ffrac{1}{2}$
for all $0 \leq r \leq 1$.

By \eqref{int_eq_simp_stieltjes} we have $\int_0^r (2 f(u)-u) \, f^\circ(u)\,
\mathrm{d}u =c r^2$ for any $0\leq r \leq 1$, thus $(2 f(r)-r)f^\circ(r)=2 c
r$ for Lebesgue-almost all $r\in [0,1]$, from which it follows that the
Radon-Nykodim derivative $f^\circ$ can be chosen to be the continuous function
$f^\circ(r)=\frac{2 c r}{2 f(r)-r}$, therefore $f$ is continuously
differentiable, $f'=f^\circ$ and \eqref{diff_eq_for_f_fix} holds.

Since $f$ is nondecreasing, (\ref{diff_eq_for_f_fix}) implies $c\geq
0$. Solving (\ref{diff_eq_for_f_fix}) with $c=0$ yields $f(1)=\ha$,
contradicting $f(1)>\ha$, so we conclude that $c>0$.

Assume, conversely, that $f$ solves (\ref{diff_eq_for_f_fix}) and
(\ref{diff_bc}). Then (\ref{diff_eq_for_f_fix})~(ii) implies
(\ref{int_eq_simp_stieltjes}) and (\ref{diff_eq_for_f_fix})~(i) yields
(\ref{stieltjes_int})~(i). Combining this with (\ref{stieltjes_int})~(ii)
and (\ref{diff_bc}), we see that $f$ solves (\ref{int_eq_simplified}).
\end{proof}

\bl[Well-defined ODE]
For\label{L:fcdef} each $c\in\half$, there exists a unique
continuous function $f_c:[0,1]\to\R$ that solves (\ref{bifix})~(i) and (ii).
\el

\bpro
Solutions to (\ref{bifix})~(i) and (ii) exist and are unique up to the first
time $\tau$ when $f(r)=\ha r$. Since solutions are nondecreasing with
$f(0)=\ha$, we have $\tau\geq 1$. If $\tau=1$ then $f(1)=\ha$ which
corresponds to the case $f(r)=\ha$ $(r\in[0,1])$, so $f$ is in any case
continuous on $[0,1]$.
\epro

\bl[Integral equation for $f_c$] \label{lemma_int_eq_for_psi} The function $f_c$ from
  Lemma~\ref{L:fcdef} satisfies
\begin{equation}\label{int_eq_for_f_psi}
   f_c(r)^2= \ffrac{1}{4} +r f_c(r) - \int_0^r f_c(u) \, \mathrm{d}u + c r^2, \quad 0 \leq r \leq 1, \quad c \in\half.
 \end{equation}
 \el
\begin{proof} \eqref{int_eq_for_f_psi} holds for $r=0$ since $f_c(0)=\ha$
and the derivatives of the two sides of
\eqref{int_eq_for_f_psi} are equal for all $0\leq r \leq 1$ by
(\ref{diff_eq_for_f_fix}).
\end{proof}

\begin{proof}[Proof of Lemma \ref{L:bifix}]
We note that the function $\ov f(r)=\ha$ $(r\in[0,1])$ solves
(\ref{diff_eq_for_f_fix}) for $r\in[0,1)$ and $c=0$.
In view of this, Lemma~\ref{map_is_homogeneous},
Remark~\ref{remark_constant_f_triv}, and Lemma~\ref{lemma_from_int_to_ode}
show that $T^{(2)}(\rho^{(2)})=\rho^{(2)}$ if and only if the function $f$
from (\ref{rhoF}) satisfies (\ref{bifix})~(i) and (ii) with $c=\frac{1}{4}+\ha
f(1)-\int_0^1f(s)\di s\geq 0$. To see that this latter condition is equivalent
to (\ref{bifix})~(iii), we insert $r=1$ into \eqref{int_eq_for_f_psi} which
yields $\ffrac{1}{4} +\frac{f_c(1)}{2}  -\int_0^1 f_c(s)\,
\mathrm{d}s=f_c(1)^2- \frac{f_c(1)}{2}- c$.
\end{proof}

\subsection{ Finding the nontrivial control parameter }
\label{subsection_proof_of_lemma_fixsol}

The goal of this subsection is to prove Lemma~\ref{L:fixsol}.
By Lemma~\ref{L:fcdef}, the ODE (\ref{bifix})~(i) with the left boundary
condition (\ref{bifix})~(ii) has a unique solution $f_c$ for all $c\geq 0$.
We need to prove the existence and uniqueness of a control parameter
$c_2>0$ for which $f_{c_2}$ also solves the right boundary
condition (\ref{bifix})~(iii). In Lemma~\ref{lemma_f_1_implicit_eq}
we solve the ODE and obtain an implicit equation for $f_c(1)$.
In Lemma \ref{lemma_h} we use this to rewrite (\ref{bifix})~(iii) as $h(c)=1$
for some explicit function $h$ (see \eqref{h_psi_def}). In Lemma
\ref{lemma_unique_fixed_point} we show that there is a unique
$c_2>0$ such that $h(c_2)=1$ holds, and $c_2\in(0,\frac{1}{4})$.

Given any $c \in (0,\infty)$, let us define  $g_+(c), g_-(c), A_+(c), A_-(c)$ by
\begin{equation}\label{g_pm_A_pm}
  g_{\pm}(c):= \frac{1}{4}\left( 1 \pm \sqrt{1+16c} \right), \qquad A_\pm(c):= -\frac{1}{2} \pm \frac{1}{2 \sqrt{1+ 16 c}}.
\end{equation}


\bl[Solution of ODE for $f_c$] \label{lemma_f_1_implicit_eq}
For any $c>0$, the function $f_c$ from Lemma~\ref{L:fcdef} is given by
$f_c(r)=rg_c(r)$ $(r\in(0,1])$, where $g_c(r)$ is the unique element of $(g_+(c), +\infty)$ that satisfies
\begin{equation}\label{g_psi_1_implicit}
  \ffrac{1}{2} \left( g_c(r)-g_+(c) \right)^{A_+(c)} \left( g_c(r)-g_-(c) \right)^{A_-(c)}=r.
\end{equation}
\el
\begin{proof} If we define
$ g_c(r):= f_c(r)/r$ for any $r \in (0,1]$,
 then we can use (\ref{bifix})~(i) to show that the function $r \mapsto g_c(r)$ solves the
 ODE
 \begin{equation}\label{g_ode}
   \frac{g_c(r)-1/2}{c-g_c(r)(g_c(r)-1/2)  }  g_c'(r) = \frac{1}{r}, \qquad r \in (0,1].
\end{equation}
 We first find the general solution of this ODE by integrating both sides of \eqref{g_ode}.
 In order to calculate the indefinite integral of the l.h.s., we perform the substitution $g=g_c(r)$ and apply the partial fraction decomposition
 \begin{equation}
  \frac{1/2-g}{g^2-g/2-c}\stackrel{(\ref{g_pm_A_pm}) }{=}\frac{A_+(c)}{g-g_+(c)}+ \frac{A_-(c)}{g-g_-(c)}.
\end{equation}
Integrating and then exponentiating both sides of \eqref{g_ode}, we obtain that
 the general solution of \eqref{g_ode} satisfies the implicit equation $R( g_c(r) )=r$ for any  $r \in (0,1]$, where
\begin{equation}\label{g_implicit_general}
 R(g):=
  \alpha^* \left( g -g_+(c) \right)^{A_+(c)} \left( g - g_-(c) \right)^{A_-(c)}, \quad g \in ( g_+(c),+\infty)
\end{equation}
for some positive constant $\alpha^*$. Note that the function $g \mapsto R(g)$
is strictly decreasing (since both $A_+(c)$ and $A_-(c)$ are negative) and
that it satisfies $\lim_{g \to g_+(c)}R( g ) = +\infty$ as well as $\lim_{g \to
  \infty}R( g ) =0$.  Therefore, the equation $R( g )=r$ has a unique solution
$g$ for any $r \in (0,1]$.  In order to identify the value of $\alpha^*$, we
  observe that (\ref{bifix})~(ii) is equivalent to $\lim_{r \to 0_+} g_c(r) r
  =\ffrac{1}{2}$, which is in turn equivalent to
\begin{equation}\label{g_asymp_to_infty}
  \lim_{g \to \infty} g \alpha^* \left( g-g_+(c) \right)^{A_+(c)} \left( g-g_-(c) \right)^{A_-(c)} =\ffrac{1}{2}.
\end{equation}
Noting that $A_+(c) + A_-(c)=-1$ (c.f.\ \eqref{g_pm_A_pm}), we obtain
$\alpha^*=\ffrac{1}{2}$ using \eqref{g_asymp_to_infty}.
\end{proof}

\bl[$f_c$ is increasing and concave]
For\label{L:fcprop} any $c>0$, the function $f_c$ from Lemma~\ref{L:fcdef} is
twice continuously differentiable with $f_c(0)=\ha$, $f'_c(r)\geq 0$, and
$f''_c(r)>0$ $\big(r\in(0,1]\big)$.
\el
\bpro
The facts that $f_c(0)=\ha$ and $f'_c(r)\geq 0$ are immediate from
(\ref{bifix})~(i) and (ii). To see that $f_c$ is twice continuously
differentiable with $f''_c(r)\geq 0$, we observe that by (\ref{bifix})~(i),
\be\label{fcc}
f''_c(r)=\dif{r}\frac{c}{r^{-1}f_c(r)-\ha}=\dif{r}\frac{c}{g_c(r)-\ha},
\ee
where $g_c$ is the function in Lemma~\ref{lemma_f_1_implicit_eq}.
Since the function in (\ref{g_implicit_general}) is strictly decreasing,
$g_c(r)$ is strictly decreasing, and hence the right-hand side of (\ref{fcc})
is strictly positive for $r>0$.
\epro

Let us define
\begin{equation}\label{h_psi_def}
 h(c):= \ffrac{1}{4} \left(\frac{\sqrt{1+32 c }-\sqrt{1+16 c}  }{ \sqrt{1+32 c }+\sqrt{1+16 c}   }\right)^{\frac{1}{\sqrt{1+16 c}}}\frac{1}{c}, \quad c \in (0, \infty).
 \end{equation}

\bl[Right boundary condition]
Let\label{lemma_h} $c \in (0,+\infty)$. The following conditions are equivalent:
\begin{align}
\label{f_one_quad_cond_for_psi}
  2c &= f_c(1)^2- \ha f_c(1) \\
  \label{solution_of_quadratic_eq}
  f_c(1) &= \ffrac{1}{4}\left(1 +\sqrt{ 1 +32 c  }\right)\\
\label{h_is_one_cond_for_psi}
  h(c)&=1
\end{align}
\el
\begin{proof} The positive solution of the quadratic equation
\eqref{f_one_quad_cond_for_psi} is \eqref{solution_of_quadratic_eq}.
By Lemma~\ref{lemma_f_1_implicit_eq}, $f_c(1)$ is the unique element of $(g_+(c), +\infty)$ that satisfies
\begin{equation}\label{f_psi_1_implicit}
  \ffrac{1}{2} \left( f_c(1)-g_+(c) \right)^{A_+(c)} \left( f_c(1)-g_-(c) \right)^{A_-(c)}=1.
\end{equation}
 Now by \eqref{f_psi_1_implicit} and $A_+(c) + A_-(c)=-1$,  \eqref{solution_of_quadratic_eq} is equivalent to
\begin{equation}\label{pure_psi_compl_identity}
   2\left(   \sqrt{1+32c} -\sqrt{1+16c}   \right)^{A_+(c)} \left(  \sqrt{1+32c} +\sqrt{1+16c}   \right)^{A_-(c)}=1.
\end{equation}
Finally, the equivalence of the condition \eqref{pure_psi_compl_identity} and \eqref{h_is_one_cond_for_psi} (c.f.\ \eqref{h_psi_def}) follows using elementary algebra.
\end{proof}

\bl[Existence and uniqueness of the positive root]
There\label{lemma_unique_fixed_point} exists exactly one $c_2 \in (0, +\infty)$ such that $h(c_2)=1$. Moreover we have
 \begin{equation}
 \label{psi_less_than_quarter}
  c_2 \in (0, \ffrac{1}{4}).
 \end{equation}
 \el
\begin{proof} Let us first observe that $\lim_{c \to 0_+} h(c)=1$, thus $h$ is a continuous function on $[0,+\infty)$ if we define $h(0):=1$.
Next we observe that
\begin{equation}\label{h_quarter}
  h(1/4) \stackrel{(\ref{h_psi_def}) }{=}  \left(\frac{3-\sqrt{5}  }{ 3+\sqrt{5}   }\right)^{\frac{1}{\sqrt{5}}} <  1.
\end{equation}
We will show that
\begin{equation}\label{unique_maximum_of_h}
  \exists \, \widetilde{c} \in (0,+\infty) \; : \; h'(c) >0 \text{ if } c \in (0,\widetilde{c}), \text{ but } h'(c) <0 \text{ if } c \in (\widetilde{c},+\infty).
\end{equation}
Once we have this, the statement of Lemma~\ref{lemma_unique_fixed_point} will
follow from the facts that $h(0)=1$ and $h(1/4)<1$.

It remains to prove \eqref{unique_maximum_of_h}. Let us define
\begin{equation}\label{k_psi_r_psi_def}
 k(c):=\ln(h(c/16)), \qquad r(c):=  (1+c)^{3/2} k'(c).
\end{equation}
Let us observe that in order to prove \eqref{unique_maximum_of_h}, it is
enough to prove
 \begin{equation}\label{unique_root_of_r}
  \exists \, \widehat{c} \in (0,+\infty) \; : \; r(c) >0 \text{ if } c \in (0,\widehat{c}), \text{ but } r(c) <0 \text{ if } c \in (\widehat{c},+\infty),
\end{equation}
where actually $\widehat{c}=16 \widetilde{c}$.
It remains to prove \eqref{unique_root_of_r}. First note that  we have
\begin{equation}\label{k_prime}
  k'(c)= \frac{-c +(1+2c)^{-1/2} -1 }{c^2+c } - \frac12 (1+c)^{-3/2} \ln\left(  \frac{\sqrt{1+2c}-\sqrt{1+c}  }{\sqrt{1+2c}+\sqrt{1+c} } \right),
\end{equation}
thus $\lim_{c \to 0_+}k'(c)=+\infty$. Also note that $ \exists \, c \, : \, k'(c)<0$, since $k(0)=0$ and $k(4)<0$ by \eqref{h_quarter} and \eqref{k_psi_r_psi_def}. These observations imply that the function
 $c \mapsto r(c)$ takes both positive and negative values. Thus in order to prove \eqref{unique_root_of_r}, it is enough to prove that $r: (0, +\infty) \to \mathbb{R}$ is a decreasing function.
 \begin{equation}
   r'(c)= \frac{ \sqrt{1+c} }{2 c^2 (2 c+1)^{3/2} }q(c), \quad \text{where} \quad   \quad q(c):=  \sqrt{2 c +1}(2-2c^2+3c)-2-6c.
 \end{equation}
 It remains to check that $q(c)<0$ for all $c>0$. This readily follows after we observe that
 $q(0)=0$, $q'(0)=-1$ and $q''(c)=\frac{-15 c}{\sqrt{2 c +1}}$ for any $c>0$.
\end{proof}

\begin{remark}
Although it is just elementary calculus, the proof of the uniqueness part
of Lemma~\ref{lemma_unique_fixed_point} is one of the trickiest of the paper.
Since ultimately, the uniqueness of the nontrivial scale invariant fixed point
of Theorem~\ref{T:scalefix} hinges on this, one would like to find a more
elegant and insightful proof. It is tempting to try and prove that the
function $h$, or the function $c\mapsto ch(c)$, are either convex or concave
on the entire positive axis, but this is \emph{not true}. The function
$c\mapsto f_c(1)^2-\ha f_c(1)$ that occurs in (\ref{bifix})~(iii) appears to be
concave, but we have been unable to prove so.
\end{remark}


\subsection{ Non-trivial solution of the bivariate RDE }
\label{subsection_p_of_l_nontriv}

\begin{proof}[Proof of Lemma \ref{L:nontriv}]
We apply Lemma~\ref{L:suffPast} to the function $f_{c_2}$.
Condition~(i) is satisfied since
\be
f_{c_2}(1)
\stackrel{(\ref{solution_of_quadratic_eq})}{=}
\ffrac{1}{4}\left(1 +\sqrt{ 1 +32 c_2  }\right)
\stackrel{(\ref{psi_less_than_quarter}) }{<}
\ffrac{1}{4}\left(1 +\sqrt{9}\right)=1.
\ee
Conditions~(ii), (iii) and (v) of Lemma~\ref{L:suffPast} are satisfied by
Lemma~\ref{L:fcprop}, so it remains to check condition~(iv), which requires
$2f'_{c_2}(1)=f_{c_2}(1)$. Using (\ref{diff_eq_for_f_fix})~(ii), we can rewrite this as
\be\label{nodiag}
\frac{2c_2}{f_{c_2}(1)-\ha}=f_{c_2}(1),
\ee
which is satisfied by Lemmas \ref{lemma_h} and  \ref{lemma_unique_fixed_point}.
\end{proof}

\begin{remark}
Formula (\ref{nodiag}) shows that condition (\ref{bifix})~(iii) is equivalent
to the statement that the measure $\rho^{(2)}$ associated with $f$ puts no
mass on the diagonal $\big\{(r,r):r\in[0,1]\big\}$.
\end{remark}

\section{Frozen percolation}\label{S:frz}

\subsection{Outline}

In the previous section, we have proved Theorem~\ref{T:scalefix}, which
implies that the RTP $(\tau_\ibf,\kappa_\ibf,Y_\ibf)_{\ibf\in\T}$
corresponding to the map $\chi$ from (\ref{chidef}) and law $\rro$ from
(\ref{rhodef}) is nonendogenous.
In the present section, we provide the proofs of our
remaining results, which are Theorems~\ref{T:orfrz}, \ref{T:main},
\ref{T:dirfrz}, and \ref{T:nonend}, as well as Lemma~\ref{L:perc},
Proposition~\ref{P:scale}, and Lemma~\ref{L:scinv}.

In Subsection~\ref{S:RDEequiv}, we show that there is a one-to-one
correspondence between solutions to the RDEs (\ref{frzRDE}) and (\ref{MRDE}),
under which the measure $\nu$ from (\ref{nudef}) corresponds to the measure
$\rho$ from (\ref{rhodef}). We also prove a correspondence between solutions
to the associated bivariate RDEs and use this to derive Theorem~\ref{T:nonend}
from Theorem~\ref{T:scalefix}.

In Subsection~\ref{S:RDEsol}, we classify all solutions to the RDE
(\ref{MRDE}). Using results from the preceding subsection, this also leads to
a description of general solutions to the RDE (\ref{frzRDE}).

Theorem~\ref{T:dirfrz} is proved in Subsections~\ref{S:MBBTfrz} and
\ref{S:binfrz}. In Subsection~\ref{S:MBBTfrz}, we use the classification of
solutions to (\ref{MRDE}) to prove a version of Theorem~\ref{T:dirfrz} for
frozen percolation on the MBBT. In Subsection~\ref{S:binfrz} this is then
translated into a result for the oriented binary tree using a coupling between
two RTPs, one for frozen percolation on the MBBT, and the other for the
oriented binary tree.

In Subsection~\ref{S:frzMBBT}, we prove Lemma~\ref{L:perc} as well as
Proposition~\ref{P:scale} and Lemma \ref{L:scinv} about scale invariance of
the MBBT. Lemma~\ref{L:scinv} allows us to identify the nontrivial solution
$\rho^{(2)}_2$ of the bivariate RDE from Theorem~\ref{T:scalefix} as
$\un\rro^{(2)}$. Using results from Subsection~\ref{S:RDEsol}, we use this
to obtain an explicit formula for $\un\nnu^{(2)}$ based on our formula for
$\un\rro^{(2)}$.

In Subsection~\ref{S:unorient} we mainly rely on arguments from \cite{Ald00}
to translate results about frozen percolation on the oriented binary tree into
results about frozen percolation on the 3-regular tree. In particular, we
derive Theorem~\ref{T:orfrz} from Theorem~\ref{T:dirfrz} and
Theorem~\ref{T:main} from Theorem~\ref{T:nonend}.

\subsection{Equivalence of RDEs}\label{S:RDEequiv}

In this subsection, we show that there is a one-to-one correspondence between
solutions to the RDEs (\ref{frzRDE}) and (\ref{MRDE}), under which the measure
$\nu$ from (\ref{nudef}) corresponds to the measure $\rho$ from
(\ref{rhodef}). We also prove a correspondence between solutions to the
associated bivariate RDEs and use this to derive Theorem~\ref{T:nonend}
from Theorem~\ref{T:scalefix}. We start with a simple observation.

\bl[No burning before the critical point]
Every\label{L:half} solution $\ngu$ to the RDE (\ref{frzRDE}) is
concentrated on $I':=[\ha,1]\cup\{\infty\}$.
\el

\bpro
If $\ngu$ solves the RDE (\ref{frzRDE}), then we can construct an RTP
$(\tau_\ibf,X_\ibf)_{\ibf\in\T}$ corresponding to the map $\ga$ from
(\ref{gadef}) and $\ngu$. Then by (\ref{XupX}),
\be
\ngu\big([0,\ha]\big)=\P[X_\wurz\leq \ha]\leq\P[X^\up_\wurz\leq \ha]
=\P\big[\wurz\percol{\T^{1/2}\beh\Fx}\infty\big]
\leq\P\big[\wurz\percol{\T^{1/2}}\infty\big]=0,
\ee
where the last equality follows from the fact that a branching process with a
binomial offspring distribution with parameters $2,\ha$ is critical and hence
dies out a.s.
\epro

The next lemma, which is the first main result of the present subsection, says
that there is a one-to-one correspondence between solutions to the RDEs
(\ref{frzRDE}) and (\ref{MRDE}). The idea behind the proof (and in
  particular the occurrence of the geometric distribution in (\ref{Tgeom}))
  will become more clear in Section~\ref{S:binfrz} below.

\bl[Equivalence of RDEs]
Let\label{L:equivRDE} $I':=[\ha,1]\cup\{\infty\}$ and let $\Hh:I\to I'$ be the
bijection defined by $\Hh(t):=1/(2-t)$ $(t\in[0,1])$ and
$\Hh(\infty):=\infty$. If $\rgo$ solves the RDE (\ref{MRDE}), then its image
under the map $\Hh$ solves the RDE (\ref{frzRDE}). Conversely, if $\nggu$
solves the RDE (\ref{frzRDE}), then its image under the map $\Hh^{-1}$ solves
the RDE (\ref{MRDE}).
\el

\bpro
Let $T_{\rm y}$ be defined as in (\ref{Tdef}) but for the map $\chi$ in
(\ref{chidef}), i.e.,
\be
T_{\rm y}(\mu):=\mbox{ the law of }\chi[\tau_\wurz,\kappa_\wurz](Y_1,Y_2)
\ee
where $Y_1,Y_2$ are i.i.d.\ with law $\mu$ and independent of
$(\tau_\wurz,\kappa_\wurz)$. Then we can write
\be
T_{\rm y}=\ha T_\Phi+\ha T_{\tt min},
\ee
where
\be
T_\Phi(\mu):=\mbox{ the law of }\Phi[\tau_\wurz](Y_1)
\quand
T_{\tt min}(\mu):=\mbox{ the law of }Y_1\wedge Y_2,
\ee
and $\Phi:[0,1]\times I\to I$ denotes the function
\be\label{Phidef}
\Phi[t](x):=\left\{\ba{ll}
x\quad&\mbox{if }x>t,\\[5pt]
\infty\quad&\mbox{if }x\leq t.\ea\right.
\ee
Note that the map $T_\Phi$ is linear, but $T_{\tt min}$ is not.
Let us define
\be\label{Tgeom}
T_{\rm z}:=\sum_{n=1}^\infty 2^{-n}T_\Phi^{n-1}T_{\tt min}.
\ee
We claim that $\mu$ is a fixed point of $T_{\rm y}$ if and only if it is a
fixed point of $T_{\rm z}$. Indeed, $T_{\rm y}(\mu)=\mu$ implies $T_{\rm
  min}(\mu)=2\mu-T_\Phi(\mu)$ and hence, using the linearity of $T_\Phi$,
\be
T_{\rm z}(\mu)=\sum_{n=1}^\infty 2^{-n}T_\Phi^{n-1}\big(2\mu-T_\Phi(\mu)\big)=\mu.
\ee
Conversely, since
\be
T_{\rm z}=\ha T_{\tt min}+\ha T_\Phi\circ T_{\rm z}
\ee
$T_{\rm z}(\mu)=\mu$ implies $\mu=\ha T_{\tt min}(\mu)+\ha T_\Phi(\mu)=T_{\rm
  y}(\mu)$.

We observe that
\be
T_{\rm z}(\mu):=\mbox{ the law of }
\Phi[\tau_1]\circ\cdots\circ\Phi[\tau_N](Y_1\wedge Y_2)
=\ga[\tau_1\vee\cdots\vee\tau_N](Y_1,Y_2),
\ee
where $(\tau_k)_{k\geq 1}$ are uniformly distributed on $[0,1]$, the r.v.'s
$Y_1,Y_2$ have law $\mu$, the r.v.\ $N$ is geometrically distributed with
$\P[N=n]=2^{-n-1}$ $(n\geq 0)$, and all r.v.'s are independent. Since
\be\label{maxgeom}
\P\big[\tau_1\vee\cdots\vee\tau_N\leq t\big]
=\sum_{n=0}^\infty 2^{-n-1}t^n=\frac{\ha}{1-\ha t}=\Hh(t)
\qquad\big(t\in[0,1]\big),
\ee
we have that $\tau:=\Hh(\tau_1\vee\cdots\vee\tau_N)$ satisfies
$\P[\tau=\ha]=\P[N=0]=\ha$ and
\be
\P\big[\tau<t\big]=
\P\big[\tau_1\vee\cdots\vee\tau_N<\Hh^{-1}(t)\big]=\Hh\big(\Hh^{-1}(t)\big)
=t\qquad\big(t\in[\ha,1]\big).
\ee
Then, using the fact that
\be\label{moninv}
\ga[\Hh(t)]\big(\Hh(x),\Hh(y)\big)=\Hh\big(\ga[t](x,y)\big)
\qquad\big(x,y\in I,\ t\in[0,1]\big)
\ee
and using also Lemma~\ref{L:half}, we see that the law $\mu$ of an $I$-valued
random variable $Y$ solves the RDE $T_{\rm z}(\mu)=\mu$ or equivalently
\be
Y\isd\ga[\tau_1\vee\cdots\vee\tau_N](Y_1,Y_2),
\ee
if and only if $X:=\Hh(Y)$, $X_1:=\Hh(Y_1)$, and $X_2:=\Hh(Y_2)$ solve the RDE
(\ref{frzRDE}).
\epro

\bl[Equivalence of special solutions]
The\label{L:rhonu} measure $\nnu$ in (\ref{nudef}) is the image of the measure
$\rro$ in (\ref{rhodef}) under the map $\Hh$.
\el

\bpro
Since $\Hh^{-1}(t)=2-1/t$ $(t\in[\ha,1])$ is the inverse of $\Hh(t):=1/(2-t)$
$(t\in[0,1])$, we see that
\be
\rro\big([0,\Hh^{-1}(t)]\big)=\ha \Hh^{-1}(t)=1-\frac{1}{2t}=\nnu\big([0,t]\big)
\qquad\big(t\in[\ha,1]\big),
\ee
which shows that $\nnu$ is the image of $\rro$ under $\Hh$.
\epro

We next turn our attention to the bivariate RDEs.

\bl[Equivalence of bivariate RDEs]
Let\label{L:bivequivRDE} $H:I\to I'$ be the map defined in
Lemma~\ref{L:equivRDE}. Let $T^{(2)}_{\rm x}$ and $T^{(2)}_{\rm y}$ be the
bivariate maps defined as in (\ref{bivdef}) for the maps $\ga$ in
(\ref{gadef}) and $\chi$ in (\ref{chidef}), respectively. Then a measure
$\mu^{(2)}\in\Pc(I^2)$ solves the bivariate RDE $T^{(2)}_{\rm
  y}(\mu^{(2)})=\mu^{(2)}$ if and only if its image $\nu^{(2)}$ under the map
$(y_1,y_2)\mapsto\big(H(y_1),H(y_2)\big)$ solves the bivariate RDE
$T^{(2)}_{\rm x}(\nu^{(2)})=\nu^{(2)}$.
\el

\bpro
Let $T_H:\Pc(I)\to\Pc(I')$ be the function that maps a measure on $I$ to its
image under the map $H$. The proof of Lemma~\ref{L:equivRDE} consisted of
showing that for any $\mu\in\Pc(I)$, one has $T_{\rm y}(\mu)=\mu$ if and only
if $T_{\rm z}(\mu)=\mu$, and moreover $T_{\rm x}T_H=T_HT_{\rm z}$. With
exactly the same proof, these statements remain true if we replace the maps
$T_{\rm x},T_{\rm y},T_{\rm z}$, and $T_H$ with their bivariate versions
$T^{(2)}_{\rm x},T^{(2)}_{\rm y},T^{(2)}_{\rm z}$, and $T^{(2)}_H$. It follows
that $\mu^{(2)}\in\Pc(I^2)$ solves $T^{(2)}_{\rm y}(\mu^{(2)})=\mu^{(2)}$ if
and only if $\nu^{(2)}:=T^{(2)}_H(\mu^{(2)})$ solves $T^{(2)}_{\rm
  x}(\nu^{(2)})=\nu^{(2)}$, which is the claim of the lemma.
\epro

Our results so far allow us to prove Theorem~\ref{T:nonend}.\med

\bpro[Proof of Theorem~\ref{T:nonend}]
By Theorem~\ref{T:scalefix}, the bivariate map $T^{(2)}_{\rm y}$ has a fixed
point $\rho^{(2)}_2\in\Pc(I^2)_\rho$ that is not concentrated on the diagonal
$\{(y,y):y\in I\}$. Let $\nu^{(2)}_2$ denote the image of $\rho^{(2)}_2$ under
the map $(y_1,y_2)\mapsto\big(H(y_1),H(y_2)\big)$. Then
$\nu^{(2)}_2\in\Pc(I^2)_\nu$ by Lemma~\ref{L:rhonu}. By
Lemma~\ref{L:bivequivRDE}, $\nu^{(2)}_2$ is a fixed point of
$T^{(2)}_{\rm x}$. Since $\nu^{(2)}_2$ is not concentrated on the diagonal,
Theorem~\ref{T:bivar}~(i) and (iii) imply that the RTP associated with $\nu$
is nonendogenous.
\epro

Each solution $\mu$ to an RDE defines an RTP, which through (\ref{unnu})
defines a special solution $\un\mu^{(2)}$ to the corresponding bivariate
RDE. In particular, we define $\un\nnu^{(2)}$ and $\un\rro^{(2)}$ in this way
starting from the measures $\rro$ and $\nnu$ defined in (\ref{rhodef}) and
(\ref{nudef}). The final result of this subsection relates these measures to
each other.

\bl[Nontrivial solutions to bivariate RDE]
Let\label{L:nurho} $(Y_\wurz,Y'_\wurz)$ be a random variable with law
$\un\rro^{(2)}$ and let $\Hh$ be the function from (\ref{sigF}).
Then $\big(\Hh(Y_\wurz),\Hh(Y'_\wurz)\big)$ has law
$\un\nnu^{(2)}$.
\el

\bpro
We will use a characterization of $\un\rro^{(2)}$ and $\un\nnu^{(2)}$ from
\cite{MSS18}. We first need some abstract definitions. Let $I$ be a
Polish space. If $\xi$ is a random probability law on $I$, and
$\eta\in\Pc(\Pc(I))$ is the law of $\xi$, then
\be
\eta^{(n)}:=\E\big[\underbrace{\xi\otimes\cdots\otimes\xi}_{\mbox{$n$ times}}\big]
\ee
is called the \emph{$n$-th moment measure} of $\eta$. In \cite[Lemma~2]{MSS18}
it was shown that for each map $T$ of the form (\ref{Tdef}), there exists a
\emph{higher level map} $\ch T:\Pc(\Pc(I))\to\Pc(\Pc(I))$ that is uniquely
characterised by
\be\label{chT}
\ch T(\eta)^{(n)}=T^{(n)}(\eta^{(n)})
\qquad\big(n\geq 1,\ \eta\in\Pc(\Pc(I))\big),
\ee
where $T^{(n)}$ is the associated \emph{$n$-variate map}.
Let $\nu$ be a solution to the RDE (\ref{RDE}) and let $\Pc(\Pc(I))_\nu$
denote the space of all $\eta\in\Pc(\Pc(I))$ with $\eta^{(1)}=\nu$.
In \cite[Prop~3]{MSS18}, it was shown that the set
$\{\eta\in\Pc(\Pc(I))_\nu:\ch T(\eta)=\eta\}$, equipped with the
\emph{convex order}, has a unique minimal element $\un\nu$ and maximal element
$\ov\nu$. Moreover, by \cite[Lemma~2 and Props~3 and 4]{MSS18}, the measures
$\un\nu^{(2)}$ and $\ov\nu^{(2)}$ from (\ref{ovnu}) and (\ref{unnu}) are the
second moment measures of $\un\nu$ and $\ov\nu$.

We now return to our special setting with $I=[0,1]\cup\{\infty\}$.  Let
$T_{\rm x},T_{\rm y}$, and $T_H$ be as in the proof of
Lemma~\ref{L:bivequivRDE}. In Lemma~\ref{L:equivRDE}, we have proved that
$\mu\in\Pc(I)$ satisfies $T_{\rm y}(\mu)=\mu$ if and only if $\nu:=T_H(\mu)$
satisfies $T_{\rm x}(\nu)=\nu$. In Lemma~\ref{L:bivequivRDE}, we have shown
that the same is true for the bivariate maps $T^{(2)}_{\rm x},T^{(2)}_{\rm
  y}$, and $T^{(2)}_H$. The argument carries over without a change for
general $n$-variate maps and therefore, by (\ref{chT}), the statement is also
true for the associated higher-level maps $\ch T_{\rm x},\ch T_{\rm y}$, and
$\ch T_H$. In particular, using also Lemma~\ref{L:rhonu}, we obtain that the
image of the set
\be
A:=\big\{\eta\in\Pc(\Pc(I))_\rho:\ch T_{\rm y}(\eta)=\eta\big\}
\ee
under the higher-level map $\ch T_H$ is the set
\be
B:=\big\{\eta\in\Pc(\Pc(I')))_\nu:\ch T_{\rm x}(\eta)=\eta\big\}.
\ee
Since by \cite[Prop~3]{MSS18}, higher-level maps are monotone w.r.t.\ the
convex order, $\ch T_H$ maps the minimal element of $A$, which is $\un\rho$,
into the minimal element of $B$, which is $\un\nu$. By (\ref{chT}), this
implies that the bivariate map $T_H^{(2)}$ maps $\un\rho^{(2)}$ to
$\un\nu^{(2)}$, which is the claim we wanted to prove.
\epro

\subsection{General solution of the RDE}\label{S:RDEsol}

In this subsection, we classify all solutions to the RDE (\ref{MRDE}). Through
Lemma~\ref{L:equivRDE}, this then also implies the form of a general solution
of the RDE (\ref{frzRDE}), significantly extending \cite[Lemma~3]{Ald00}, who
only considered solutions without atoms in $[0,1]$.

Let $O\sub(0,1]$ be open. Then $O$ is a countable union of disjoint open
intervals $(O_k)_{0\leq k<n+1}$ for some $0\leq n\leq\infty$ (with
$\infty+1:=\infty$). Without loss
of generality we can assume that $\emptyset\neq O_k\sub(0,1)$ for all
$1\leq k<n+1$ while either $1\in O_0$ or $O_0=\emptyset$.
We let $x_k\in(0,1)$ and $c_k>0$ denote the center and radius
of $O_k$, respectively, i.e., $O_k=(x_k-c_k,x_k+c_k)$, and we choose
$x_0\in(0,1]\cup\{2\}$ and $c_0>0$ such that
$O_0=(x_0-c_0,x_0+c_0)\cap(0,1]$. We define a measure $\rgo$ on $[0,1]$ by
\be\label{gensol}
\rgo(\di t):=\ha 1_{[0,1]\beh O}(t)\di t+1_{\{x_0\neq 2\}}c_0\de_{x_0}(\di t)
+\sum_{k=1}^n c_k\de_{x_k}(\di t).
\ee
It is easy to see that $\rgo([0,1])\leq 1$, so we can unambiguously extend
$\rgo$ to a probability measure on $I=[0,1]\cup\{\infty\}$. We will prove the
following result.

\bp[General solution to RDE]
The\label{P:genRDE} probability measure $\rgo$ defined in (\ref{gensol})
solves the RDE (\ref{MRDE}), and conversely, every solution of (\ref{MRDE}) is
of this form.
\ep

We need one preparatory lemma.

\bl[RDE for MBBT]
A\label{L:rhoRDE} probability measure $\rgo$ on $I$ solves the RDE
(\ref{MRDE}) if and only if
\be\label{rhoRDE}
\int_{[0,t]}\!\rgo(\di s)\,s=\rgo\big([0,t]\big)^2\qquad\big(t\in[0,1]\big).
\ee
\el

\bpro
Let $\Phi$ be the function defined in (\ref{Phidef}). Then
\be
\chi[\tau,1](x,y)=\Phi[\tau](x)
\quand
\chi[\tau,2](x,y)=x\wedge y.
\ee
Using this and the fact that the function $F(t):=\rgo\big([0,t]\big)$
$(t\in[0,1])$ uniquely characterizes $\rgo$, we see that (\ref{MRDE}) is
equivalent to
\bc
\dis F(t)=\P[Y_\wurz\leq t]
&=&\dis\ha\int_0^1\!\di s\,\P\big[\Phi[s](Y_1)\leq t\big]
+\ha\P[Y_1\wedge Y_2\leq t]\\[5pt]
&=&\dis\ha\int_0^1\!\di s\,\P[s<Y_1\leq t]
+\ha\big(1-\P[Y_1>t]^2\big)\\[5pt]
&=&\dis\ha\int_0^t\!\di s\,\big\{F(t)-F(s)\big\}
+\ha\big(1-(1-F(t))^2\big)\\[5pt]
&=&\dis\ha tF(t)-\ha\int_0^t\!\di s\,F(s)
+F(t)-\ha F(t)^2,
\ec
which can be rewritten as
\be\label{FMRDE}
\big(t-F(t)\big)F(t)=\int_0^t\!\di s\,F(s)\qquad\big(t\in[0,1]\big).
\ee
Using the fact that
\be\label{dF}
tF(t)=\int_{[0,t]}\di(sF(s))=\int_{[0,t]}s\,\di F(s)+\int_{[0,t]}F(s)\,\di s,
\ee
we can rewrite (\ref{FMRDE}) as (\ref{rhoRDE}).
\epro

\bpro[Proof of Proposition~\ref{P:genRDE}]
We first prove that the measure in (\ref{gensol}) solves (\ref{rhoRDE}).
With $x_k$ and $c_k$ as in (\ref{gensol}), we will prove that the measure
$\rgo'$ on $\half$ defined by
\be
\rgo'(\di t):=\ha 1_{\half\beh O}(t)\di t+\sum_{k=0}^n c_k\de_{x_k}(\di t)
\ee
solves (\ref{rhoRDE}) for all $t\geq 0$. Restricting $\rgo'$ to $[0,1]$ we
then see that $\rgo$ satisfies (\ref{rhoRDE}) for all $t\in[0,1]$.

If $\rgo'(\di s):=\ha\di s$ then the left-hand side of (\ref{rhoRDE}) is
$\ha\cdot\ha t^2$ while the right-hand side is $(\ha t)^2$, so (\ref{rhoRDE})
holds. Next, if we modify $\rgo'$ by concentrating all the mass in an interval
of the form $(x-c,x+c)$ in the middle of that interval, then (\ref{rhoRDE})
remains true for all $t\leq x-c$ and $t\geq x+c$. Applying this observation
inductively and taking the limit, we see that $\rgo'$
solves (\ref{rhoRDE}) for all $t\in\half\beh O$. But the left- and right-hand
sides of (\ref{rhoRDE}) are constant on the intervals $[x_k-c_k,x_k)$ and
$[x_k,x_k+c_k]$ $(k\geq 0)$ so (\ref{rhoRDE}) holds for all $t\geq 0$.

The proof that all solutions of (\ref{rhoRDE}) are of the form (\ref{gensol})
goes in a number of steps. Taking increasing limits, we observe that
(\ref{rhoRDE}) implies
\be\label{lrhoRDE}
\int_{[0,t)}\!\rgo(\di s)\,s=\rgo\big([0,t)\big)^2\qquad\big(t\in(0,1]\big).
\ee
We next claim that:
\be\label{rhoclaim1}
\mbox{If $\rgo$ solves (\ref{rhoRDE}) and $\rgo\big([0,t)\big)=\ha u$ with
$0\leq t\leq u$, then $\rgo\big([t,u]\big)=0$.}
\ee
Indeed, we obtain from (\ref{rhoRDE}) that
\be
\int_{[0,t)}\!\rgo(\di s)\,s+\int_{[t,u]}\!\rgo(\di s)\,s
=\big[\rgo\big([0,t)\big)+\rgo\big([t,u]\big)\big]^2,
\ee
which using (\ref{lrhoRDE}) and our assumption that $\rgo\big([0,t)\big)=\ha u$
yields
\be
\rgo\big([t,u]\big)^2=\int_{[t,u]}\!\rgo(\di s)\,s-u\rgo\big([t,u]\big)\leq 0.
\ee
Our next claim is that:
\be\ba{l}\label{rhoclaim2}
\dis\mbox{If $\rgo$ solves (\ref{rhoRDE}) and $c:=\rgo(\{t\})>0$ for some
  $t\in[0,1]$,}\\
\dis\mbox{then $c=2\big[\ha t-\rgo\big([0,t)\big)\big]$.}
\ec
Indeed, (\ref{rhoRDE}) implies
\be
\int_{[0,t)}\!\rgo(\di s)\,s+ct=\big[\rgo\big([0,t)\big)+c\big]^2,
\ee
which using (\ref{lrhoRDE}) implies
\be
ct=2c\rgo\big([0,t)\big)+c^2.
\ee
Using our assumption that $c>0$, we arrive at (\ref{rhoclaim2}).
Let $F$ denote the function $F(t):=\rgo\big([0,t]\big)$ $(t\in[0,1])$.
We need one more claim, which says that:
\be\ba{l}\label{rhoclaim3}
\mbox{If $\rgo$ solves (\ref{rhoRDE}) and has no atoms in $[s,u)$,}\\
\mbox{then $\rgo\big([0,s)\big)<\ha s$ implies $\rgo\big([s,u)\big)=0$.}
\ec
Indeed, if $\rgo$ has no atoms in $[s,u)$, then 
the function $F(t):=\rgo\big([0,t]\big)$ $(t\in[0,1])$ solves
\be
t\rgo(\di t)=t\di F(t)\stackrel{(\ref{rhoRDE})}{=}\di\big(F(t)^2\big)
=2F(t)\di F(t)=2F(t)\rgo(\di t)
\ee
on $[s,u)$, which shows that the restriction of $\rgo$ to $[s,u)$ is
concentrated on $\{t\in[s,u):F(t)=\ha t\}$. Now if (\ref{rhoclaim3})
would not hold, then $\tau:=\inf\{t\in[s,u):F(t)=F(s)+\eps\}$ would satisfy
$s<\tau<u$ for some $\eps>0$. But then $\rgo\big([s,\tau]\big)=0$ and hence
$F(\tau)=F(s)$, which is a contradiction.

Claim (\ref{rhoclaim1}) says that if $F(t)>\ha t$, then $F$ must stay constant
until the next time when $F(t)=\ha t$. Claim (\ref{rhoclaim3}) says that if
$F(t)<\ha t$, then $F$ must stay constant until the next time when it makes a
jump. Claim (\ref{rhoclaim2}) says that if $F$ makes a jump at time $t$, then
it jumps from $\ha t-\ha c$ to $\ha t+\ha c$ for some $c>0$. Using these
facts, it is easy to see that $\rgo$ must be of the form (\ref{gensol}).
\epro

\subsection{Frozen percolation on the MBBT}\label{S:MBBTfrz}

In this subsection, we prove a version of Theorem~\ref{T:dirfrz} for frozen
percolation on the MBBT, from which in the next subsection we will derive
Theorem~\ref{T:dirfrz}. We first need some definitions concerning general RTPs
corresponding to the RDE (\ref{MRDE}), similar to those introduced in
Subsection~\ref{S:nonend} for general RTPs corresponding to the RDE
(\ref{frzRDE}).

Let $(\tau_\ibf,\kappa_\ibf,Y_\ibf)_{\ibf\in\T}$ be an RTP corresponding to
the map $\chi$ from (\ref{chidef}) and a general solution $\rgo$ to the RDE
(\ref{MRDE}). Generalising the definition in (\ref{Sdef}), we set
\be\label{Sidef}
\S_\ibf:=\big\{\ibf j_1\cdots j_n\in\T:
j_m\leq\kappa_{\ibf j_1\cdots j_{m-1}}\ \forall 1\leq m\leq n\big\}.
\ee
Modifying the definition of $\T^t$ in (\ref{Fxdef}), in the present context,
we set
\be\ba{c}\label{StF}
\T^t:=\big\{\ibf\in\T:\kappa_\ibf=2\mbox{ or }\tau_\ibf\leq t\big\},\quad
\S^t_\ibf:=\T^t\cap\S_\ibf,\\[5pt]
\quand\F_{\rm y}:=\big\{\ibf\in\T:\kappa_\ibf=1,\ \tau_\ibf\geq Y_{\ibf 1}\big\}.
\ec
Similar to (\ref{Xupdef}), we define $I$-valued random variables
$(Y^\up_\ibf)_{\ibf\in\T}$ by
\be\label{Yupdef}
Y^\up_\ibf:=\inf\big\{t\in[0,1]:\ibf\percol{\S^t_\ibf\beh\F_{\rm y}}\infty\big\},
\ee
with $\inf\emptyset:=\infty$. Note that if $\ibf\in\S$, then in (\ref{Yupdef})
we can equivalently replace $\S^t_\ibf$ by $\S^t_\wurz=:\S^t$. At time
$t\in[0,1]$, we call points in $\T^t\beh\F_{\rm y}$ \emph{open}, points in
$\T^t\cap\F_{\rm y}$ \emph{frozen}, and all other points in $\T$
\emph{closed}. We call $\tau_\ibf$ the \emph{activation time} of $\ibf$ and
refer to $Y_\ibf$ and $Y^\up_\ibf$ as its \emph{burning time} and
\emph{percolation time}, respectively. Note that our modified definition of
$\T^t$ has the effect that \emph{branching points}, i.e., points $\ibf$ for
which $\kappa_\ibf=2$, are always open. The remaining \emph{blocking points},
i.e., points $\ibf$ for which $\kappa_\ibf=1$ are initially closed. At its
activation time, a blocking point $\ibf$ either freezes or opens, depending on
whether at that moment $\ibf 1$ is burnt or not.

It follows from the inductive relation (\ref{Yind}) that if $\kappa_\ibf=1$,
then $Y_\ibf>\tau_\ibf$, i.e., a blocking point can only burn after its
activation time. We see from the definition of $\F_{\rm y}$ in (\ref{StF}) and
the definition of the map $\chi$ in (\ref{chidef}) that if a blocking point
$\ibf$ burns at some time $Y_\ibf\in[0,1]$, then $\ibf$ must be open at that
time. Formula (\ref{chidef}) moreover implies that if a point $\ibf\in\T$
burns at some time $Y_\ibf\in[0,1]$, then starting at $\ibf$ there must be a
ray in $\S_\ibf$ consisting of points that burn at the same time as $\ibf$. By
our earlier remark and since branching points are always open, such a ray must
be open, which proves that (compare (\ref{XupX}))
\be\label{YupY}
Y^\up_\ibf\leq Y_\ibf\quad{\rm a.s.}\quad(\ibf\in\T).
\ee
The next proposition says that the opposite inequality holds only if $\rgo$ is
the special solution $\rro$ to the RDE defined in (\ref{rhodef}).

\bp[Percolation probability]
Let\label{P:perprob} $(\tau_\ibf,\kappa_\ibf,Y_\ibf)_{\ibf\in\T}$ be an RTP
corresponding to the map $\chi$ from (\ref{chidef}) and a solution $\rgo$ to
the RDE (\ref{MRDE}). Then
\be\label{perprob}
\P\big[Y^\up_\ibf\leq t\big]
=F(t)\vee\big(t-F(t)\big)
\qquad\big(t\in[0,1]\big),
\ee
where $F(t):=\rgo\big([0,t]\big)$ $(t\in[0,1])$. Moreover, one has
$Y^\up_\wurz=Y_\wurz$ a.s.\ if and only if $\rgo$ is the measure $\rro$ in
(\ref{rhodef}).
\ep

The proof of Proposition~\ref{P:perprob} needs some preparations. We will be
interested in the law of the open connected component of the root conditional
on the root not being burnt. In the next lemma we condition on the origin not
being burnt and calculate the probability that (i) the root is a branching
point, (ii) the root is a blocking point and its descendant is not burnt,
(iii) the root is a blocking point and its descendant is burnt. We show that
conditional on the event (ii), the activation time of the root is uniformly
distributed.

\bl[Law conditioned on not being burnt]
Let\label{L:unfroz} $(\tau_\ibf,\kappa_\ibf,Y_\ibf)_{\ibf\in\T}$ be an RTP
corresponding to the map $\chi$ from (\ref{chidef}) and a solution $\rgo$ to
the RDE (\ref{MRDE}). Then
\be\ba{l@{\quad}r@{\,}c@{\,}l}\label{unfroz}
{\rm(i)}&\dis\P\big[\kappa_\wurz=2\,\big|\,Y_\wurz>t\big]
&=&\dis\ha\big(1-F(t)\big),\\[5pt]
{\rm(ii)}&\dis\P\big[\kappa_\wurz=1,\ Y_1>t\,\big|\,Y_\wurz>t\big]
&=&\dis\ha,\\[5pt]
{\rm(iii)}&\dis\P\big[\kappa_\wurz=1,\ Y_1\leq t\,\big|\,Y_\wurz>t\big]
&=&\dis\ha F(t),
\ec
where $F(t):=\rgo\big([0,t]\big)$ $(t\in[0,1])$. Moreover,
\be\label{stilluni}
\P\big[\tau_\wurz\leq s\,\big|\,\kappa_\wurz=1,\ Y_1>t,\ Y_\wurz>t\big]=s
\qquad\big(s,t\in[0,1]\big).
\ee
\el

\bpro
One has
\bc
\dis\P\big[\kappa_\wurz=2,\ Y_\wurz>t\big]
&=&\dis\P\big[\kappa_\wurz=2,\ Y_1>t,\ Y_2>t\big]=\ha\big(1-F(t)\big)^2,\\[5pt]
\dis\P\big[\kappa_\wurz=1,\ Y_1>t,\ Y_\wurz>t\big]
&=&\ha\P[Y_1>t]=\ha\big(1-F(t)\big).
\ec
Dividing by $\P[Y_\wurz>t]=1-F(t)$ yields (\ref{unfroz})~(i) and (ii), and
the remaining formula follows since the total probability is one.
Since $\kappa_\wurz=1$ and $Y_1>t$ a.s.\ imply $Y_\wurz>t$,
and since $\tau_\wurz$ is independent of $Y_1,\kappa_\wurz$ and uniformly
distributed, we also obtain (\ref{stilluni}).
\epro

For $t\in[0,1]$, we inductively define $(\Ob^t_n)_{n\geq 0}$ by
$\Ob^t_0:=\{\wurz\}$ and
\be
\Ob^t_n:=\big\{\ibf j:
\ibf\in(\Ob^t_{n-1}\cap\T^t)\beh\F_{\rm y},\ 1\leq j\leq\kappa_\ibf\big\}.
\ee
We call $\Ob^t:=\bigcup_{n=0}^\infty\Ob^t_n$ the \emph{open component of the
  root}. Note that $\Ob^t_n$ consists of all descendants of open elements of
$\Ob^t_{n-1}$, while elements of $\Ob^t_{n-1}$ that are closed or frozen
produce no offspring. As a result, the root percolates at time $t\in[0,1]$ if
and only if $\Ob^t$ is infinite. The next lemma says that conditional on the
event that the root is not burnt, $(\Ob^t_n)_{n\geq 0}$ is a branching process
that can be subcritical, critical, or supercritical, depending on $t$ and our
choice of the solution $\rgo$ to the RDE (\ref{MRDE}).

\bl[The open unburnt component of the root]
Fix\label{L:unburnt} $t\in[0,1]$ and write $\Ob^t_n=\{\ibf
j:\ibf\in\Ob^t_{n-1},\ 1\leq j\leq\la^t_\ibf\}$ with
$\la^t_\ibf\in\{0,1,2\}$. If $(\U_k)_{0\leq k<n}$ is a possible realization of $(\Ob^t_k)_{0\leq k<n}$, then
conditional on the event
$\Ai^t:=\{Y_\wurz>t,\ (\Ob^t_k)_{0\leq k<n}=(\U_k)_{0\leq k<n}\}$,
the random variables $(\la^t_\ibf)_{\ibf\in\U_{n-1}}$ are i.i.d.\ with law
\be\ba{c}\label{unburnt}
\dis\P[\la^t_\ibf=0\,|\,\Ai^t]=\ha\big(1-t+F(t)\big),
\quad\P[\la^t_\ibf=1\,|\,\Ai^t]=\ha t,\\[5pt]
\dis\quad\P[\la^t_\ibf=2\,|\,\Ai^t]=\ha\big(1-F(t)\big),
\ec
where $F(t):=\rgo\big([0,t]\big)$ $(t\in[0,1])$.
\el

\bpro
Fix $t\in[0,1]$. We claim that $Y_\wurz>t$ implies $Y_\ibf>t$ for all
$\ibf\in\Ob^t$. Indeed, if $\ibf\in\Ob^t_{n-1}$ is open and not burnt, then
all its descendants must be unburnt, while elements that are not open have no
descendants in $\Ob^t_n$, so the claim follows by induction.

Fix $(\U_k)_{0\leq k<n}$ and define $\Ai^t$ as in the lemma, which by what we
have just proved is the same as the event
\be
\Ai^t=\big\{Y_\ibf>t\ \forall\ibf\in\U,
\ (\Ob^t_k)_{0\leq k<n}=(\U_k)_{0\leq k<n}\big\},
\ee
where $\U:=\bigcup_{0\leq k<n}\U_k$. By Lemma~\ref{L:unfroz}, independently
for each $\ibf\in\U_{n-1}$,
\be\ba{l@{\quad}r@{\,}c@{\,}l}
{\rm(i)}&\dis\P\big[\kappa_\ibf=2\,\big|\,\Ai^t\big]
&=&\dis\ha\big(1-F(t)\big),\\[5pt]
{\rm(ii)}&\dis\P\big[\kappa_\ibf=1,\ \tau_\ibf\leq t,\ Y_{\ibf 1}>t
\,\big|\,\Ai^t\big]
&=&\dis\ha t,\\[5pt]
{\rm(iii)}&\dis
\P\big[\kappa_\ibf=1,\ \tau_\ibf>t,\ Y_{\ibf 1}>t\,\big|\,\Ai^t\big]
&=&\dis\ha(1-t),\\[5pt]
{\rm(iv)}&\dis\P\big[\kappa_\ibf=1,\ Y_{\ibf 1}\leq t\,\big|\,\Ai^t\big]
&=&\dis\ha F(t),
\ec
which are the conditional probabilities that (i) $\ibf$ is a branching point,
(ii) $\ibf$ is an open blocking point, (iii) $\ibf$ is a closed blocking point
and its descendant is not burnt, (iv) $\ibf$ is a blocking point and its
descendant is burnt, which is only possible if $\ibf$ is closed or frozen.
Since $\la^t_\ibf=2$ in case (i), $\la^t_\ibf=1$ in case (ii), and
$\la^t_\ibf=0$ in the remaining cases, the lemma follows.
\epro

\bpro[Proof of Proposition~\ref{P:perprob}]
By (\ref{YupY}),
\bc\label{FmF}
\dis\P\big[Y^\up_\wurz\leq t\big]
&=&\dis\P\big[Y_\wurz\leq t\big]+\P\big[Y_\wurz>t\big]
\P\big[Y^\up_\wurz\leq t\,\big|\,Y_\wurz>t\big]\\[5pt]
&=&\dis F(t)+\big(1-F(t)\big)
\P\big[\Ob^t_n\neq\emptyset\ \forall n\geq 0\,\big|\,Y_\wurz>t\big].
\ec
By Lemma~\ref{L:unburnt}, the probability
\be
p:=\P\big[\Ob^t_n\neq\emptyset\ \forall n\geq 0\,\big|\,Y_\wurz>t\big]
\ee
is the survival probability of a branching process with offspring distribution
as in (\ref{unburnt}). It is well-known \cite[Thm~III.4.1]{AN72} that the
survival probability is the largest solution in $[0,1]$ of the equation
$\Psi(p)=p$, where (compare formula (\ref{Psidef}) in the appendix)
\be
\Psi(p):=\ha\big(1-F(t)\big)p(1-p)-\ha\big(1-t+F(t)\big)p.
\ee
Assuming that $F(t)<1$, it follows that
\be
p=0\vee\Big\{1-\frac{1-t+F(t)}{1-F(t)}\Big\}=0\vee\frac{t-2F(t)}{1-F(t)}.
\ee
Inserting this into (\ref{FmF}) we arrive at (\ref{perprob}). This argument
does not work if $F(t)=1$, which by Proposition~\ref{P:genRDE} is only
possible if $t=1$ and $\rgo=\de_1$. In this case, no freezing takes place
until at time $t=1$ all $\ibf\in\T$ are open, so the left- and right-hand
sides of (\ref{perprob}) are both trivially equal to one.

Formula (\ref{YupY}) says that $Y^\up_\wurz\leq Y_\wurz$ a.s., so
we have $Y^\up_\wurz=Y_\wurz$ a.s.\ if and only if
\be
\P\big[Y^\up_\wurz\leq t\big]
=\P[Y_\wurz\leq t]=F(t)\qquad\big(t\in[0,1]\big),
\ee
which by (\ref{perprob}) happens if and only if $F(t)\geq\ha t$ $(t\in[0,1])$.
By Proposition~\ref{P:genRDE}, the only solution to the RDE (\ref{MRDE}) with
this property is the measure $\rro$ in (\ref{rhodef}).
\epro

\subsection{Frozen percolation on the binary tree}\label{S:binfrz}

In this subsection we derive Theorem~\ref{T:dirfrz} from
Proposition~\ref{P:perprob}. Our main tool is a coupling between, one the one
hand, an RTP $(\tau_\ibf,\kappa_\ibf,Y_\ibf)_{\ibf\in\T}$ corresponding to the
map $\chi$ from (\ref{chidef}), and on the other hand, an RTP
$(\tau_\ibf,X_\ibf)_{\ibf\in\T}$ corresponding to the map $\ga$ from
(\ref{gadef}). We first describe the main idea of the construction and then
fill in the technical details.

It is easy to see that for an RTP corresponding to the map $\chi$ from
(\ref{chidef}), the number of blocking points between two consecutive
branching points is geometrically distributed with parameter $1/2$. Imagine,
for the moment, that instead there would always be exactly one blocking point
between two consecutive branching points. Then, comparing (\ref{gadef}) and
(\ref{chidef}), one can check that the inductive relation satisfied by the
burning times $(Y_\ibf)_{\ibf\in\S,\ \kappa_\ibf=1}$ of blocking points would
be exactly the same as the inductive relation satisfied by the burning times
$(X_\ibf)_{\ibf\in\T}$ of arbitrary points in an RTP corresponding to the map
$\ga$ from (\ref{gadef}). Inspired by this, starting from an RTP corresponding
to the map $\chi$ from (\ref{chidef}), we will construct an associated RTP
corresponding to the map $\ga$ from (\ref{gadef}) along the following steps:
\begin{enumerate}
\item If there are two or more blocking points between two consecutive
  branching points, then we replace them by one point, whose new activation
  time is the maximum of the activation times of the blocking points it
  replaces.
\item If there are no blocking points between two consecutive branching
  points, then we add one such point, and assign it an activation time that is
  uniformly distributed on $[-1,0]$.
\item We transform the activation times that we obtain by this procedure using
  a monotone mapping from $[-1,1]$ to $[0,1]$, which has the result that the
  transformed times are uniformly distributed on $[0,1]$.
\end{enumerate}

We now formulate this a little more precisely. Let
\be
1_n:=\underbrace{1\cdots 1}_{\mbox{$n$ times}}
\ee
denote the word of length $n\geq 0$ that contains only 1's. For each
$\ibf\in\S$, we set
\be
b(\ibf):=\ibf 1_{N_\ibf}\quad\mbox{with}\quad
N_\ibf:=\inf\{n\geq 0:\kappa_{\ibf 1_n}=2\}.
\ee
In words, $b(\ibf)$ is the next branching point above $\ibf$ (which may
be $\ibf$ itself). We inductively define a map $\psi:\T\to\S$ by
$\psi(\wurz)=\wurz$ and
\be\label{psidef}
\psi(\ibf j):=b\big(\psi(\ibf)\big)j\quad(\ibf\in\T,\ j=1,2).
\ee
Note that points of the form $\psi(\ibf)$ with $\ibf\in\T\beh\{\wurz\}$ are
direct descendants of branching points, and $N_{\psi(\ibf)}$ is the number of
steps we have to walk up from $\psi(\ibf)$ to reach the next branching point.

We let $(\ti\tau_\ibf)_{\ibf\in\T}$ be an i.i.d.\ collection of uniformly
distributed $[-1,0]$-valued random variables, independent of everything else.
For each $\ibf\in\T$, we define
\be\label{sigdef}
\sig_\ibf:=\left\{\ba{ll}
\dis\max\big\{\tau_{\psi(\ibf)1_n}:0\leq n\leq N_{\psi(\ibf)}-1\big\}
\quad&\mbox{if }N_{\psi(\ibf)}\geq 1,\\[5pt]
\dis\ti\tau_\ibf\quad&\mbox{otherwise,}
\ea\right.
\ee
i.e., $\sig_\ibf$ is the maximum of the activation times of blocking points
that lie directly below the branching point $b(\psi(\ibf))$, if there are any,
and $\sig_\ibf=\ti\tau_\ibf$ otherwise. For each $\ibf\in\T$, the number
$N_{\psi(\ibf)}$ of blocking points that lie below the branching point
$b(\psi(\ibf))$ is geometrically distributed with parameter $1/2$, and the
values of their activation times are i.i.d.\ uniformly distributed on $[0,1]$
and independent of $N_{\psi(\ibf)}$. These quantities are moreover independent
for different $\ibf\in\T$. As a result, the $(\sig_\ibf)_{\ibf\in\T}$ are
i.i.d.\ with distribution function
\be\label{sigF}
\P[\sig_\ibf<s]=\Hh(s):=\left\{\ba{ll}
\dis\ha(1+s)\quad&\mbox{if }s\in[-1,0],\\[5pt]
\dis\frac{1}{2-s}\quad&\mbox{if }s\in[0,1],
\ea\right.
\ee
where we have used the calculation in (\ref{maxgeom}) and we extend the
function $\Hh:[0,1]\to[\ha,1]$ from Lemma~\ref{L:equivRDE} into a function
$\Hh:[-1,1]\to[0,1]$.

\bp[Coupling of RTPs]
Let\label{P:RTPcoup} $(\tau_\ibf,\kappa_\ibf,Y_\ibf)_{\ibf\in\T}$ be an RTP
corresponding to the map $\chi$ from (\ref{chidef}) and any solution to the RDE
(\ref{MRDE}). Let $(\ti\tau_\ibf)_{\ibf\in\T}$ be an independent
i.i.d.\ collection of uniformly distributed $[-1,0]$-valued random variables,
and let $\psi:\T\to\T$, $(\sig_\ibf)_{\ibf\in\T}$, and $\Hh:[-1,1]\to[0,1]$ be
defined as in (\ref{psidef}), (\ref{sigdef}), and (\ref{sigF}). Then setting
\be\label{RTPcoup}
\ov\tau_\ibf:=\Hh(\sig_\ibf)\quand X_\ibf:=\Hh(Y_{\psi(\ibf)})\qquad(\ibf\in\T)
\ee
defines an RTP $(\ov\tau_\ibf,X_\ibf)_{\ibf\in\T}$ corresponding to the map
$\ga$ from (\ref{gadef}). Moreover, any RTP corresponding to $\ga$ is equal
in distribution to an RTP constructed in this way. Finally, one has
\be\label{percXY}
X^\up_\ibf:=\Hh(Y^\up_{\psi(\ibf)})\qquad(\ibf\in\T),
\ee
where $X^\up_\ibf$ is defined in (\ref{Xupdef}) and $Y^\up_{\psi(\ibf)}$ is defined in (\ref{Yupdef}).
\ep

\bpro
We claim that $(Y_{\psi(\ibf)})_{\ibf\in\T}$ satisfy the inductive
relation
\be\label{Ypsind}
Y_{\psi(\ibf)}=\ga[\sig_\ibf]\big(Y_{\psi(\ibf 1)},Y_{\psi(\ibf 2)}\big)
\qquad(\ibf\in\T),
\ee
where we define $\ga[t](x,y)$ as in (\ref{gadef}) also for negative $t$.
Indeed, if $N_{\psi(\ibf)}=0$, then $\sig_\ibf\leq 0$ while $Y_{\psi(\ibf
  1)},Y_{\psi(\ibf 2)}>0$ a.s., and
\be\label{Yps1}
Y_{\psi(\ibf)}=\chi[2]\big(Y_{\psi(\ibf)1},Y_{\psi(\ibf)2}\big)
=Y_{\psi(\ibf 1)}\wedge Y_{\psi(\ibf 2)}.
\ee
On the other hand, if $N_{\psi(\ibf)}\geq 1$, then
\bc\label{Yps2}
\dis Y_{\psi(\ibf)}
&=&\dis\chi[\tau_{\psi(\ibf)},1]\circ\cdots\circ
\chi[\tau_{\psi(\ibf)1_{N_{\psi(\ibf)}}},1]
\circ\chi[2]\big(Y_{\psi(\ibf)1_{N_{\psi(\ibf)}}1},Y_{\psi(\ibf)1_{N_{\psi(\ibf)}}2}\big)\\[5pt]
&=&\dis\ga[\tau_{\psi(\ibf)}\vee\cdots\vee\tau_{\psi(\ibf)1_{N_{\psi(\ibf)}}}]
\big(Y_{\psi(\ibf 1)},Y_{\psi(\ibf 2)}\big).
\ec
Using (\ref{Ypsind}) and (\ref{moninv}),
we conclude that $(X_\ibf)_{\ibf\in\T}$ satisfy the inductive
relation (\ref{Xind}). By (\ref{sigF}), the random variables
$(\ov\tau_\ibf)_{\ibf\in\T}$ are i.i.d.\ and uniformly distributed on $[0,1]$.
Moreover, for any finite rooted subtree $\U\sub\T$, the
r.v.'s $(X_\ibf)_{\ibf\in\pa\U}$ are independent of
$(\ov\tau_\ibf)_{\ibf\in\pa\U}$ and i.i.d.

This completes the proof that
$(\ov\tau_\ibf,X_\ibf)_{\ibf\in\T}$ is an RTP corresponding to the map $\ga$
from (\ref{gadef}). Using Lemma~\ref{L:equivRDE}, we see that every RTP
$(\tau_\ibf,X_\ibf)_{\ibf\in\T}$ corresponding to the map $\ga$ from
(\ref{gadef}) and some solution $\ngu$ to the RDE (\ref{frzRDE})
is equal in distribution to an RTP constructed in this way.

To prove also (\ref{percXY}), we observe that the frozen set $F_{\rm x}$ from
(\ref{Fxdef}) for the RTP $(\ov\tau_\ibf,X_\ibf)_{\ibf\in\T}$ is given by
\be\ba{l}
\Fx
=\big\{\ibf\in\T:\ov\tau_\ibf\geq X_{\ibf 1}\wedge X_{\ibf 2}\big\}
=\big\{\ibf\in\T:\sig_\ibf\geq Y_{\psi(\ibf 1)}\wedge Y_{\psi(\ibf 2)}\big\}\\[5pt]
\dis\quad=\big\{\ibf\in\T:N_{\psi(\ibf)}\geq 1,
\ \tau_{\psi(\ibf)1_n}\geq Y_{\psi(\ibf)1_{n+1}}
\mbox{ for some }0\leq n<N_{\psi(\ibf)}\big\}\\[5pt]
\dis\quad=\big\{\ibf\in\T:N_{\psi(\ibf)}\geq 1,\ \psi(\ibf)1_n\in\F_{\rm y}
\mbox{ for some }0\leq n<N_{\psi(\ibf)}\big\},
\ec
and hence at time $t$ there exists a ray in $\S^t_{\psi(\ibf)}\beh\F_{\rm y}$
starting at $\psi(\ibf)$ if and only if at time $s:=\Hh(t)$ there exists a ray
in $\T^s\beh\Fx$ starting at $\ibf$.
\epro

\bpro[Proof of Theorem~\ref{T:dirfrz}]
By Lemma~\ref{L:half}, $\mu$ is concentrated on $I'=[\ha,1]\cup\{\infty\}$.
Let $\mu'$ be the image of $\mu$ under the inverse of the map $H:I\to I'$
defined in Lemma~\ref{L:equivRDE}. Then $\mu'$ solves the RDE (\ref{MRDE}).
Let $(\tau_\ibf,\kappa_\ibf,Y_\ibf)_{\ibf\in\T}$ be the RTP corresponding to the
map $\chi$ from (\ref{chidef}) and the measure $\mu'$. We couple this RTP to
$(\tau_\ibf,X_\ibf)_{\ibf\in\T}$ as in Proposition~\ref{P:RTPcoup}.
Since the function $\Hh$ is strictly increasing, we see that
$X^\up_\wurz=X_\wurz$ a.s.\ if and only if $Y^\up_\wurz=Y_\wurz$ a.s.
By Proposition~\ref{P:perprob} this is equivalent to $\mu'$ being the measure
$\rro$ in (\ref{rhodef}), which by Lemma~\ref{L:rhonu} is equivalent to
$\ngu$ being the measure $\nnu$ in (\ref{nudef}).
\epro

\subsection{Scale invariance of the MBBT}\label{S:frzMBBT}

The aim of the present subsection is to prove Proposition~\ref{P:scale}
and Lemma~\ref{L:scinv} about scale invariance of (frozen percolation on) the
MBBT. Lemma~\ref{L:scinv}, in particular, allows us to identify the nontrivial
fixed point $\rho^{(2)}_2$ from Theorem~\ref{T:scalefix} as $\un\rho^{(2)}$.
Combining this with Lemma~\ref{L:nurho}, we also obtain an explicit expression
for $\un\nu^{(2)}$. As a preparation for this, we first prove
Lemma~\ref{L:perc}.\med

\bpro[Proof of Lemma~\ref{L:perc}]
It is well-known \cite[Thm~III.4.1]{AN72} that the survival probability is the
largest solution in $[0,1]$ of the equation $\Psi(p)=p$, where (compare
formula (\ref{Psidef}) in the appendix)
\be\label{binPsi}
\Psi(p)=\{(1-p)-(1-p)^2\}+(1-t)\{(1-p)-(1-p)^0\}=p(1-p)-(1-t)p.
\ee
Since $\Psi(p)=0$ has two roots, $p=0$ and $p=t$, we conclude that the
survival probability is $t$.
%
%
\epro

We next turn our attention to the proof of Proposition~\ref{P:scale}.  Let
$(\Ti,\Pi)$ be the MBBT. If we cut $\Ti$ at points in $\Pi_t$, then the
connected component of the root is the family tree of a continuous-time
branching process where particles split into two with rate one and die with
rate $1-t$. The tree $\Ti'$ defined in (\ref{Tiac}) is the \emph{skeleton} of
this process. It is well-known that $\Ti'$ is the family tree of a branching
process, which is known as the \emph{skeletal process}. There exist standard
ways to find the skeletal process associated with a given branching process.
Using these, it is easy to check that $\Ti'$ is the family tree of a
binary branching process with branching rate $t$. In Appendix~\ref{A:skel}, we
outline a proof of this fact along these lines, with references to the
relevant literature.

To prove Proposition~\ref{P:scale}, we need a bit more, however, since we need
to determine the joint law of $\Ti'$ and $\Pi'$. To prove also
Lemma~\ref{L:scinv}, we will moreover need a scaling property of RTPs
corresponding to the map $\chi$ in (\ref{chidef}) and law $\rro$ from
(\ref{rhodef}). In view of this, we find it more convenient to give
self-contained proofs of Proposition~\ref{P:scale} and Lemma~\ref{L:scinv},
not referring to the abstract theory of skeletal processes.

Let $(\tau_\ibf,\kappa_\ibf,Y_\ibf)_{\ibf\in\T}$ be the RTP corresponding to
the map $\chi$ from (\ref{chidef}) and law $\rro$ from (\ref{rhodef}), and let
$(\ell_\ibf)_{\ibf\in\T}$ be an independent i.i.d.\ collection of
exponentially distributed random variables with mean $1/2$. As in
Subsection~\ref{S:MBBT}, we use the random variables
$(\tau_\ibf,\kappa_\ibf,\ell_\ibf)_{\ibf\in\T}$ to define an MBBT
$(\Ti,\Pi)$. In particular, $\Ti$ is the family tree of a branching process
$(\nab\S_h)_{h\geq 0}$ where $\S$, defined in (\ref{Sdef}), is the collection
of all individuals that will ever live.

We fix $0<t\leq 1$ and define
\be\label{Yast}
Y^\ast_\ibf:=\left\{\ba{ll}
t^{-1}Y_\ibf\quad&\mbox{if }Y_\ibf\leq t,\\[5pt]
\infty\quad&\mbox{otherwise.}
\ea\right.
\ee
We also define $(\Ti',\Pi')$ as in (\ref{Tiac}) and define
$(\Ti^\ast,\Pi^\ast)$ by
\be
\Ti^\ast:=\big\{(x,th):(x,h)\in\Ti'\big\},\quad
\Pi^\ast:=\big\{(x,th,t^{-1}\tau_{(x,h)}):(x,h,\tau_{(x,h)})\in\Pi'\big\}.
\ee
As in Proposition~\ref{P:scale}, we view $(\Ti',\Pi')$ and
$(\Ti^\ast,\Pi^\ast)$ as marked metric spaces, i.e., we do not care about the
precise labeling of elements of $\Ti'$ or
$\Ti^\ast$. Proposition~\ref{P:scale} can be rephrased by saying that the
conditional law of $(\Ti^\ast,\Pi^\ast)$ given
$\wurz\percol{\Ti\beh\Pi_t}\infty$ is equal to the original law of
$(\Ti,\Pi)$. The following lemma says that in a sense, $(\Ti^\ast,\Pi^\ast)$
contains all relevant information about $Y^\ast_\wurz$.

\bl[Relevant information]
One\label{L:relinfo} has
\be
\P\big[Y^\ast_\wurz\in\,\cdot\,\big|\,(\tau_\ibf,\kappa_\ibf)_{\ibf\in\T}\big]
=\P\big[Y^\ast_\wurz\in\,\cdot\,\big|\,(\Ti^\ast,\Pi^\ast)\big]
\quad{\rm a.s.}
\ee
\el

The following proposition extends Proposition~\ref{P:scale} to a scaling
property of the joint law of $(Y^\ast_\wurz,\Ti^\ast,\Pi^\ast)$. In
particular, this implies Proposition~\ref{P:scale}.

\bp[Scaling of the joint law]
One\label{P:jointscale} has
\be
\P\big[(Y^\ast_\wurz,\Ti^\ast,\Pi^\ast)\in\,\cdot\,\big|\,
\Ti^\ast\neq\emptyset\big]=\P\big[(Y_\wurz,\Ti,\Pi)\in\,\cdot\,\big].
\ee
\ep

Before we prove Lemma~\ref{L:relinfo} and Proposition~\ref{P:jointscale}, we
first show how they imply Lemma~\ref{L:scinv}.\med

\bpro[Proof of Lemma~\ref{L:scinv}]
Conditional on $(\tau_\ibf,\kappa_\ibf)_{\ibf\in\T}$, let $(Y'_\ibf)_{\ibf\in\T}$
be an independent copy of $(Y_\ibf)_{\ibf\in\T}$. Then, according to the
definitions in (\ref{ovnu}) and (\ref{unnu})
\be
\ov\rro^{(2)}=\P\big[(Y_\wurz,Y_\wurz)\in\,\cdot\,\big]
\quand
\un\rro^{(2)}=\P\big[(Y_\wurz,Y'_\wurz)\in\,\cdot\,\big].
\ee
Clearly, these measures are symmetric and their one-dimensional marginals are
given by $\rro$. It remains to show that they have the scaling property
(\ref{scinv}). The claim for $\ov\rro^{(2)}$ follows easily from the
fact that $Y_\wurz$ has the law $\rro$ in (\ref{rhodef}). It remains to prove
the statement for $\un\rro^{(2)}$.

Fix $r,s,t\in[0,1]$. Since $Y_\wurz=\infty$ a.s.\ on the complement of the event
$\wurz\percol{\Ti\beh\Pi_t}\infty$, we have
\be\ba{l}
\dis\P\big[(Y_\wurz,Y'_\wurz)\in[0,tr]\times[0,ts]\big]\\[5pt]
\dis\quad=\P\big[(Y_\wurz,Y'_\wurz)\in[0,tr]\times[0,ts]
\,\big|\,\wurz\percol{\Ti\beh\Pi_t}\infty\big]
\,\P\big[\wurz\percol{\Ti\beh\Pi_t}\infty\big].
\ec
Here $\P\big[\wurz\percol{\Ti\beh\Pi_t}\infty\big]=t$ by Lemma~\ref{L:perc},
so to show that $\un\rro^{(2)}$ has the scaling property (\ref{scinv}), it
suffices to show that
\be\label{scaleproof}
\P\big[(Y_\wurz,Y'_\wurz)\in[0,tr]\times[0,ts]
\,\big|\,\wurz\percol{\Ti\beh\Pi_t}\infty\big]
=\P\big[(Y_\wurz,Y'_\wurz)\in[0,r]\times[0,s]\big].
\ee
Since $Y_\wurz$ and $Y'_\wurz$ are conditionally independent given the
\si-field generated by $(\tau_\ibf,\kappa_\ibf)_{\ibf\in\T}$, and since the
event that $\wurz\percol{\Ti\beh\Pi_t}\infty$ is measurable w.r.t.\ this
\si-field, we can rewrite the left-hand side of (\ref{scaleproof}) as
\be\ba{l}\label{condind}
\E\Big[\P\big[Y_\wurz\in[0,tr]\,\big|\,
\wurz\percol{\Ti\beh\Pi_t}\infty,\;(\tau_\ibf,\kappa_\ibf)_{\ibf\in\T}\big]\\
\dis\phantom{\E\Big[}\cdot
\P\big[Y_\wurz\in[0,ts]\,\big|\,\wurz\percol{\Ti\beh\Pi_t}\infty,
\;(\tau_\ibf,\kappa_\ibf)_{\ibf\in\T}\big]\Big]\\[5pt]
\dis\quad\stackrel{1}{=}\E\Big[\P\big[Y^\ast_\wurz\in[0,r]\,\big|\,
\Ti^\ast\neq\emptyset,\;(\Ti^\ast,\Pi^\ast)\big]
\P\big[Y^\ast_\wurz\in[0,s]\,\big|\,\Ti^\ast\neq\emptyset,
\;(\Ti^\ast,\Pi^\ast)\big]\Big]\\[5pt]
\dis\quad\stackrel{2}{=}\E\Big[\P\big[Y_\wurz\in[0,r]\,\big|\,(\Ti,\Pi)\big]
\P\big[Y_\wurz\in[0,s]\,\big|\,(\Ti,\Pi)\big]\Big]\\[5pt]
\dis\quad\stackrel{3}{=}\E\Big[\P\big[Y_\wurz\in[0,r]\,\big|\,
(\tau_\ibf,\kappa_\ibf)_{\ibf\in\T}\big]
\P\big[Y_\wurz\in[0,s]\,\big|\,
(\tau_\ibf,\kappa_\ibf)_{\ibf\in\T}\big]\Big],
\ec
which equals the right-hand side of (\ref{scaleproof}). Here, in step~1, we
have used the definition of $Y^\ast_\wurz$ in (\ref{Yast}), as well as the
fact that the event $\{\wurz\percol{\Ti\beh\Pi_t}\infty\}$ is the same as the
event $\{\Ti^\ast\neq\emptyset\}$, which is measurable with respect to the
\si-fields generated by $(\tau_\ibf,\kappa_\ibf)_{\ibf\in\T}$ and
$(\Ti^\ast,\Pi^\ast)$, and we have applied Lemma~\ref{L:relinfo}. Step~2
follows from Proposition~\ref{P:jointscale}. In step~3 we have again applied
Lemma~\ref{L:relinfo} but this time for $t=1$, in which case
$(Y^\ast_\wurz,\Ti^\ast,\Pi^\ast)=(Y_\wurz,\Ti,\Pi)$.
\epro

\bpro[Proof of Propositions~\ref{P:scale} and \ref{P:jointscale}]
Let $(\tau_\ibf,\kappa_\ibf,Y_\ibf)_{\ibf\in\T}$ be the RTP corresponding to
the map $\chi$ from (\ref{chidef}) and law $\rro$ from (\ref{rhodef}), and let
$(\ell_\ibf)_{\ibf\in\T}$ be an independent i.i.d.\ collection of
exponentially distributed random variables with mean $1/2$. Fix
$t\in(0,1]$. For any $A\sub\T$ and $\ibf \prec \jbf\in\T$, we
write $\ibf\percol{A}\jbf$ if there exist $\ibf_0,\ldots,\ibf_n\in A$, $n\geq
0$, such that $\ibf_0=\ibf$, $\ibf_n=\jbf$, and $\lvec\ibf_k=\ibf_{k-1}$
$(k=1,\ldots,n)$. Let us say that $\ibf\in\T$ is \emph{active} if it is either
open or frozen, i.e., if $\kappa_\ibf=2$ or $\tau_\ibf\leq t$, and let
\be
\A:=\big\{\ibf\in\T:\wurz\percol{\S^t}\ibf\percol{\S^t}\infty\big\},
\ee
with $\S^t$ as in (\ref{StF}) denote the collection of points that lie on an
active ray in $\S$ starting at the root. Note that by Lemma~\ref{L:perc}, the
probability that $\A$ is not empty is $t$. We give each $\ibf\in\A$ a
\emph{type} $\om_\ibf\in[0,t)\cup\{1,2\}$, which is defined as follows:
\be\label{om12}
\om_\ibf:=\left\{\ba{ll}
\dis\tau_\ibf\quad&\mbox{if }\kappa_\ibf=1,\\
\dis 1\quad&\mbox{if }\kappa_\ibf=2\mbox{ and $\{\ibf 1,\ibf 2\}\cap\A$ has
  precisely one element},\\
\dis 2\quad&\mbox{if }\kappa_\ibf=2\mbox{ and $\ibf 1,\ibf 2$ are both elements
  of $\A$.}
\ea\right.
\ee
Let $\A_n:=\{\ibf\in\A:|\ibf|=n\}$. We claim that conditional on the event
that $\A\neq\emptyset$, the process $(\A_n)_{n\geq 0}$ with the types
assigned to its elements is a multitype branching process with the following
description. In each generation, we first assign types to the particles that
are alive in an i.i.d.\ fashion according to the law
\be\label{offtype}
\Pb[\om\leq s]:=\ha s\quad\big(s\in[0,t]\big),\quad
\Pb[\om=1]:=1-t,\quand\Pb[\om=2]=\ha t,
\ee
and then let particles of type 2 produce two offspring while all other particles
produce one offspring. To see this, observe that by Lemma~\ref{L:perc}, for
each $\ibf\in\T$ and $s\in[0,t]$,
\bc\label{unnorm}
\dis\P\big[\kappa_\ibf=1,\ \tau_\ibf\leq s,\ \ibf 1\percol{\S^t_\ibf}\infty\big]
&=&\dis\ha st,\\[5pt]
\dis\P\big[\kappa_\ibf=2,\ \ibf 1\percol{\S^t_\ibf}\infty
\mbox{ or }\ibf 2\percol{\S^t_\ibf}\infty\mbox{ but not both}\big]
&=&\dis t(1-t),\\[5pt]
\dis\P\big[\kappa_\ibf=2,\ \ibf 1\percol{\S^t_\ibf}\infty
\mbox{ and }\ibf 2\percol{\S^t_\ibf}\infty\big]&=&\dis\ha t^2.
\ec
If we condition on $(\A_k)_{0\leq k\leq n}$ and also on the types of
particles in generations $0,\ldots,n-1$, then the types of particles in the
$n$-th generation are i.i.d.\ and their law is the distribution in
(\ref{unnorm}) normalised to make it a probability law, which is the
distribution $\Pb$ in (\ref{offtype}).

Let $(\Ti,\Pi)$ be the MBBT constructed as in Subsection~\ref{S:MBBT} from the
random variables $(\tau_\ibf,\kappa_\ibf,\ell_\ibf)_{\ibf\in\T}$, and let
$(\Ti',\Pi')$ be as in (\ref{Tiac}). Then $(\Ti',\Pi')$ is uniquely determined
by the branching process $\A$ and the types $\om_\ibf$ and lifetimes
$\ell_\ibf$ of elements $\ibf\in\A$. However, $\A$ contains, in a sense, too
much information, since points $\ibf\in\A$ with type $\om_\ibf=1$ are not
visible in $(\Ti',\Pi')$. To remedy this, we need a procedure to remove these
points, which we describe now.

For $\ibf\in\A$ with $\om_\ibf\neq 2$, let $f(\ibf):=\ibf j$ where $j$ is the
unique element of $\{1,2\}$ such that $\ibf j\in\A$, and let
\be\label{nb}
b(\ibf):=f^{n(\ibf)}(\ibf)\quad\mbox{with}\quad
n(\ibf):=\inf\{k\geq 0:\om_{f^k(\ibf)}\neq 1\}
\ee
denote the next point above $\ibf$ that is not of type 1. Let
$\B:=\{\ibf\in\A:\om_\ibf\neq 1\}$. We inductively define a map
$\psi:\B\to\T$ by $\psi\big(b(\wurz)\big):=\wurz$ and
\be\ba{r@{\,}c@{\,}ll}
\dis\psi\big(b(\ibf j)\big)&:=&\dis\psi(\ibf)j\quad(j=1,2)
\quad&\dis\mbox{if }\om_\ibf=2,\\[5pt]
\dis\psi\big(b(\ibf 1)\big)&:=&\dis\psi(\ibf)1
\quad&\dis\mbox{if }\om_\ibf\in[0,t).
\ec
We let $\S'$ denote the image of $\B$ under the map $\psi$ and assign types to
the elements of $\S'$ by
\be
\om'_{\psi(\ibf)}:=\om_\ibf\qquad(\ibf\in\B).
\ee
We also define new lifetimes by
\be
\ell'_{\psi(\ibf)}:=\sum_{k=0}^{n(f(\ibf))}\ell_{f^k(\ibf)}
\ee
where $n(\ibf)$ is defined as in (\ref{nb}). Then the set $\S'$ and the random
variables $(\om'_\ibf,\ell'_\ibf)_{\ibf\in\S}$ contain precisely the
information needed to construct $(\Ti',\Pi')$, and nothing more.

Let $\S'_n:=\{\ibf\in\S':|\ibf|=n\}$. The process $(\S'_n)_{n\geq 0}$ inherits
the branching property from the process $(\A_n)_{n\geq 0}$. To get the new
generation, we first assign i.i.d.\ types to the particles in the present
generation according to the law
\be\label{offtype2}
\Pb'[\om\leq s]:=\frac{s}{2t}\quad\big(s\in[0,t]\big),\quad
\Pb'[\om=2]=\ha,
\ee
which is the law in (\ref{offtype}) conditioned on $\om\neq 1$, and then let
particles with type in $[0,t)$ and $\{2\}$ produce one or two offspring,
respectively. Each lifetime $\ell'_\ibf$ is the sum of a geometric number of
exponentially distributed random variables. From this, it is easy to see that
conditional on $\S'$ and the types, the lifetimes $(\ell'_\ibf)_{\ibf\in\S'}$
are i.i.d.\ and exponentially distributed with mean $\ha t^{-1}$.
Since the random tree $\Ti'$ is the family tree of the branching process
$(\S'_n)_{n\geq 0}$ with the lifetimes $(\ell'_\ibf)_{\ibf\in\S'}$, and the
Poisson set $\Pi'$ records points with type $\om_\ibf\in[0,t)$ together with
their activation times $\tau'_\ibf:=\om_\ibf\in[0,t)$, this completes proof of
Proposition~\ref{P:scale}.

We could have obtained Proposition~\ref{P:scale} faster by referring to the
the abstract theory of skeletal processes (see Appendix~\ref{A:skel}). The
advantage of our explicit construction, however, is that it also easily
yields the stronger statement of Proposition~\ref{P:jointscale}. To see this,
we define $(Y'_\ibf)_{\ibf\in\S'}$ by
\be
Y'_{\psi(\ibf)}:=\left\{\ba{ll}
\dis Y_\ibf\quad&\dis\mbox{if }Y_\ibf\leq t\\[5pt]
\dis\infty\quad&\dis\mbox{otherwise.}\ea\right.
\qquad(\ibf\in\B).
\ee
Since we started from an RTP corresponding to the law $\rro$ from
(\ref{rhodef}), and since $Y_\ibf>t$ a.s.\ on the complement of the event
$\ibf\percol{\S_t}\infty$, we see that conditional on $(\S'_k)_{0\leq k\leq n}$
and the types of particles in generations $0,\ldots,n-1$, the random variables
$(Y'_\ibf)_{\ibf\in\S'_n}$ are i.i.d.\ with law $\P[Y'_\ibf\leq s]=\ha s/t$
$(s\in[0,t])$. We claim that they satisfy the inductive relation
\be\label{Yind2}
Y'_\ibf=\chi[\om'_\ibf](Y'_{\ibf 1},Y'_{\ibf 2})\qquad(\ibf\in\S'),
\ee
where (compare (\ref{chidef}))
\be\label{chidef2}
\chi[\om](x,y):=\left\{\ba{ll}
x\quad&\mbox{if }\om\in[0,t),\ x>\om,\\[5pt]
\infty\quad&\mbox{if }\om\in[0,t),\ x\leq\om,\\[5pt]
x\wedge y\quad&\mbox{if }\om=2.\ea\right.
\ee
Note that $\ibf 2\not\in\S'$ if $\om_\ibf\in[0,t)$, but since in this case,
$\chi[\om_\ibf](x,y)$ does not depend on $y$, (\ref{Yind2}) is unambiguous.
Indeed, (\ref{Yind2}) follows from the fact that the original random variables
$(Y_\ibf)_{\ibf\in\T}$ satisfy the inductive relation (\ref{Yind}) and,
in view of (\ref{om12}), $Y_\ibf=Y_{\ibf 1}$ if $\ibf\in\A$ is of type
$\om_\ibf=1$.

These observations imply the statement of Proposition~\ref{P:jointscale}.
Indeed, if we set
\be
Y^\ast_\ibf:=t^{-1}Y'_\ibf,\quad
\om^\ast_\ibf:=\left\{\ba{ll}
t^{-1}\om'_\ibf\quad&\mbox{if }\om'_\ibf\in[0,t),\\
2\quad&\mbox{if }\om'_\ibf=2,\ea\right.\quad
\ell^\ast_\ibf:=t^{-1}\ell'_\ibf,
\ee
then the random variables $\S'$ and
$(\om^\ast_\ibf,\ell^\ast_\ibf)_{\ibf\in\S'}$ define a marked tree
$(\Ti^\ast,\Pi^\ast)$ such that the joint law of
$(Y^\ast_\wurz,\Ti^\ast,\Pi^\ast)$, conditioned on $\Ti^\ast\neq\emptyset$, is
equal to the joint law of $(Y_\wurz,\Ti,\Pi)$.
\epro

\bpro[Proof of Lemma~\ref{L:relinfo}]
We use notation as in the proof of Propositions~\ref{P:scale} and
\ref{P:jointscale}. We adapt the proof of \cite[Lemma~46]{MSS20} to our
present setting. We set $\T_{(n)}:=\{\ibf\in\T:|\ibf|<n\}$ and let
$\ov\Fi_{(n)}$ and $\ov\Fi$ be the \si-fields generated by the random
variables $\tau_\ibf,\kappa_\ibf$ with $\ibf\in\T_{(n)}$ and $\ibf\in\T$,
respectively. We also set $\S'_{(n)}:=\S'\cap\T_{(n)}$, we let $\Fi_{(n)}$ be
the \si-field generated by the random variables $\S'_{(n)}$ and
$(\om'_\ibf)_{\ibf\in\S'_{(n)}}$, and we define $\Fi$ similarly, with
$\S'_{(n)}$ replaced by $\S'$. We observe that $\Fi_{(n)}\sub\ov\Fi_{(n)}$
$(n\geq 1)$.

The inductive relation (\ref{Yind2}) shows that
conditional on $\Fi_{(n)}$, the state at the root $Y'_\wurz$ is a
deterministic function of $(Y'_\ibf)_{\ibf\in\T_n}$. Since
$(Y'_\ibf)_{\ibf\in\T_n}$ are independent of $\ov\Fi_{(n)}$, it follows that
$Y'_\wurz$ is conditionally independent of $\ov\Fi_{(n)}$ given $\Fi_{(n)}$,
i.e.,
\be
\P\big[Y'_\wurz\in A\,\big|\,\ov\Fi_{(n)}\big]
=\P\big[Y'_\wurz\in A\,\big|\,\Fi_{(n)}\big]\quad{\rm a.s.}
\ee
for any measurable $A\sub\R$. Letting $n\to\infty$, using martingale
convergence
and observing that $Y'_\wurz$ contains the same information as $Y^\ast_\wurz$
while $(\Ti^\ast,\Pi^\ast)$ contains the same information as $\Fi$, the claim
follows.
\epro

\begin{remark}
It follows from Lemma~\ref{L:scinv} that $\un\rho^{(2)}=\rho^{(2)}_2$, the
nontrivial scale-invariant fixed point from Theorem~\ref{T:scalefix}.
Therefore, combining Lemma \ref{L:nurho} with formula
(\ref{finite_both_coord}), we obtain a formula for $\un\nnu^{(2)}$. Indeed,
\be\label{explit}
\un\nnu^{(2)}\big([0,r]\times[0,s]\big)
=2-\frac{1}{2r}-\frac{1}{2s}
-\Big(2-\frac{1}{s\vee r}\Big)f_{c_2}
\left(\frac{2-\frac{1}{s\wedge r}}{2-\frac{1}{s\vee r}}\right)
\ee
$(\ha<r,s\leq 1)$,
where $f_{c_2}$ is the function defined in Theorem~\ref{T:scalefix}.
\end{remark}

\subsection{Frozen percolation on the 3-regular tree}\label{S:unorient}

In this subsection, we use methods from \cite{Ald00} to derive
Theorems~\ref{T:orfrz} and \ref{T:main}, which are concerned with the
unoriented 3-regular tree, from Theorems~\ref{T:dirfrz} and \ref{T:nonend},
which are concerned with the oriented binary tree. We start with a preparatory
lemma.

Let $(\Ui,\vec F)$ satisfy properties (i)--(iii) of Theorem~\ref{T:orfrz} and
let $F$ be defined in terms of $\vec F$ as in that theorem. Recall that we
call edges in $E_t\beh F$ open, edges in $E_t\cap F$ frozen, and all other
edges closed. A similar convention applies in the oriented setting. For each
$w\in T$ and $t\in[0,1]$, let $C_t(w)$ resp.\ $\vec C_t(w)$ denote the set of
vertices that can at time $t$ be reached by an open unoriented resp.\ oriented
path starting at $w$.

\bl[Finite unoriented clusters]
Almost\label{L:finclust} surely, for all $t\in[0,1]$, if $C_t(w)$ is finite,
then $C_t(w)=\vec C_t(w)$.
\el

\bpro
Clearly $C_t(w)\sub\vec C_t(w)$ regardless of whether $C_t(w)$ is finite or not.
To see that equality holds if $C_t(w)$ is finite, assume the converse.
Then there must be $x\in C_t(w)$ and $y\not\in C_t(w)$ such that the
oriented edge $(x,y)$ is open at time $t$. Among all such edges, we can choose
the unique one for which $s:=U_{\{x,y\}}$ is minimal. Since $y\not\in C_t(w)$,
the oriented edge must have frozen at time $s$, so by property~(i) of
Theorem~\ref{T:orfrz}, at time $s$ there must be
an open ray starting at $x$ not using $y$. Such a ray must
use an oriented edge to leave $C_t(w)$ that is open at time $s$ and hence also
at the later time $t$, contradicting the minimality of $U_{\{x,y\}}$.
\epro

\bpro[Proof of Theorem~\ref{T:orfrz}]
We first prove uniqueness. Assume that $\vec F$ satisfies properties
(i)--(iii). For each $(v,w)\in\vec E$, let
\bc\label{XFdef}
\dis X_{(v,w)}&:=&\dis\inf\big\{t\in[0,1]:\exists\,\mbox{ray $(v_n,w_n)_{n\geq 0}$
  starting with $(v_0,w_0)=(v,w)$}\\
&&\dis\phantom{\inf\big\{t\in[0,1]:\;}
\mbox{such that }(v_n,w_n)\in\vec F\ \forall n\geq 0\big\},
\ec
with $\inf\emptyset:=\infty$. Let $\ga$ be the map in
(\ref{gadef}). Property~(i) implies that
\be\label{3ind}
X_{(x,v)}:=\ga[U_{\{x,v\}}](X_{(v,y)},X_{(v,z)})
\ee
whenever $v\in T$ and $x,y,z$ are the three neighbours of $v$.
Let $S$ be a finite subtree of $(T,E)$. Then, for each $(v,w)\in\pa S$, the
set $\vec E_{(v,w)}$ is naturally isomorphic to the oriented binary tree $\T$.
Formula (\ref{3ind}) and properties~(ii) and (iii) imply that
$(\Ui_{\{x,y\}},X_{(x,y)})_{(x,y)\in\vec E_{(v,w)}}$ is an RTP corresponding
to the map $\ga$ and some solution $\mu$ to the RDE (\ref{frzRDE}).
Property~(i) and Theorem~\ref{T:dirfrz} imply that $\mu=\nu$, the measure
defined in (\ref{nudef}). By property~(iii), the RTPs corresponding to
different $(v,w)\in\pa S$ are independent. By (\ref{3ind}), these RTPs
uniquely determine $X_{(x,y)}$ for each $(x,y)\in\vec E$. This shows that the
joint law of $\Ui=(U_{\{x,y\}})_{\{x,y\}\in E}$ and $(X_{(x,y)})_{(x,y)\in E}$ is
uniquely determined. Since
\be\label{vecF}
(x,v)\in\vec F\mbox{ if and only if }U_{\{x,v\}}\geq X_{(v,y)}\wedge X_{(v,z)}
\ee
whenever $v\in T$ and $x,y,z$ are the three neighbours of $v$, the joint law
of $(\Ui,\vec F)$ is also uniquely determined.

As Aldous already showed in \cite{Ald00}, existence follows basically from the
same argument. We fix a finite subtree $S$ of $(T,E)$, construct independent
RTPs corresponding to $\ga$ and $\nu$ for each $(v,w)\in\pa S$, inductively
define $X_{(x,y)}$ for each $(x,y)\in\vec E$ by (\ref{3ind}), and then define
$\vec F$ by (\ref{vecF}). It follows from the properties of RTPs that if we
add a vertex to $S$ or remove a vertex, then the law of the object we have
just constructed does not change. As a result, our construction is independent
of the choice of $S$, the law of $(\Ui,\vec F)$ is invariant under
automorphisms of the tree, and property~(iii) holds for general
$S$. Property~(i) now follows from Theorem~\ref{T:dirfrz}, completing the
proof that an object satisfying (i)--(iii) exists.

It is clear that $(\Ui,F)$, defined in terms of $(\Ui,\vec F)$, is invariant
under automorphisms of the tree. To see that it also satisfies property~(i) of
Theorem~\ref{T:frz}, we observe that by the way $F$ has been defined in terms
of $\vec F$ and property~(i) of Theorem~\ref{T:orfrz}, $\{v,w\}\not\in F$ if
and only if for each $t<U_{\{v,w\}}$, the oriented clusters $\vec C_t(v)$ and
$\vec C_t(w)$ are both finite. By Lemma~\ref{L:finclust}, this is equivalent
to $C_t(v)$ and $C_t(w)$ being finite, proving property~(i) of
Theorem~\ref{T:frz}.
\epro

The following simple abstract lemma prepares for the proof of
Theorem~\ref{T:main}.

\bl[Almost surely not equal]
Let\label{L:asdif} $(\om_\ibf,X_\ibf)_{\ibf\in\T}$ be a nonendogenous RTP,
where $\T$ denotes the space of all finite words made up from the alphabet
$\{1,\ldots,d\}$, with $d\geq 2$. Let $(X'_\ibf)_{\ibf\in\T}$ be a copy of
$(X_\ibf)_{\ibf\in\T}$, conditionally independent given
$(\om_\ibf)_{\ibf\in\T}$. Then $(X_\ibf)_{\ibf\in\T}\neq(X'_\ibf)_{\ibf\in\T}$
a.s.
\el

\bpro
Let $\nu$ denote the solution of the RDE used to construct the RTP.
Let $\T_n:=\{\ibf\in\T:|\ibf|=n\}$. Then $(X_\ibf,X'_\ibf)_{\ibf\in\T_n}$ are
i.i.d.\ with common law $\un\nu^{(2)}$ as in (\ref{unnu}). By
Theorem~\ref{T:bivar}, $\un\nu^{(2)}\neq\ov\nu^{(2)}$, which implies that $p:=\P[
X_\ibf\neq X'_\ibf]>0$ and hence
\be
\P\big[(X_\ibf)_{\ibf\in\T}\neq(X'_\ibf)_{\ibf\in\T}\big]\leq
\P\big[X_\ibf=X'_\ibf\mbox{ for all }\ibf\in\T_n\big]\leq(1-p)^{d^n}.
\ee
Since $d\geq 2$ and $n$ is arbitrary, the claim follows.
\epro

\bpro[Proof of Theorem~\ref{T:main}]
We use the construction of $(\Ui,\vec F)$ in the proof of
Theorem~\ref{T:orfrz}. We fix a finite subtree $S$ of $(T,E)$. Independently
for each $(v,w)\in\pa S$, we construct an RTP
$(\Ui_{\{x,y\}},X_{(x,y)})_{(x,y)\in\vec E_{(v,w)}}$ corresponding to the map
$\ga$ in (\ref{gadef}) and measure $\nu$ in (\ref{nudef}), and we let
$(X'_{(x,y)})_{(x,y)\in\vec E_{(v,w)}}$ be a copy of
$(X_{(x,y)})_{(x,y)\in\vec E_{(v,w)}}$, conditionally independent given
the random variables $(\Ui_{\{x,y\}})_{\{x,y\}\in E_{(v,w)}}$. Using (\ref{3ind}), we inductively
define $X_{(x,y)}$ and $X'_{(x,y)}$ for all $(x,y)\in\vec E$ and in terms of
these random variables we define $\vec F$ and $\vec F'$ as in (\ref{vecF}),
which are finally used to define $F$ and $F'$ as in Theorem~\ref{T:orfrz}.
Then $F$ and $F'$ are conditionally independent given $\Ui$.

It follows from Theorem~\ref{T:nonend} and Lemma~\ref{L:asdif} that
a.s.\ $X_{(x,y)}\neq X'_{(x,y)}$ for some $(x,y)\in\vec F$. By (\ref{XFdef}),
this implies that $\vec F\neq\vec F'$ a.s. By Lemma~\ref{L:finclust} and
property~(i) of Theorem~\ref{T:orfrz}, the set $\vec F$ is a.s.\ determined by
the pair $(\Ui,F)$, and likewise $\vec F'$ is a.s.\ determined by $(\Ui,F')$,
so $\vec F\neq\vec F'$ a.s.\ implies $F\neq F'$ a.s.
\epro

\appendix

\section{Skeletal branching processes}\label{A:skel}

Informally speaking, the \emph{skeletal process} of a branching process is the
process consisting of those particles whose offspring will never die out. It
is well-known that the skeletal process of a branching process is itself a
branching process. For discrete time processes, a proof can be found in
\cite[Thm~I.12.1]{AN72}. There is also an extensive literature about skeletal
processes of superprocesses, see \cite{EKW15} and references therein.
In this appendix, we show how the skeletal process of a continuous-time
branching process can be calculated, and use this to sketch an alternative
proof that $\Ti'$, defined in (\ref{Tiac}), is the family tree of a binary
branching process with branching rate~$t$.

Generalising our set-up, let $(Z_h)_{h\geq 0}$ be a continuous-time branching
processes in which each particle is with rate $\rbf(k)$ replaced by $k$ new
particles. A sufficient condition for $(Z_h)_{h\geq 0}$ to be well-defined and
nonexplosive is that $\sum_k\rbf(k)k<\infty$. A convenient tool is the
\emph{generating semigroup} $(U_h)_{h\geq 0}$ defined as $U_h\phi:=u_h$
$(\phi\in[0,1])$, where $(u_h)_{h\geq 0}$ is the unique solution with initial
state $u_0=\phi$ to the differential equation
\be\ba{l}\label{Psidef}
\dis\dif{h}u_h=\Psi(u_h)\quad(h\geq 0)\\[5pt]
\dis\quad\mbox{with}\quad
\Psi(u):=\sum_{k\geq 0}\rbf(k)\big\{(1-u)-(1-u)^k\big\}\quad\big(u\in[0,1]\big).
\ec
The generating semigroup uniquely determines the transition probabilities of
$(Z_h)_{h\geq 0}$ through the relation
\be\label{bradual}
\E\big[(1-\phi)^{Z_h}\big]=\E\big[(1-U_h\phi)^{Z_0}\big]
\qquad(h\geq 0).
\ee
This can be deduced, for example, from \cite[Sect.~III.3]{AN72}, although the
notation there is quite different.

Let $p$ be the survival probability of $(Z_h)_{h\geq 0}$, which is the largest
root in $[0,1]$ of the equation $\Psi(p)=0$. Then we claim that setting
\be\label{skeleton}
U'_h\phi:=p^{-1}U_h(p\phi)\qquad\big(\phi\in[0,1]\big)
\ee
defines a generating semigroup, which corresponds to the skeletal process
$(Z'_h)_{h\geq 0}$ of $(Z_h)_{h\geq 0}$. For discrete time processes, a proof
can be found in \cite[Thm~I.12.1]{AN72}. The statement for continuous-time
processes can easily be derived from this by adding independent exponentially
distributed lifetimes to the discrete time process. In particular, if
$r(0)=1-t$, $r(2)=1$, and all other rates are zero, then the differential
equation in (\ref{Psidef}) reads
\be\label{difbin}
\dif{h}u_h=\Psi(u_h)=u_h(1-u_h)-(1-t)u_h\qquad(h\geq 0),
\ee
and $(U'_h)_{h\geq 0}$ is given by the solutions to the differential equation
\be\label{difskel}
\dif{h}v_h=t^{-1}\Psi(tv_h)=t^{-1}\big(tv_h(1-tv_h)-(1-t)tv_h\big)
=tv_h(1-v_h)\qquad(h\geq 0),
\ee
which we recognise as the generating semigroup of a branching process where
particles split into two with rate $t$ and never die.

The transformation in (\ref{skeleton}) can be traced back to \cite{Har48}
while the interpretation in terms of the skeletal process dates back to
\cite[Thm~I.12.1]{AN72}. See also \cite[Thm~9]{FS04} for a statement in the
context of superprocesses. It is possible to go further and write $(Z_h)_{h\geq
  0}$ as the union of skeletal and non-skeletal particles, which then form a
two-type branching process. This sort of statements date back to \cite{OCo93}
and have been developed and exploited in a superprocess setting; see
\cite{EKW15} and references therein.

\end{document}